\documentclass[11pt,namelimits,sumlimits,a4paper]{amsart}
\usepackage{cite}
\usepackage{comment}

\usepackage{amssymb,amsmath}
\usepackage[mathscr]{eucal}
\usepackage{slashed}

\usepackage{graphicx}
\usepackage{amsmath,amsthm}
\usepackage{graphics}
\usepackage{color}
\usepackage{epsfig}
\usepackage{amssymb,amsmath}
\usepackage[mathscr]{eucal}

\usepackage[latin1]{inputenc}
\usepackage[T1]{fontenc}
\usepackage[UKenglish]{babel}
\usepackage{amsfonts}
\usepackage{fancyhdr}
\usepackage{graphicx}
\usepackage{amsmath}
\usepackage{amsthm}
\usepackage{amsmath,amscd}
\usepackage{latexsym}
\usepackage{cite}
\usepackage{amssymb,amsmath}

\usepackage{comment}
\usepackage[mathscr]{eucal}
\usepackage{slashed}
\usepackage{times,psfrag}
\usepackage{mathtools}
\usepackage{hyperref}
\usepackage{multirow}
\usepackage{longtable}
\usepackage{bm}
\usepackage[all]{xypic}
\usepackage{fancyhdr}
\usepackage{float}

\usepackage{esint}

\usepackage[]{graphicx}
\usepackage{graphics}
\usepackage{color}
\usepackage{epsfig}

\usepackage{layout}
\usepackage{enumitem}

\usepackage{lmodern}
\usepackage[rgb]{xcolor}
\usepackage[draft,author={Lorenzo Foscolo}]{pdfcomment}

\numberwithin{equation}{section}
\newtheorem{theorem}[equation]{Theorem}
\newtheorem{lemma}[equation]{Lemma}
\newtheorem{prop}[equation]{Proposition}
\newtheorem{corollary}[equation]{Corollary}

\theoremstyle{definition}
\newtheorem{definition}[equation]{Definition}

\theoremstyle{remark}
\newtheorem{remark}[equation]{Remark}
\newtheorem*{remark*}{Remark}

\newtheorem*{assumption}{Assumption}

\newcounter{mtheorem}

\newtheoremstyle{mystyle}
  {}
  {}
  {\itshape}
  {}
  {\bfseries}
  {.}
  { }
  {}

\theoremstyle{mystyle}
\newtheorem{mtheorem}[mtheorem]{Theorem}

\newcommand{\ie}{\emph{i.e.} }

\newcommand{\beq}{\begin{equation}}
\newcommand{\eeq}{\end{equation}}
\newcommand{\bea}{\begin{eqnarray}}
\newcommand{\eea}{\end{eqnarray}}

\newcommand{\C}{\mathbb{C}}

\newcommand{\R}{\mathbb{R}}

\newcommand{\Z}{\mathbb{Z}}
\newcommand{\N}{\mathbb{N}}

\newcommand{\HH}{\mathbb{H}}

\newcommand{\PP}{\mathbb{P}}
\newcommand{\Sph}{\mathbb{S}}
\newcommand{\ra}{\rightarrow}

\newcommand{\dvol}{\operatorname{dv}}
\newcommand{\Real}{\operatorname{Re}}
\newcommand{\Imag}{\operatorname{Im}}

\newcommand{\Lie}[1]{\mathfrak{#1}}

\newcommand{\BC}{\textup{BC}}
\newcommand{\tu}[1]{\textup{#1}}

\newcommand{\gtwo}{\ensuremath{\textup{G}_2}}

\newcommand{\unitary}[1]{\textup{U$(#1)$}}

\newcommand{\sunitary}[1]{\textup{SU$(#1)$}}

\newcommand{\sorth}[1]{\textup{SO$(#1)$}}
\newcommand{\Sp}[1]{\textup{Sp$(#1)$}}

\newcommand{\spins}{\ensuremath{\textup{Spin}(7)}}

\newcommand{\triple}[1]{\boldsymbol{#1}}
\newcommand{{\isomgtc}}{\ensuremath{\sunitary{2}^3 \rtimes S_3}}

\def\co{\colon\thinspace}

\begin{document}

\title[Complete non-compact Spin$(7)$--manifolds from self-dual Einstein $4$-orbifolds]{Complete non-compact Spin$(7)$--manifolds from self-dual Einstein $4$-orbifolds}
\author{Lorenzo Foscolo}
\maketitle

\begin{abstract}
We present an analytic construction of complete non-compact $8$-dimensional Ricci-flat manifolds with holonomy $\spins$. The construction relies on the study of the adiabatic limit of metrics with holonomy $\spins$ on principal Seifert circle bundles over asymptotically conical $\gtwo$--orbifolds. The metrics we produce have an asymptotic geometry, so-called ALC geometry, that generalises to higher dimensions the geometry of 4-dimensional ALF hyperk\"ahler metrics.

We apply our construction to asymptotically conical $\gtwo$--metrics arising from self-dual Einstein $4$-orbifolds with positive scalar curvature. As illustrative examples of the power of our construction, we produce complete non-compact $\spins$--manifolds with arbitrarily large second Betti number and infinitely many distinct families of ALC $\spins$--metrics on the same smooth $8$-manifold.
\end{abstract}

\section{Introduction}

In this paper we provide a new analytic construction of complete non-compact Ricci-flat $8$-manifolds with holonomy $\spins$ and non-maximal volume growth. The starting point of the construction is an asymptotically conical (AC) $7$-dimensional orbifold $(B,g_0)$ with holonomy $\gtwo$ together with a suitable principal circle orbibundle $\pi\co M\ra B$ with total space $M$ a smooth $8$-manifold (we will then say that $\pi\co M\ra B$ is a Seifert bundle). Our method then produces a $1$-parameter family $\{ g_\epsilon\}_{\epsilon>0}$ of circle invariant $\spins$--metrics on $M$ such that $(M,g_\epsilon)$ collapses back to the orbifold $(B,g_0)$ as $\epsilon\ra 0$. The metric $g_\epsilon$ has controlled asymptotic geometry, so-called ALC (asymptotically locally conical) geometry: along the (unique) end of $M$ the metric $g_\epsilon$ approaches a Riemannian submersion with base a conical (orbifold) metric and circle fibres of fixed finite length.

\begin{mtheorem}\label{thm:Main}
Let $M^8$ be a smooth non-compact $8$-manifold with an almost-free circle action, \ie such that the quotient space $B=M/S^1$ is an orbifold. Assume that $B$ carries an AC orbifold metric $g_0$ with holonomy $\gtwo$ and that the principal circle orbibundle $M\ra B$ satisfies the topological condition
\[
c_1^{orb}(M)\cup [\varphi_0] = 0 \in H^5_{orb}(B).
\]
Here $\varphi_0$ is the closed and coclosed $\gtwo$ $3$-form on $B$ inducing the $\gtwo$--metric $g_0$.

Then for every $\epsilon>0$ sufficiently small there exists a circle-invariant ALC $\spins$--metric $g_\epsilon$ on $M$ such that the sequence $(M,g_\epsilon)$ collapses to $(B,g_0)$ with bounded curvature as $\epsilon\ra 0$.
\end{mtheorem}
We refer the reader to Theorem \ref{thm:Main:Technical} later in the paper for a more precise statement. 

\subsubsection*{Motivation and applications} 

In \cite{FHN:ALC:G2:from:AC:CY3}, in collaboration with Haskins and Nordstr\"om, we developed a similar construction of highly collapsed ALC $\gtwo$--holonomy metrics on suitable principal circle bundles over \emph{smooth} AC Calabi--Yau $3$-folds. The construction of \cite{FHN:ALC:G2:from:AC:CY3} allowed us to exploit recent progress on the existence of Calabi--Yau cone metrics \cite{GMSW,Futaki:Ono:Wang,Collins:Szekelyhidi} and AC Calabi--Yau metrics  \cite{VanCoevering:Existence,Goto:Crepant,Conlon:Hein:I} to produce infinitely many complete non-compact $\gtwo$--manifolds and complete $\gtwo$--metrics depending on an arbitrarily large number of parameters. Only a handful of complete non-compact $\gtwo$--manifolds was previously known.

The existence of an analogous construction of $\spins$--metrics from AC $\gtwo$--manifolds as in Theorem \ref{thm:Main} is therefore not in itself surprising. The fact that such a construction can be used to produce significant results in $\spins$--geometry, however, is a priori much less clear: the naive generalisation of \cite{FHN:ALC:G2:from:AC:CY3} to the $\spins$--setting using only smooth manifolds would be a theorem that currently applies to only one example! The simultaneous extension of \cite{FHN:ALC:G2:from:AC:CY3} to the orbifold setting is the crucial new ingredient that makes Theorem \ref{thm:Main} useful. Indeed, in contrast to the Calabi--Yau case, our current knowledge of smooth AC $\gtwo$--manifolds is extremely limited: in 1989 Bryant--Salamon \cite{Bryant:Salamon} constructed three (explicit) examples of AC $\gtwo$--metrics; only very recently \cite[Theorem C]{ALC:G2:coho1} an infinite family of new simply connected examples has been found. On the other hand, Bryant--Salamon's construction of AC $\gtwo$--metrics yields an AC \emph{orbifold} $\gtwo$--metric on the total space of the orbibundle of anti-self-dual $2$-forms over any self-dual Einstein $4$-orbifold with positive scalar curvature. This construction yields a large supply of AC $\gtwo$--orbifolds since many self-dual Einstein orbifold metrics can be constructed using the quaternionic K\"ahler quotient construction of Galicki--Lawson \cite{Galicki:Lawson}. By using such orbifolds we are able to produce a wealth of new complete non-compact $\spins$--manifolds.

\begin{mtheorem}\label{thm:Examples:Betti}
For every $k\geq 1$ there exists a smooth non-compact $8$-manifold that retracts onto $\sharp_k (S^2\times S^3)$ and carries a family of complete ALC $\spins$--metrics. In particular, there exist complete non-compact $\spins$--manifolds with arbitrarily large second Betti number.
\end{mtheorem}
Only a handful of complete non-compact $\spins$--metrics was previously known \cite{Bryant:Salamon,CGLP:A8:B8,CGLP:B8,Bazaikin:C8,Bazaikin:B8,Gukov:Sparks,Kovalev:ACyl:Spin7}. As a further illustration of the power of our construction, we also find a smooth non-compact $8$-manifold that can be described as a circle orbibundle over an AC $\gtwo$--orbifold in infinitely many different ways.

\begin{mtheorem}\label{thm:Examples:Families}
The non-trivial rank-$3$ real vector bundle over $S^5$ carries infinitely many families of complete ALC $\spins$--metrics. Different families are distinguished by their (unique) tangent cone at infinity.
\end{mtheorem}
In other words there are infinitely many inequivalent circle actions on the $8$-manifold $M$ in question such that the orbit space $M/S^1$ is the orbibundle of anti-self-dual $2$-forms over a self-dual Einstein $4$-orbifold with positive scalar curvature.

The analytic framework introduced in this paper to work on AC orbifolds can also be exploited to extend the construction of complete ALC $\gtwo$--metrics in \cite{FHN:ALC:G2:from:AC:CY3} to the orbifold setting. In \cite{FHN:ALC:G2:from:AC:CY3} examples of ALC $\gtwo$--metrics arose from AC Calabi--Yau metrics on crepant resolutions of Calabi--Yau cones. Often it is natural to consider only \emph{partial} resolutions of Calabi--Yau cones, which replace the singularity of the cone with simpler (albeit non-isolated) orbifold singularities. Combining the techniques of this paper with \cite{FHN:ALC:G2:from:AC:CY3} allows us to construct complete $\gtwo$--metrics on suitable circle orbibundles over these orbifold partial resolutions. As an illustration of the possible complete $\gtwo$--metrics arising from this construction, Theorem \ref{thm:ALC:G2} establishes an analogue of Theorem \ref{thm:Examples:Families} in the $\gtwo$ setting by constructing infinitely many distinct families of ALC $\gtwo$--metrics on $S^3\times\R^4$.

Considering sequences of $\spins$--metrics collapsing to $\gtwo$--orbifolds and not only smooth manifolds is also very natural from the point of view of the theory of Riemannian collapse. Indeed, Fukaya \cite[Proposition 11.5]{Fukaya} has shown that Gromov--Hausdorff limits of Riemannian manifolds that collapse with bounded curvature in codimension $1$ must be orbifolds. Thus, besides extending it to the $\spins$ setting, Theorem \ref{thm:Main} extends the construction of \cite{FHN:ALC:G2:from:AC:CY3} to its most general context.

There are three main aspects to the proof of Theorem \ref{thm:Main} and of its corollaries Theorems \ref{thm:Examples:Betti} and \ref{thm:Examples:Families}. The general strategy of the proof of Theorem \ref{thm:Main} relies on the adiabatic limit of $\spins$--metrics with a circle symmetry. The strategy is partially motivated by known families of cohomogeneity one ALC $\spins$--metrics and the duality between M theory and Type IIA String theory in theoretical physics. Successfully implementing this strategy requires a refined knowledge of closed and coclosed forms on AC manifolds and orbifolds. Note that the orbifolds we consider in this paper are non-compact and have a singular set that is allowed to extend all the way to infinity. Describing the analytic framework to work on such orbifolds is an important technical aspect of this paper. Finally, the third aspect of this work is the search for self-dual Einstein $4$-orbifolds with positive scalar curvature that give rise to concrete examples of AC $\gtwo$--orbifolds to feed into Theorem \ref{thm:Main}. In the rest of this Introduction we discuss each of these three aspects.

\subsubsection*{The adiabatic limit of circle invariant $\spins$--metrics}

In \cite{Bryant:Salamon} Bryant--Salamon constructed the first known example of a complete metric with holonomy $\spins$. The Bryant--Salamon metric is an explicit AC $\spins$--metric on the spinor bundle $\slashed{S}(\Sph^4)$ of the round $4$-sphere. The metric is asymptotic at infinity to the Riemannian cone over the squashed Einstein metric on $S^7$ \cite{Jensen}.

The Bryant--Salamon metric on $\slashed{S}(\Sph^4)$ is in fact invariant under the natural cohomogeneity one action of $\tu{Sp}(2)$ (recall that $\tu{Sp}(2)\simeq \tu{Spin}(5)$). In \cite{CGLP:A8:B8} Cveti\v{c}--Gibbons--L\"u--Pope studied the ODE system describing general $\tu{Sp}(2)$--invariant $\spins$--metrics. They found a new explicit example and argued that it moved in a $1$-parameter family up to scale, the existence of which was further studied in \cite{CGLP:B8} and rigorously proved in \cite{Bazaikin:B8}. The asymptotic behaviour of the metrics in this family, labelled $\mathbb{B}_8$ in the physics literature, is different from the Bryant--Salamon metric: the $\mathbb{B}_8$ metrics have non-maximal volume growth $r^7$; at infinity the metric approaches a Riemannian submersion with base metric the Riemannian cone over the nearly K\"ahler metric on $\C\PP^3$ and circle fibres with fixed finite length $\ell$. In \cite{CGLP:A8:B8} the acronym ALC (asymptotically locally conical) was introduced to describe this asymptotic geometry: the ALC asymptotic geometry is analogous to the asymptotic geometry of ALF (asymptotically locally flat) $4$-dimensional manifolds, except that the tangent cone at infinity is not necessarily flat. Up to scale, the asymptotic length $\ell$ of the circle fibres can be taken as the parameter that distinguishes different members of the $\mathbb{B}_8$ family. In fact, the asymptotic circle action extends to a global symmetry of the $\mathbb{B}_8$ metrics, which are not only invariant under the left action of $\tu{Sp}(2)$, but also under the circle acting on the fibres of $\slashed{S}(\Sph^4)$ as the Hopf circle action. This circle action is not free since it fixes the zero-section, but the quotient space is still a smooth manifold, $\Lambda^-T^\ast\Sph^4$. As $\ell\ra 0$ the family of ALC $\mathbb{B}_8$ metrics collapses to the Bryant--Salamon AC $\gtwo$--metric on $\Lambda^-T^\ast\Sph^4$ \cite{Bryant:Salamon}, which is asymptotic to the cone over the nearly K\"ahler metric on $\C\PP^3$. As $\ell\ra 0$ the curvature of the $\mathbb{B}_7$ metrics blows up along the zero-section $\Sph^4$, the fixed locus of the circle action on $\slashed{S}(\Sph^4)$.

These first examples lead to an explosion of activity in the physics and, later, mathematics literature discussing further (conjectural) families of ALC manifolds with exceptional holonomy. An explicit ALC $\spins$--metric on $\R^8$ was found in \cite{CGLP:A8:B8} and a new family of $\tu{Sp}(2)$--invariant ALC $\spins$--metrics on the canonical line bundle of $\C\PP^3$ was studied numerically in \cite[Section 2]{CGLP:Coho1} and later constructed rigorously in \cite{Bazaikin:C8}. Further work concentrated on the case of cohomogeneity one $\sunitary{3}$--invariant $\spins$--metric with principal orbits the Aloff--Wallach spaces $\sunitary{3}/U(1)_{k,l}$, where the integers $k,l$ determine the embedding of $\unitary{1}$ in the maximal torus of $\sunitary{3}$: the discovery of some explicit solutions, numerical investigations of the relevant ODE systems and a rigorous study of local solutions defined in a neighbourhood of the possible singular orbits were carried out by various authors \cite{CGLP:Coho1,Gukov:Sparks,Gukov:Sparks:Tong,Kanno:Yasui:I,Kanno:Yasui:II,Reidegeld}. In general, however, the existence of complete solutions remains open.

From the physics perspective, the interest in ALC metrics with exceptional holonomy arises from the equivalence between M theory and Type IIA String theory in the limit of weak string coupling constant. Kaluza--Klein reduction of supersymmetric M-theory solutions along a circle of small radius proportional to the string coupling constant corresponds geometrically to the study of sequences of manifolds with exceptional holonomy collapsing in codimension $1$ along a (degenerate) circle fibration. For instance, the collapse of the $\mathbb{B}_8$ family of ALC $\spins$--metrics to the Bryant--Salamon AC $\gtwo$--metric in the limit $\ell\ra 0$ realises the duality between M theory ``compactified'' on $\slashed{S}(\Sph^4)$ and Type IIA String theory on $\Lambda^-T^\ast\Sph^4$ with a D6-brane wrapping the zero-section. The geometric interpretation of the latter physical jargon is that $\slashed{S}(\Sph^4)\setminus\Sph^4$ can be regarded as a principal circle bundle over $\Lambda^-T^\ast\Sph^4\setminus\Sph^4$ with first Chern class evaluating to one on the $2$-sphere linking the zero-section. The collapse with bounded curvature exhibited by the families of ALC $\spins$--metrics in Theorem \ref{thm:Main} corresponds instead to the physical statement that the weak coupling limit of Type IIA theory on the $\gtwo$--orbifold $B$ with Ramond--Ramond $2$-form flux representing $c_1^{orb}(M)$ is equivalent to the low energy limit of M-theory on the total space $M$ of the circle orbibundle.

The idea of proof of Theorem \ref{thm:Main} is close to this physical interpretation, see for example \cite{CGLP:KK:reduction,Kaste:et:al}. We consider a $\spins$--manifold $(M,g)$ with an isometric circle action. Denote by $B$ the orbit space $M/S^1$ and assume for simplicity it is a smooth manifold. A $\spins$--metric is uniquely determined by a closed $4$-form $\Phi$ on $M$ satisfying certain point-wise nonlinear algebraic constraints. In the presence of a Killing field $\xi$ (that we assume preserves also $\Phi$) we can write the metric $g$ on $M$ as $g= h^{\frac{1}{3}}g_B + h^{-1}\theta^2$, where $\theta$ is an $S^1$--invariant $1$-form on $M$ dual to $\xi$, \ie a connection $1$-form on the principal circle bundle $M\ra B$, and $h$ and $g_B$ are a positive function and a Riemannian metric on $B$. We can then formulate the holonomy reduction of $g$ as a system of nonlinear partial differential equations $\Psi (\theta, h, \varphi)=0$. Here $\varphi = \xi \lrcorner\Phi$ defines a $\gtwo$--structure on $B$ inducing the metric $g_B$. In dimension $4$, the analogous dimensional reduction of hyperk\"ahler metrics along the orbits of a triholomorphic vector field yields the famous Gibbons--Hwaking Ansatz \cite{Gibbons:Hawking}, which reduces the existence of the hyperk\"ahler metric to a \emph{linear} equation on the orbit space $\R^3$. The system of equations $\Psi (\theta, h, \varphi)=0$ arising from the dimensional reduction of $\spins$--metrics is nonlinear and in general it is not at all clear how to study existence of solutions.

We then employ the strategy of deforming $\Psi (\theta, h, \varphi)=0$ to a different equation we can handle better. The most natural geometric degeneration is to consider families of $S^1$--invariant $\spins$--metrics with circle orbits of smaller and smaller length. We introduce a small parameter $\epsilon>0$ and consider a sequence of $S^1$--invariant metrics $g_\epsilon = h^{\frac{1}{3}}g_{B} + \epsilon^2 h^{-1}\theta$. The metric $g_\epsilon$ has $\spins$--holonomy if and only if $(\theta,h,\varphi)$ satisfy $\Psi (\epsilon\,\theta, h, \varphi)=0$. Here $\varphi=\xi\lrcorner\Phi_\epsilon$, where $\Phi_\epsilon$ is the $4$-form inducing $g_\epsilon$. For $\epsilon>0$ the equation $\Psi (\epsilon\,\theta, h, \varphi)=0$ is equivalent to $\Psi (\theta,h,\varphi)=0$ by scaling, but at $\epsilon=0$ the equation simplifies: solutions are of the form $(0,1,\varphi_0)$ where $\varphi_0$ is a torsion-free $\gtwo$--structure on the orbit space $B$, \ie the limiting metric $g_B$ induced by $\varphi_0$ has holonomy $\gtwo$. It is important to note that the equation $\Psi (\epsilon\,\theta, h, \varphi)=0$ depends smoothly on $\epsilon$ up to and including $\epsilon=0$. In Theorem \ref{thm:Main} we assume that an AC $\gtwo$--metric on $B$ is given and then try to perturb the solution $(0,1,\varphi_0)$ into a solution of $\Psi (\epsilon\,\theta, h, \varphi)=0$ for $\epsilon>0$.

The first step is to understand elements in the kernel of the linearisation $\mathcal{L}$ of $\Psi$ at $(0,1,\varphi_0)$, since they correspond to formal tangent vectors to curves of solutions to $\Psi (\epsilon\,\theta, h, \varphi)=0$ for $\epsilon \in [0,\epsilon_0)$. A dichotomy arises at this stage: it is geometrically meaningful to consider bounded solutions to the linearised problem as well as unbounded solutions with prescribed singularities in codimension $4$. In this paper we only consider the former case, which corresponds to sequences of $\spins$--metrics collapsing with bounded curvature; the case of codimension-$4$ singularities, related to collapse with unbounded curvature along the fibres of a circle fibration which degenerates in codimension $4$, is more involved and will be treated elsewhere. It turns out that bounded solutions $(\theta_0,h_0,\rho_0)$ to the linearised problem $\mathcal{L}(\theta_0,h_0,\rho_0)=0$ are completely determined by the choice of a principal circle bundle $M\ra B$ with $c_1 (M)=[d\theta_0]$. The topological constraint $c_1(M)\cup [\varphi_0]=0\in H^5(B)$ arises at this stage as the necessary and sufficient condition for solving the linearised problem.

We can now imagine reconstructing a curve of solutions to $\Psi (\epsilon\,\theta,h,\varphi)=0$ for $\epsilon\geq 0$ sufficiently small by deforming away from the initial solution $(0,1,\varphi_0)$ in the direction of $(\theta_0,h_0,\rho_0)$ via an application of the Implicit Function Theorem. The key step is the study of the mapping properties of the linear operator $\mathcal{L}$. Now, $\mathcal{L}$ is not obviously elliptic as it involves a combination of differential, codifferential and decomposition of differential forms into different types induced by the representation theory of $\gtwo$ (analogous to the $(p,q)$-type decomposition on complex manifolds). It is therefore not immediately obvious how to identify the cokernel of $\mathcal{L}$. In the construction of ALC $\gtwo$--metrics from AC Calabi--Yau $3$-folds in \cite{FHN:ALC:G2:from:AC:CY3} the linearised problem was complicated enough that we were only able to prove existence of solutions by solving the analogue of the equation $\Psi (\epsilon\,\theta,h,\varphi)=0$ as a power series in $\epsilon$, exploiting special cancellations that were only evident by solving the equation order-by-order in $\epsilon$. In this paper we are instead able to set up a direct argument using the Implicit Function Theorem. The key difference with respect to \cite{FHN:ALC:G2:from:AC:CY3} is that the space of $\gtwo$--structures on $\R^7$ is an open set in a linear space, while the space of $\sunitary{3}$--structures on $\R^6$ is cut out by nonlinear constraints and thereofore a further choice of ``exponential map'' is necessary. In order to understand the mapping properties of the linearised operator $\mathcal{L}$ and therefore prove Theorem \ref{thm:Main}, we need to exploit the interplay between the Laplacian, the Dirac operator and type decomposition of differential forms on AC $\gtwo$--orbifolds and, crucially, the fact that we restrict to variations of the $\gtwo$--structure $\varphi$ in the same cohomology class as $\varphi_0$.     

\vspace{-0.05cm}
\subsubsection*{Analysis on AC orbifolds}

The discussion so far has been formal. In order to implement the strategy we have just outlined we need to develop analytic tools to work on AC orbifolds. There are two main issues to take into account: the fact that we work on non-compact spaces and the orbifold singularities. The analysis of linear elliptic operators on \emph{smooth} AC manifolds using weighted Sobolev and H\"older spaces is well established. Analysis on \emph{compact} orbifolds has also been used in many geometric applications. However, the orbifolds we consider in this paper are non-compact and in view of the applications we have in mind we cannot insist that the singular set be compact. The simultaneous presence of an AC end and orbifold singularities allowed to extend to infinity seems not to have been considered before in the literature. We therefore felt it was necessary to include a self-contained exposition of the geometric and analytic tools we need. Since the orbifolds we consider arise as global quotients of smooth non-compact manifolds by a circle action, we develop the theory in a way that makes crucial use of this assumption. Instead of working on the AC orbifold $B$ itself we work on the smooth total space $M$ of the circle orbibundle over $B$: elliptic operators acting on differential forms on $B$ are replaced with transversally elliptic operators acting on basic forms on $M$. Note that since every Riemannian orbifold arises as the quotient of a smooth manifold by the action of a compact Lie group (the orthogonal frame orbibundle of an orbifold is always a smooth manifold, see \cite[Corollary 1.24]{Ruan:Orbifolds:Book}), a similar strategy can be (and has been) applied more generally. The case of \emph{Seifert circle bundles} (\ie principal circle orbibundles with smooth total space) allows us to give a particularly clean exposition.     

The central object in our exposition is the so-called \emph{adapted connection} $\nabla$ of a Riemannian foliation, see \cite[Definition 1.7]{Bismut} and \cite[Definition 3.13]{Reinhart}: a certain metric connection with torsion on $TM$ that preserves the splitting of the tangent bundle of $M$ into vertical and horizontal sub-bundles. We use $\nabla$ instead of the Levi--Civita connection of $M$ to define natural elliptic operators acting on basic sections of appropriate vector bundles. For instance, the exterior differential and codifferential acting on differential forms are replaced by the \emph{covariant} differential and codifferential induced by $\nabla$. Restricting these ``adapted'' operators to basic forms allows us to develop the linear theory of elliptic operators acting on weighted Banach spaces on AC orbifolds exactly as in the case of smooth AC manifolds. Once the right language has been developed, the only new analytic and geometric ingredient is Parker's Equivariant Sobolev Inequality \cite{Parker}.     

As non-experts in the theory of Riemannian foliations, we are unable to evaluate the originality in our treatment and how much our clean exposition depends on the restriction to the simple case of foliations with totally geodesic $1$-dimensional leafs. For example other authors use different ``adapted'' connections for different purposes instead of our uniform approach using $\nabla$, see \cite{Tondeur:Book}. An original contribution of this paper is a calculation of all the topological contributions to the weighted $L^2$--cohomology of AC manifolds and orbifolds. The $L^2$--cohomology of smooth AC manifolds is well known, see \cite[Theorem 1.A]{HHM} and \cite[Example 0.15]{Lockhart}. For geometric applications, however, it is often important to work with differential forms that are not necessarily square-integrable. For instance, there are many examples of higher-dimensional AC manifolds with special holonomy that are asymptotic to their tangent cone at infinity with a non-$L^2$ rate of decay. In Theorem \ref{thm:L2:cohomology} we apply the Fredholm theory we develop in Section \ref{sec:AC:Orbifolds} to give a complete description of the topological contributions to the weighted $L^2$--cohomology of AC manifolds. Our elementary proof immediately generalises to the case of AC orbifolds. Special cases of our result have been derived by other authors (see for example \cite[Section 4.5]{Karigiannis:Lotay}), but as far as we are aware a proof in arbitrary dimension is not currently available in the literature.

\subsubsection*{Self-dual Einstein $4$-orbifolds and special holonomy}

The analytic tools we develop in Section \ref{sec:AC:Orbifolds}, including our results about weighted $L^2$--cohomology of AC orbifolds, allow us to implement the adiabatic limit strategy and prove our main abstract existence result Theorem \ref{thm:Main}. The final part of the paper is devoted to the  study of concrete examples produced by this construction. All the examples we consider in the paper are obtained by applying Theorem \ref{thm:Main} to Bryant--Salamon's AC $\gtwo$--metrics arising from suitable self-dual Einstein $4$-orbifolds with positive scalar curvature.

It is well known that self-dual Einstein metrics with positive scalar curvature in dimension $4$ generate many different geometries related to special holonomy, see \cite[\S 13.4]{Boyer:Galicki}. If $Q$ is a self-dual Einstein $4$-manifold (or orbifold) with positive scalar curvature then its twistor space $Z$, the unit sphere bundle in the (orbi)bundle of anti-self-dual $2$-forms, carries two Einstein metrics with positive scalar curvature: a K\"ahler--Einstein metric \cite{Salamon:QK} and a nearly K\"ahler metric \cite{Eells:Salamon:Twistor}. The Konishi bundle $S$, the principal $\sunitary{2}$ or $\sorth{3}$ bundle associated with $Z$, also has two Einstein metrics: a $3$-Sasaki metric \cite{Konishi,Swann:Bundle} and a (strict) nearly parallel $\gtwo$--metric \cite{Galicki:Salamon,FKMS:nearly:parallel:G2}. Except for the K\"ahler--Einstein metric, these higher dimensional compact Einstein spaces carry real Killing spinors and are therefore related to special holonomy via the cone construction \cite{Bar}: the cone over the nearly K\"ahler metric on $Z$ has holonomy $\gtwo$, the cone over the $3$-Sasaki metric on $S$ is hyperk\"ahler and the one over the nearly parallel $\gtwo$--metric on $S$ has holonomy $\spins$. Furthermore, vector-bundle constructions of Ricci-flat metrics on (orbi)bundles over $Q$ can be used to produce non-compact spaces with special holonomy (partially) desingularising these cones. For example, a  well known seminal construction by Calabi \cite{Calabi:AC:CY} yields an AC Calabi--Yau metric on the canonical line bundle over the K\"ahler--Einstein $3$-fold $Z$; this metric is asymptotic to the cone over a finite quotient of $S$ by a cyclic group that only preserves one Sasaki structure in the $3$-Sasaki structure. Bryant--Salamon's construction \cite{Bryant:Salamon} of a (unique up to scale) AC $\gtwo$--metric on the (orbi)bundle of anti-self-dual $2$-forms on $Q$ plays a distinguished role in this paper.

There are also known constructions of $\spins$--holonomy metrics from self-dual Einstein $4$-manifolds: in \cite{Bryant:Salamon} Bryant--Salamon construct an AC $\spins$--metric on the spinor bundle of a spin self-dual Einstein $4$-manifold with positive scalar curvature; in \cite{Bazaikin:B8} Baza\u{\i}kin shows that the Bryant--Salamon AC metric is in fact the limit of a $1$-parameter family of ALC $\spins$--metrics on the same $8$-manifold; Baza\u{\i}kin \cite{Bazaikin:C8} also constructs families of ALC $\spins$--metrics arising as deformations of Calabi's AC Calabi--Yau metrics on $K_Z$. By a result of Hitchin \cite{Hitchin:SDE+} the only smooth self-dual Einstein $4$-manifolds with positive scalar curvature are $\Sph^4$ and $\C\PP^2$. As a consequence, there are only three smooth $\spins$--manifolds produced by these constructions. Bryant--Salamon's and Baza\u{\i}kin's constructions immediately generalise to self-dual Einstein $4$-orbifolds $Q$ to produce many \emph{singular} $\spins$--metrics. The reasons these metrics are never complete is that the self-dual Einstein $4$-orbifold $Q$ or its twistor space $Z$ (which is always singular if $Q$ is) are always embedded in the resulting spaces. If one considers principal orbibundles instead of vector bundles, however, it is instead often possible to obtain smooth manifolds. For example, Boyer--Galicki and their collaborators constructed infinitely many smooth $3$-Sasaki manifolds using self-dual Einstein $4$-orbifolds \cite{Boyer:Galicki:Inventiones}. Similarly, many smooth Sasaki--Einstein manifolds arise as circle orbibundles over K\"ahler--Einstein Fano orbifolds, see \cite[Chapter 11]{Boyer:Galicki}. Theorem \ref{thm:Main} allows us to obtain smooth $\spins$--manifolds from self-dual Einstein $4$-orbifolds in an analogous way.

Now, a self-dual Einstein $4$-orbifold $Q$ with positive scalar curvature yields complete $\spins$--metrics via Theorem \ref{thm:Main} if and only if there exists a smooth $8$-manifold $M$ arising as a circle orbibundle over $B=\Lambda^-T^\ast Q$. Indeed, note that $H^5_{orb}(B)=0$ since $B$ retracts onto $Q$ and therefore the necessary topological condition in Theorem \ref{thm:Main} is vacuous. We prove in Lemma \ref{lem:Spin7:admissible} that $B$ is the circle quotient of a smooth $8$-manifold if and only if $Q$ itself is the circle quotient of a smooth $5$-manifold; whenever this happens we say that $Q$ is \emph{$\spins$--admissible}. Theorem \ref{thm:Main} is useful only if we can find a large supply of $\spins$--admissible self-dual Einstein $4$-orbifolds with positive scalar curvature.

Infinitely many self-dual Einstein $4$-orbifolds with positive scalar curvature are known thanks to the quaternionic K\"ahler quotient construction of Galicki--Lawson \cite{Galicki:Lawson}. For example, infinitely many self-dual Einstein $4$-orbifolds with positive scalar curvature arise as quaternionic K\"ahler reductions of quaternionic projective space $\HH\PP^n$ by a subgroup of $\tu{Sp}(n+1)$. In \cite{Galicki:Lawson} Galicki--Lawson illustrate their quotient construction by considering self-dual Einstein metrics on weighted complex projective planes $\mathbb{W}\C\PP^2[q_1,q_2,q_3]$ arising as quotients of $\HH\PP^2$ by a circle. All these orbifolds are clearly $\spins$--admissible since weighted projective planes are all circle quotients of $S^5$. Theorem \ref{thm:Examples:Families} follows from applying our main existence result Theorem \ref{thm:Main} to these Galicki--Lawson examples.

It is likely that many more examples of self-dual Einstein $4$-orbifolds are \spins--admissible. For example, all toric self-dual Einstein $4$-orbifolds, \ie $4$-orbifolds with a $T^2$--symmetry, must arise as quaternionic K\"ahler quotients of $\HH\PP^n$ by an $(n-1)$-dimensional torus \cite{Calderbank:Singer}. The geometry of these toric orbifolds is then completely encoded in the combinatorics of the embedding of the Lie algebra of $T^{n-1}$ into the Lie algebra of the maximal torus of $\tu{Sp}(n+1)$. It is likely that combinatorial conditions characterising \spins--admissibility can be given in the same way that clear combinatorial criteria characterise the existence of a smooth $3$-Sasaki Konishi bundle \cite[Theorem 2.14]{Boyer:Galicki:Inventiones}.  Instead of pursuing such a systematic combinatorial approach, however, in this paper we construct by hand an explicit family of examples with unbounded second orbifold Betti number. In the proof of Theorem \ref{thm:Examples:Betti} we use an infinite list of self-dual Einstein $4$-orbifolds with positive scalar curvature arising from ALE gravitational instantons of type $A_n$ via the hyperk\"ahler/quaternionic K\"ahler correspondence. This correspondence associates to each hyperk\"ahler metric with a circle action that preserves only one complex structure in the twistor sphere an $S^1$--invariant quaternionic K\"ahler space of the same dimension. The examples of self-dual Einstein $4$-orbifolds we consider (originally considered by Galicki--Nitta \cite{Galicki:Nitta} without any reference to the hyperk\"ahler/quaternionic K\"ahler correspondence) give rise to $8$-manifolds that are rank-$3$ real vector bundles over $\sharp_{k}(S^2\times S^3)$ for any $k\geq 1$. Joyce's analytic constructions of compact $\spins$--manifolds \cite{Joyce:Spin7:I,Joyce:Spin7:II} can also be adapted to produce complete non-compact $\spins$--metrics by desingularising orbifold quotient singularities of non-compact flat orbifolds \cite[\S\S 13.1 and 15.1]{Joyce:Book} or asymptotically cylindrical Calabi--Yau $4$-folds \cite{Kovalev:ACyl:Spin7}. However, the variety of examples produced by Theorem \ref{thm:Examples:Betti} is new.

\subsubsection*{ALC $\gtwo$--manifolds from AC Calabi--Yau orbifolds}

The geometric and analytic framework to work on AC orbifolds we introduce in this paper allows us to extend the construction of ALC $\gtwo$--manifolds from AC Calabi--Yau $3$-folds in \cite{FHN:ALC:G2:from:AC:CY3} to the orbifold case. While \cite{FHN:ALC:G2:from:AC:CY3} already yields infinitely many examples of complete non-compact $\gtwo$--manifolds, as a simple application of our orbifold extension we produce infinitely many distinct families of ALC $\gtwo$--metrics on $S^3\times\R^4$, see Theorem \ref{thm:ALC:G2}. As in Theorem \ref{thm:Examples:Families} the families are distinguished by their tangent cones at infinity. In order to prove Theorem \ref{thm:ALC:G2} we show that there are infinitely many ways of realising $S^3\times\R^4$ as a circle orbibundle over a small orbifold partial resolution of a Gorenstein toric K\"ahler cone. For example, there is an infinite list of $S^1$--actions on $S^3\times\R^4$ labelled by two coprime positive integers $p,q$ such that $B=S^3\times\R^4/S^1_{p,q}$ is a small partial resolution of the Calabi--Yau cone over the so-called $Y^{p,q}$ Sasaki--Einstein $5$-manifold \cite{GMSW}. AC Calabi--Yau metrics on $B$ are constructed by Martelli--Sparks \cite{Martelli:Sparks:Partial:Resolutions} using the formalism of Hamiltonian $2$-forms. The construction of \cite{FHN:ALC:G2:from:AC:CY3}, suitably extended to the orbifold setting using the analysis of Section \ref{sec:AC:Orbifolds} in this paper, then immediately yields families of highly collapsed ALC $\gtwo$--metrics on $S^3\times\R^4$.    

\subsubsection*{Plan of the paper} The rest of the paper is organised in three main sections corresponding to the three different aspects of the proof of Theorem \ref{thm:Main} and of its applications. Section \ref{sec:AC:Orbifolds} develops the necessary geometric and analytic framework to work on AC orbifolds and includes the proof of Theorem \ref{thm:L2:cohomology} about weighted $L^2$--cohomology of AC orbifolds. Section \ref{sec:Adiabatic} contains the proof of Theorem \ref{thm:Main} implementing the adiabatic limit strategy we have outlined. Finally, Section \ref{sec:Examples} presents the concrete examples of Theorems \ref{thm:Examples:Betti} and \ref{thm:Examples:Families}.

\subsection*{Acknowledgements}
The author wishes to thank the Royal Society for the support of his research under a University Research Fellowship and NSF for the partial support of his work under grant DMS-1608143. The author would also like to thank Mark Haskins for comments on an early version of this paper, for introducing him to \cite{Martelli:Sparks:Partial:Resolutions} and for discussions over the years about self-dual Einstein $4$-orbifolds and their twistor spaces and about collapse of ALC metrics with exceptional holonomy. The author also thanks Olivier Biquard and Paul Gauduchon for suggesting the connection with the HK/QK correspondence in Remark \ref{rmk:HK/QK}. 

\section{Asymptotically conical orbifolds}\label{sec:AC:Orbifolds}

In this section we develop the necessary geometric and analytic framework to work on AC orbifolds. For the geometric applications of the paper it will suffice to consider orbifolds arising as quotients of smooth manifolds by a circle action. While we will restrict to this situation for ease of exposition, note that every (effective) orbifold arises as the quotient of a smooth manifold by an effective \emph{almost-free} (\ie with finite stabilisers) action of a compact Lie group \cite[Corollary 1.24]{Ruan:Orbifolds:Book}. In Section \ref{sec:AC:Orbifolds:Basics} we collect preliminary materials on orbifolds and foliations and introduce the language we are going to use in the rest of the paper. In Section \ref{sec:Transversally:AC} we develop a good Fredholm theory for linear elliptic operators on AC orbifolds. We apply this theory in Section \ref{sec:L2:cohomology} to provide a computation of the weighted $L^2$--cohomology of AC orbifolds.

\subsection{Orbifolds and principal Seifert circle bundles}\label{sec:AC:Orbifolds:Basics}

As a preliminary, we collect the facts about orbifolds and foliations that we are going to use throughout the section. We use the traditional notion of an orbifold, \ie an \emph{effective} orbifold in the sense of \cite[Definitions 1.1 and 1.2]{Ruan:Orbifolds:Book}, and avoid almost completely the language of groupoids. Indeed, the orbifolds we will consider all arise as quotients of an effective \emph{almost-free} (\ie with finite stabilisers) circle action on a smooth manifold. Our exposition uses this fact in an essential way.

Let $\xi$ be a nowhere-vanishing vector field on a smooth manifold $M$ of dimension $n+1$. Assume that the orbits of $\xi$ are all closed, \ie (possibly after an appropriate normalisation) $\xi$ generates an effective almost-free circle action on $M$. The orbit space $B=M/S^1$ has a natural \emph{orbifold} structure and $\pi\co M\ra B$ is a principal circle \emph{orbibundle}. We refer to \cite[Chapters 1 and 2]{Ruan:Orbifolds:Book} and \cite[Chapter 4]{Boyer:Galicki} for basics on orbifolds and orbibundles. Since its total space is smooth, $\pi\co M\ra B$ is a \emph{Seifert fibration} in the sense of \cite[Definition 1.2]{Goette}. An alternative viewpoint is that the vector field $\xi$ defines a foliation on $M$. We are going to use various notions from the theory of foliations \cite{Molino,Reinhart,Tondeur:Book}.

We fix a Riemannian metric $g$ on $M$ such that (a) $\xi$ has unit length, and (b) the orbits of $\xi$ are geodesics. By \cite[Proposition 6.7]{Tondeur:Book} such a metric exists if and only if there exists a $1$-form $\theta$ on $M$ such that $\theta(\xi)=1$ and $\mathcal{L}_\xi\theta=0$. Denote by $\mathcal{H}$ the \emph{horizontal} bundle $\ker\theta=\xi^\perp$. Observe that $\mathcal{H}$ can be identified with the pull-back to $M$ of the orbifold tangent bundle of $B$. The restriction of $g$ to $\mathcal{H}$ will be denoted by $g_B$ since it defines a Riemannian metric on the orbifold $B$. We will refer to the data $(M,\pi,\theta,g_B)$ as a \emph{Riemannian principal Seifert (circle) bundle}. 

We will always assume that $M$ is oriented with volume form $\dvol_g = \theta\wedge \dvol_B$, where $\dvol_B=\xi\lrcorner\dvol_g$ is a nowhere-vanishing section of $\Lambda^n\mathcal{H}^\ast$ satisfying $\mathcal{L}_\xi\dvol_B=0$. We will denote by $\ast_M$ the Hodge-star operator of $(M,g,\dvol_g)$.

\subsubsection{Projectable bundles and connections}

Let $\pi\co M\ra B$ be a Riemannian principal Seifert circle bundle and let $P\ra M$ be a principal $G$--bundle, where $G$ is a compact Lie group. We say that $P$ is \emph{projectable} if the circle action on $M$ lifts to a circle action on $P$ commuting with the $G$--action. Projectable principal bundles on $M$ are in one-to-one correspondence with principal orbibundles on the orbifold $B$. In the theory of foliations there is a weaker notion of a foliated bundle \cite[\S 2.6]{Molino}, where one only assumes that the vector field $\xi$ lifts to a vector field $\tilde{\xi}$ on $P$. The restriction of a foliated bundle $P$ to an orbit $\mathcal{O}\simeq S^1$ of $\xi$ is a trivial principal $G$--bundle endowed with a flat connection. If $P$ is projectable then this flat connection has trivial holonomy.

Let $V$ be a $G$--representation and consider the associated vector bundle $E=P\times_G V\ra M$. If $P$ is foliated with lift $\tilde{\xi}$ of $\xi$, we say that a section $s\co M\ra E$ of $E$ is \emph{basic} if $\mathcal{L}_{\tilde{\xi}}\tilde{s}=0$, where $\tilde{s}\co P\ra V$ is the $G$--equivariant function corresponding to $s$. If $P$ is projectable we can interpret basic sections as sections of the orbibundle $E/S^1\ra B$. For this reason we will denote the space of basic smooth sections of $E$ by $C^\infty (B;E)$. When the circle action on $M$ is free and $B$ is a manifold then $C^\infty (B;E)$ coincides with the space of smooth sections of the bundle $E/S^1\ra B$. Spaces of basic sections with lower regularity (for example $L^2$ sections) are defined in a similar way.

A connection $A$ on $P$, thought of as a Lie algebra-valued $1$-form on $P$, is \emph{projectable} if $\tilde{\xi}\lrcorner A=0=\tilde{\xi}\lrcorner F_A$. Doing analysis on $M$ with projectable connections acting on basic sections is a replacement for doing analysis on the orbifold $B$ without worrying about its singularities.

The oriented orthonormal frame bundle of $M$ is projectable but the Levi--Civita connection $\nabla^{\tu{LC}}$ of $g$ is not. Following \cite[Definition 3.13]{Reinhart} (see also \cite[Definition 1.7]{Bismut} in the case where $B$ is smooth and $M\ra B$ is an arbitrary fibration) we will introduce an adapted connection $\nabla$ which is better suited to the Seifert fibration structure than the Levi-Civita connection. Let $\xi$ be the vertical vector field generating the circle action on $M$ and let $\theta$ be its dual $1$-form. Let $X,Y$ denote vectors in $\mathcal{H}$. The Levi-Civita connection $\nabla^{\tu{LC}}$ of $g$ is
\[
\nabla^{\tu{LC}} _\xi \xi=0, \quad \nabla^{\tu{LC}} _X\xi=\tfrac{1}{2}(X\lrcorner d\theta)^\sharp, \quad \nabla^{\tu{LC}} _\xi Y=\tfrac{1}{2}(Y\lrcorner d\theta)^\sharp +[\xi,Y], \quad \nabla^{\tu{LC}} _X Y = -\tfrac{1}{2}d\theta (X,Y)\, \xi + \left( \nabla^{\tu{LC}} _X Y\right)_{\mathcal{H}}.
\]
Note that since $\theta([\xi,X])=-d\theta (\xi,X)=0$ we have $[\xi,\mathcal{H}]\subset\mathcal{H}$. We now define
\begin{equation}\label{eq:Adapted}
\nabla_\xi \xi=0, \quad \nabla_X\xi=0, \quad \nabla_\xi Y=[\xi,Y], \quad \nabla_X Y = \left( \nabla^{\tu{LC}} _X Y\right)_{\mathcal{H}}.
\end{equation}
We will refer to $\nabla$ as the \emph{adapted connection} of the Seifert bundle $\pi\co M\ra B$. The adapted connection is a projectable metric connection, but it has non-vanishing torsion $T(U,V)=d\theta (U,V)\, \xi$. A vector field $X$ on $M$ is called \emph{basic} if it is a basic section of $\mathcal{H}$, \ie $X\in\mathcal{H}$ and $[\xi,X]=0$. Note that $\nabla_\xi X=0$ for every basic vector field. Basic vector fields are identified with vector fields on the orbifold $B$ and under this identification the adapted connection $\nabla$ corresponds to the Levi--Civita connection of $g_B$.

\subsubsection{Transverse elliptic operators}

Let $E\ra M$ be a projectable metric bundle endowed with a projectable metric connection $A$. Combing $A$ with the adapted connection $\nabla$ on $M$ we obtain a projectable connection, still denoted by $\nabla$, on any tensor bundle with values in $E$. We can then use $\nabla$ to define differential operators on $M$ acting on $E$--valued tensors.

A \emph{basic tensor} is an $S^1$--invariant section of $\bigotimes^r\mathcal{H}\otimes\bigotimes^s\mathcal{H}^\ast$. Since $\mathcal{L}_\xi$ coincide with $\nabla_\xi$, the adapted connection preserves basic tensors. Hence if $P$ is a differential operator defined using the adapted connection and acting on sections of (a sub-bundle of) $\bigotimes^r TM\otimes \bigotimes ^s T^\ast M\otimes E$, then the restriction of $P$ to basic $E$--valued tensors is well-defined. We will refer to the restriction of $P$ to basic tensors as a \emph{transverse} (or basic) operator. A basic operator is \emph{elliptic} if its extension as an operator acting on arbitrary $E$--valued tensors is elliptic. 

We are particularly interested in ``basic versions'' of $d+d^\ast$, the rough Laplacian and the Dirac operator $D$ acting on differential forms and spinors on $M$ with values in $E$. Fix an orthonormal frame $e_1,\dots,e_{n+1}$ for $(M,g)$. We will assume that $\{ e_1,\dots,e_{n+1}\}$ is an \emph{adapted frame}, \ie $e_1,\dots, e_n$ are basic vector fields and $e_{n+1}=\xi$. We then define
\begin{equation}\label{eq:Adapted:d+dstar}
d_{\nabla} =\sum_{i=1}^{n+1}{e_i\wedge \nabla_{e_i}}, \qquad d_{\nabla}^{\ast} = -\sum_{i=1}^{n+1}{e_i\lrcorner\nabla}_{e_i}, \qquad \nabla^\ast\nabla = -\sum_{i=1}^{n+1}{\nabla_{e_i}\nabla_{e_i}}
\end{equation}
acting on $E$--valued differential forms and arbitrary $E$--valued tensors, respectively.

As the notation suggests, $d_\nabla^\ast$ is the formal $L^2$--adjoint of $d_\nabla$, where the $L^2$--inner product on forms is defined using the metric $g$ and the volume form $\theta\wedge\dvol_B$. We want to understand the restriction of $d_\nabla$ and $d_\nabla^\ast$ to basic forms. According to our definition of basic tensors, a differential form $\gamma$ on $M$ is basic if and only if $\xi\lrcorner\gamma = 0=\mathcal{L}_\xi\gamma$. Let $\Omega^\bullet (B)$ denote the space of smooth basic forms. We now define a transverse Hodge-star operator $\ast$ by
\begin{equation}\label{eq:Transverse:Hodge:star}
\ast\gamma = \ast_M (\theta\wedge\gamma),
\end{equation}
for every basic form $\gamma$. Note that we also have $\ast_M\gamma =(-1)^k\,  \theta\wedge\ast\gamma$ if $\gamma\in\Omega^k(B)$. The following lemma follows from straightforward manipulations of \eqref{eq:Adapted:d+dstar} using the relations between $\ast_M$ and $\ast$ and between $\nabla^{\tu{LC}}$ and the adapted connection $\nabla$.

\begin{lemma}\label{lem:d:dast:basic}
For every $\gamma\in \Omega^k (M)$ we have
\[
d_\nabla\gamma = d\gamma - d\theta\wedge (\xi\lrcorner\gamma), \qquad d_\nabla^\ast\gamma = d_M^{\ast}\gamma - (-1)^{k(n+1-k)}\,\theta\wedge\ast_M (d\theta\wedge\ast_M\gamma).
\]
Here $d_M^{\ast}$ denotes the codifferential on $(M,g,\dvol_g)$. In particular, if $\gamma\in \Omega^k (B)$ is basic then
\[
d_\nabla\gamma = d\gamma, \qquad d_\nabla^\ast\gamma = (-1)^{n(k-1)+1}\, \ast\, d\, \ast\gamma.
\]
\proof
The formulas for $d_\nabla$ and its restriction to basic forms are immediate. The formula for $d_\nabla^\ast$ is deduced from the formula for $d_\nabla$ using the fact that $d_\nabla^\ast$ is the formal $L^2$--adjoint of $d_\nabla$. The description of the restriction of $d_\nabla^\ast$ to basic forms uses the fact that
\[
d^\ast_M \gamma = (-1)^{n(k-1)+1}\ast d\ast \gamma + (-1)^{k(n+1)}\theta\wedge \ast (d\theta\wedge \ast \gamma),\qquad \ast_M (d\theta\wedge\ast_M\gamma) = (-1)^{k}\ast (d\theta\wedge\ast\gamma)
\] 
if $\gamma$ is a basic $k$-form.
\endproof
\end{lemma}
By abuse of notation, in the rest of the paper we use the notation $d=d_\nabla|_{\Omega^\bullet(B)}$ and $d^\ast = d_\nabla^\ast|_{\Omega^\bullet(B)}$. In particular, we will say that a basic form is \emph{coclosed} if $d_\nabla^\ast\gamma=0$. Similarly, we will denote by $\triangle$ the restriction of $d_\nabla d_\nabla^\ast + d^\ast_\nabla d_\nabla$ to basic forms and say that a basic form $\gamma$ is \emph{harmonic} if $\triangle\gamma=0$.

\begin{remark*}
When $M$ is closed every basic harmonic form is closed and coclosed (in the sense we have just defined), but this is not necessarily the case if $M$ is not compact. We will therefore always keep a distinction between harmonic and closed and coclosed (basic) forms. 
\end{remark*}

In order to define the adapted Dirac operator $\slashed{D}$ we need to assume that $M$ is spin. The spin structure might not be projectable but the associated $\tu{Spin}^c$--structure always is (since the frame bundle is projectable). Since every complex representation of $\tu{Spin}(n)$ is also a representation of $\tu{Spin}^c (n)$, we can always define a twisted Dirac operator acting on spinors with values in a Hermitian vector bundle by 
\begin{equation}\label{eq:Adapted:Dirac}
\slashed{D}\psi = \sum_{i=1}^{n+1}{\gamma (e_i)\nabla_{e_i}\psi},
\end{equation}
where $\gamma$ denotes Clifford multiplication. The restriction of $\slashed{D}$ to basic (complex) spinors plays the role of the Dirac operator of the orbifold $B$.

The fact that we defined basic elliptic operators as the restriction to basic tensors of elliptic operators on $M$ allows us to extend basic properties of elliptic operators on compact manifolds (elliptic regularity estimates, properties of the spectrum, etc.) to transversally elliptic operators.

\begin{prop}\label{prop:Basic:Eigensections}
Let $M$ be a closed Riemannian principal Seifert circle bundle with orbit space $B$. Let $P\co C^\infty (B;E)\ra C^\infty (B;F)$ be a self-adjoint basic elliptic operator. Then $P$ has discrete spectrum $\lambda_1 \leq \lambda_2 \leq \dots$ and there exists an orthonormal basis of $L^2(B;E)$ consisting of eigensections of $P$. Moreover, every eigensection of $P$ is smooth.  
\end{prop}

The second advantage of introducing operators based on the adapted connection is that local computations with basic tensors coincide with local computations in the standard Riemannian case as if $B$ were smooth. In particular, one can define \emph{transverse (or basic) curvature tensors} and prove Weitzenb\"ock formulas relating squares of Dirac-type operators (such as $\slashed{D}$ and $d_{\nabla}+d^\ast_{\nabla}$) and the rough Laplacian $\nabla^\ast\nabla$. Vanishing results and eigenvalue estimates based on positivity properties of curvature terms in these Weitzenb\"ock formulas are then deduced in the usual way; see \cite[Chapter 8]{Tondeur:Book} for details about this technique. The following proposition is an example of the results obtained using this method.

\begin{remark*}
Strictly speaking Tondeur replaces the Levi--Civita connection of $g$ with a different choice of connection \cite[Equation (8.1)]{Tondeur:Book} than the adapted connection \eqref{eq:Adapted}. However, both connections satisfy the key property of \cite[Proposition 8.6]{Tondeur:Book}.  
\end{remark*}

\begin{prop}\label{prop:Eigenvalue:Estimates}
Let $M^{n+1}$ be a closed Riemannian principal Seifert circle bundle with orbit space $B$ and assume that the transverse Ricci-curvature $\tu{Ric}(g_B)$ satisfies $\tu{Ric}(g_B)\geq (n-1) g_B$.
\begin{enumerate}
\item The first non-zero eigenvalue of the Laplacian $\triangle$ acting on basic functions is greater than or equal to $n$.
\item The first eigenvalue of the Laplacian $\triangle$ acting on coclosed basic $1$-forms is greater than or equal to $2(n-1)$ and the eigenspace with eigenvalue $2(n-1)$ consists of basic $1$-forms dual to basic Killing vector fields which are also eigenvectors for the transverse Ricci-curvature with eigenvalue $2(n-1)$.
\end{enumerate}
\proof
When $B$ is smooth, part (i) is the classical Lichnerowicz--Obata Theorem and part (ii) is less well known but also classical, see  \cite[Theorem 7.6]{Cheeger:Tian:AC:Ricci:flat} or \cite[Lemma B.2]{Hein:Sun}. In light of the remarks before the proposition the same proof extends to the case where $B$ is singular.  
\endproof
\end{prop}

\begin{remark}\label{rmk:Obata}
There is also an analogue of Obata's rigidity result for the eigenvalue estimate in (i): if there is a non-trivial eigenfunction with eigenvalue $n$ then Shioya's orbifold version of Obata's Theorem \cite[Theorem 1.1]{Shioya} implies that $B$ is isometric to a finite quotient of the round $n$-sphere.
\end{remark}

\subsubsection{Basic cohomology}

Absolute de Rham cohomology of a (not necessarily closed) manifold $M$ is the cohomology of the differential complex $(\Omega^\bullet (M), d)$, where $\Omega^\bullet(M)$ is the space of smooth differential forms on $M$. Basic cohomology is similarly defined using the complex of basic forms. Indeed, note that $d$ preserves basic forms and therefore $(\Omega^\bullet (B), d)$ is a differential chain complex, whose cohomology is called the basic cohomology of $M$. (Note that $d_\nabla|_{\Omega^\bullet (B)}= d|_{\Omega^\bullet (B)}$ by Lemma \ref{lem:d:dast:basic} so there is no ambiguity here on the meaning of $d$.) We will denote the basic cohomology of $M$ by $H^\bullet( B)$, since when $B$ is smooth it coincides with the de Rham cohomology of the quotient manifold $B$. We define the compactly supported basic cohomology $H^\bullet_c(B)$ of $M$ in an analogous way.

\begin{remark}\label{eq:Basic:Gysin}
There is a natural Gysin sequence relating the cohomology of $M$ with its basic cohomology \cite[Theorem 6.13]{Tondeur:Book}. Indeed, there is an exact sequence of complexes
\[
0\ra\Omega^\bullet (B)\ra \Omega^\bullet_{S^1}(M)\ra \Omega^{\bullet-1} (B)\ra 0,
\]
where $\Omega^\bullet_{S^1}(M)$ is the space of $S^1$--invariant forms on $M$ and the second map is contraction with $\xi$. The long exact sequence in cohomology replaces the Gysin sequence of a circle fibration since the cohomology of $(\Omega^\bullet_{S^1}(M),d)$ is isomorphic to the standard de Rham cohomology of $M$ by averaging along the (compact) orbits of $\xi$.
\end{remark}

The following proposition discusses the topological interpretation of basic cohomology.

\begin{prop}\label{prop:Basic:Singular:Coho}
Let $M$ be a Riemannian principal Seifert circle bundle with orbit space $B$. The basic cohomology $H^\bullet(B)$ of $M$ is isomorphic to the equivariant cohomology of $M$, denoted by $H^\bullet_{orb}(B;\R)$. In particular, if the least common multiple of the orders of the finite stabilisers of points in $M$ is finite, then $H^\bullet(B)$ is isomorphic to the singular cohomology with real coefficients of the topological space $B$.
\proof
The proposition is a chain of isomorphisms between different cohomology theories for orbifolds and manifolds with a group action. First of all, by H. Cartan's generalised Chern--Weil theory the basic cohomology of $M$ is equivalent to the Cartan model for the equivariant cohomology of $M$ \cite[Chapter 5]{Guillemin:Sternberg}. The equivariant version of the de Rham Theorem \cite[Chapters 1--4]{Guillemin:Sternberg} states that the Cartan model is equivalent to Borel's topological construction of the equivariant cohomology of $M$ as the singular cohomology with real coefficients of $ES^1\times_{S^1}M$, which we denote by $H^\bullet_{orb}(B;\R)$. In order to explain the notation note that the (Haefliger) orbifold cohomology (with arbitrary coefficients) $H^\bullet_{orb}(B)$ of an orbifold $B$ is defined as the singular cohomology of the classifying space of the (unique up to Morita equivalence) groupoid associated to $B$, see \cite[p. 38]{Ruan:Orbifolds:Book} and \cite[Definition 4.3.6]{Boyer:Galicki}; the classifying space of a global quotient orbifold $M/G$ is equivalent to $EG\times_G M$ \cite[Example 1.53]{Ruan:Orbifolds:Book}. The last statement uses the Leray spectral sequence of the fibration $ES^1\times_{S^1}M\ra M/S^1=B$ \cite[Corollary 4.3.8]{Boyer:Galicki}.
\endproof
\end{prop}

We continue this topological parenthesis with two further observations. First of all, the isomorphism classes of principal circle orbibundles (not necessarily Seifert bundles) $\pi\co M\ra B$ are classified by the orbifold first Chern class $c^{orb}_1(M)\in H^2_{orb}(B;\Z)$, see \cite[Theorem 4.3.15]{Boyer:Galicki}. Secondly, we can define orbifold homotopy groups $\pi^{orb}_i(B)$ of an orbifold $B$ as the homotopy groups of the classifying space of the associated groupoid. A result of Thurston yields an interpretation of the orbifold fundamental group $\pi^{orb}_1(B)$ as the group of deck transformations of an orbifold universal cover, see \cite[Theorem 4.3.19]{Boyer:Galicki}. Note also that there is an exact sequence of homotopy groups associated with a Seifert circle bundle $\pi\co M\ra B$  \cite[Theorem 4.3.18]{Boyer:Galicki}
\begin{equation}\label{eq:Homotopy}
\dots \ra \pi^{orb}_2 (B)\ra \Z \ra \pi_1 (M) \ra \pi^{orb}_1(B)\ra 1. 
\end{equation}
Here the map $\pi^{orb}_2 (B)\ra \Z\simeq \pi_1 (S^1)$ is determined by the image of $c^{orb}_1(M)$ in $H^2_{orb}(B;\R)$ (\ie the orbifold first Chern class modulo torsion). 

\begin{remark}\label{rmk:Myers}
There is an orbifold version of Bonnet--Myers' Theorem: if the transverse Ricci-curvature of a complete Riemannian Seifert bundle $\pi\co M\ra B$ is strictly positive then $\pi^{orb}_1(B)$ is finite. Indeed, by \cite[Corollary 21]{Myers:Orbifolds} the diameters of $B$ and its orbifold universal cover must then be bounded. 
\end{remark}

Returning to basic cohomology, we conclude with a discussion of the basic version of Hodge theory. Exploiting the fact that $d_\nabla + d^\ast_\nabla$ is an elliptic operator, the same reasoning that lead us to Proposition \ref{prop:Basic:Eigensections} allows us to deduce a basic Hodge theorem, see \cite[Theorem 7.22]{Tondeur:Book}.

\begin{prop}\label{prop:Basic:Hodge}
Let $M$ be a closed Riemannian principal Seifert circle bundle with orbit space $B$. Denote by $\mathcal{H}^k(B)$ the space of basic closed and coclosed forms, \ie
\[
\mathcal{H}^k(B)= \{ \gamma\in\Omega^k(B)\, |\, d\gamma=0=d^\ast\gamma\}.
\]
Then the map that assigns to each closed and coclosed basic form its basic cohomology class is an isomorphism $\mathcal{H}^k(B)\simeq H^k(B)$.
\end{prop}

\subsection{Seifert bundles that are transversally AC}\label{sec:Transversally:AC}

Let $N$ be a connected, oriented, closed $n$-manifold and $\pi_\infty\co N\ra \Sigma$ be a Riemannian principal Seifert circle bundle over a closed $(n-1)$-orbifold $\Sigma$. Denote by $\xi_\infty,\theta_\infty$ and $g_\Sigma$ the choice of a nonsingular vector field, dual $1$-form and horizontal metric on $N$. In particular, $N$ is endowed with the Riemannian metric $g_N = g_\Sigma + \theta_\infty^2$.

The Riemannian cone $\tu{C}(N)$ over a Riemannian manifold $(N,g_N)$ is $\R^+\times N$ endowed with the (incomplete) Riemannian metric $g_{\tu{C}}=dr^2 + r^2 g_N$. Instead of this conical metric, exploiting the Seifert bundle structure of $N$ we will consider $\tu{BC}(N)=\tu{BC}(N,\pi,\theta_\infty,g_\Sigma)=\R^+\times N$ endowed with the \emph{b}undle-like transversally \emph{c}onical metric
\begin{equation}\label{eq:BC:metric}
g_{\tu{BC}} = dr^2 + r^2 g_\Sigma + \theta_\infty^2.
\end{equation}

Let $(M,g)$ be a complete connected oriented Riemannian manifold with only one end. We assume that $\pi\co M\ra B$ is a Riemannian principal Seifert circle bundle with generating vector field $\xi$, connection $1$-form $\theta$ and horizontal metric $g_B$. Here $B$ denotes the orbit space $M/S^1$.

\begin{definition}\label{def:Transv:AC}
We say that $(M,\pi,\theta,g_B)$ is \emph{transversally asymptotically conical (AC)} asymptotic to $\tu{BC}(N)$ with rate $\nu<0$ if there exists a compact set $K\subset M$, a positive number $R>0$ and a diffeomorphism
\[
f\co \BC(N) \cap \{ r> R\} \ra M\setminus K
\]
such that for all $j\geq 0$
\begin{equation}\label{eq:Transv:AC}
|\nabla_{g_\BC}^j (f^\ast g - g_\BC) |_{g_\BC}=O(r^{\nu-j}).
\end{equation}
\end{definition}

\begin{remark*}
Since $|d\theta_\infty|_{g_{\BC}}=O(r^{-2})$, here we can compute covariant derivatives using either the adapted or the Levi--Civita connection of $\BC (N)$ obtaining equivalent definitions. 
\end{remark*}

By averaging along the circle orbits, the diffeomorphism $f\co \BC(N) \cap \{ r> R\} \ra M\setminus K$ can be assumed to intertwine the circle actions, \ie $f_\ast\xi_\infty=\xi$. Indeed, since $g$ and $g_{\BC}$ are circle invariant this averaging procedure does not destroy the asymptotic decay of $f^\ast g$ to $g_{\BC}$. Note that the decay condition \eqref{eq:Transv:AC} then is equivalent to
\begin{equation}
|\nabla_{g_\BC}^j (\theta_\infty - f^\ast \theta) |_{g_\BC} + |\nabla_{g_\BC}^j (dr^2 + r^2 g_\Sigma - f^\ast g_B) |_{g_\BC}=O(r^{\nu-j}).\tag{\theequation${}^\prime$}
\end{equation}
Here the decay of $f^\ast\theta$ to $\theta_\infty$ allows to compare the horizontal metrics since the horizontal spaces are isomorphic for $r$ sufficiently large.

\subsubsection{Weighted Banach spaces}

Let $(E_\infty,h_\infty)$ be a projectable metric vector bundle on $N$ endowed with a projectable metric connection $\nabla_\infty$. Here $h_\infty$ is an $S^1$--invariant metric on the bundle $E_\infty$ and $\nabla_\infty$ preserves it. By abuse of notation we will use the same symbols to denote the pull-back of $(E_\infty,h_\infty,\nabla_\infty)$ to $\BC (N)$. 

\begin{definition}\label{def:AC:admissible:connection}
Let $(M,\pi,\theta,g_B)$ be a transversally AC principal Seifert circle bundle asymptotic to $\tu{BC}(N)$ with rate $\nu<0$. Let $(E,h,\nabla)$ be a projectable metric bundle and connection over $M$. We say that $(E,h,\nabla)$ is \emph{admissible} if, under the identification $f\co (R,\infty)\times N\ra M\setminus K$ of Definition \ref{def:Transv:AC}, there exists an $S^1$--equivariant bundle isomorphism $f^\ast E\simeq E_\infty$ such that $f^\ast h =h_\infty + h'$ and $f^\ast\nabla = \nabla_\infty + a$, where $(h',a)$ satisfy
\[
|\nabla_\infty ^j h'|_{g_\tu{BC}\otimes h_\infty}=O(r^{\nu-j}), \qquad |\nabla_\infty ^j a|_{g_\tu{BC}\otimes h_\infty}=O(r^{\nu-1-j}).
\]
\end{definition}

We will mostly be interested in (sub)bundles of $\bigotimes^r\mathcal{H}\otimes\bigotimes^s\mathcal{H}^\ast$, where $\mathcal{H}$ is the horizontal bundle of $M$. By \eqref{eq:Transv:AC}, any such bundle together with the metric induced by $g$ and the connection induced by the adapted connection \eqref{eq:Adapted} of $g$ is admissible.

\begin{remark*}
In Definition \ref{def:AC:admissible:connection} we used the same rate of decay $\nu$ of the transversally AC Seifert bundle for ease of exposition and because we will mostly be working with tensor bundles and connections induced by the adapted connection. This restriction is of course unnecessary.
\end{remark*}

Fix once and for all an extension of the radial function $f^\ast r$ from $M\setminus K$ to the interior of $M$. By abuse of notation we will denote this extension by $r$ and assume that $r\geq 1$ and $|\nabla^k r|$ is uniformly bounded for all $k\geq 1$.

\begin{definition}\label{def:AC:weighted:spaces}
Let $(E,h,\nabla)$ be an admissible bundle. For all $p\geq 1, k\in\N_0, \alpha\in (0,1)$ and $\nu\in\R$ we define the weighted Sobolev space $L^p_{k,\nu}(B)$ and the weighted H\"older space $C^{k,\alpha}_\nu(B)$ of basic sections of $E$ as the closure of $C^\infty _c(B;E)$ with respect to the norms
\[
\| u\|_{L^p_{k,\nu}} = \left( \sum_{j=0}^k{\| r^{-\frac{n}{p}-\nu+j}\nabla^j u\|^p_{L^p}}\right)^{\frac{1}{p}}, \qquad \| u\|_{C^{k,\alpha}_\nu} = \sum_{j=0}^k{\|r^{-\nu+j}\nabla^j u\|_{C^0} + [r^{-\nu+k}\nabla^k u]_\alpha}.
\]
By dropping the H\"older seminorm $[r^{-\nu+k}\nabla^k u]_\alpha$ in the definition of the $C^{k,\alpha}_\nu$--norm, we obtain the definition of the space of basic sections of $E$ of class $C^k_\nu$. Finally, define $C^\infty_\nu (B) = \bigcap_{k\geq 0}{C^k_\nu} (B)$.
\end{definition}

A standard technique to work with weighted Banach spaces on AC manifolds is the scaling and covering argument of \cite[Theorem 1.2]{Bartnik}. The same technique can be used in the more general context of transversally AC Seifert bundles. Decompose the region $\{ r\geq R \}$ in $\tu{BC}(N)$ into the union of ``annuli'' $\mathcal{A}_{2^k R, 2^{k+1}R}=\{ 2^k R \leq r \leq 2^{k+1}R\}$. Each region $\mathcal{A}_{2^k R, 2^{k+1}R}$ can be covered with the same finite number of open subsets of the form $(U\times S^1)/\Gamma$ for some finite group $\Gamma$. The fact that the number of open subsets is independent of the radius of the annulus follows from the fact that $\BC (N)$ (and all the structure it carries) is the radial extension of the compact Seifert bundle $\pi_\infty\co N\ra \Sigma$. Up to a factor of $(2^k R)^{-\nu}$, on each annulus the weighted Sobolev/H\"older norms of \emph{basic} sections are equivalent (with constants independent of $R$ and $k$) to the standard Sobolev/H\"older norms of basic sections on the fixed annulus $\{ 1\leq r\leq 2\}$. Then estimates on the exterior region $\{ r\geq R\}$ in $\BC (N)$ (and, via the identification $f$ of Definition \ref{def:Transv:AC}, in $M$) can be obtained by applying standard estimates for basics sections on these rescaled regions, rescaling back and summing/taking supremums over $k\in\Z_{\geq 0}$. Combined with interior estimates for basic sections on a compact set in $M$, this method yields estimates for basic sections on the whole manifold. 

The following embedding theorem can be proved using this strategy. The necessary local interior estimate is Parker's Equivariant Sobolev Embedding Theorem \cite[\S 1]{Parker}, which states that for basic sections the Sobolev inequalities work as if we were in dimension $n$ rather than $n+1$.

\begin{theorem}\label{thm:Weighted:Embedding}
Let $M$ be an $(n+1)$-dimensional transversally AC manifold. All Banach spaces below are spaces of basic sections of an admissible vector bundle $E$.
\begin{enumerate}
\item If $k\geq h\geq 0$, $k-\frac{n}{p}\geq h-\frac{n}{q}$, $p\leq q$ and $\nu\leq \nu'$ there is a continuous embedding $L^p_{k,\nu}(B)\subset L^q_{h,\nu'}(B)$. Moreover, if $k>h$, $k-\frac{n}{p}> h-\frac{n}{q}$ and $\nu<\nu'$ then the embedding is compact.
\item If $k-\frac{n}{p}\geq h+\alpha$ then there is a continuous embedding $L^p_{k,\nu}(B) \subset C^{h,\alpha}_\nu(B)$.
\end{enumerate}
\end{theorem}

\begin{remark}\label{rmk:Weighted:Embedding}
There are further results about embeddings and products that follow more easily from Definition \ref{def:AC:weighted:spaces}. For example the first statement below only uses H\"older's inequality.
\begin{enumerate}
\item If $k\geq h\geq 0$, $k-\frac{n}{p}\geq h-\frac{n}{q}$, $p> q$ and $\nu< \nu'$ there is a continuous embedding $L^p_{k,\nu}(B)\subset L^q_{h,\nu'}(B)$. Moreover, if $k>h$ and $k-\frac{n}{p}> h-\frac{n}{q}$ then the embedding is compact.
\item If $\nu\leq \nu'$ and $k+\alpha \geq h+\beta$ then there are continuous embeddings $C^{k+1}_\nu(B)\subset C^{k,\alpha}_\nu(B)\subset C^{h,\beta}_{\nu'}(B)\subset C^h_{\nu'}(B)$. Moreover, if $\nu<\nu'$ the embedding $C^{k,\alpha}_\nu(B)\subset C^h_{\nu'}(B)$ is compact.
\item If $\nu<\nu'$ there is a continuous embedding $C^{h,\alpha}_\nu(B) \subset L^q_{h,\nu'}(B)$ for every $q\geq 1$.
\item If $\nu_1 + \nu_2 \leq \nu$ then the product $C^{k,\alpha}_{\nu_1}(B)\times C^{k,\alpha}_{\nu_2}(B)\ra C^{k,\alpha}_\nu(B)$ is continuous.
\end{enumerate}
In the next section we will use (iv) to control the nonlinearities in the equations describing circle invariant $\spins$--metrics. 
\end{remark}

\subsubsection{Admissible operators}

Consider a transversally elliptic operator $P_\infty$ of order $k$ on $\BC (N)$. The fact that the conical metric $dr^2 +r^2 g_\Sigma$ is conformal to the cylindrical metric $dt^2 + g_\Sigma$, where $r=e^t$, motivates us to consider the rescaled operator $r^k P_\infty$. Using the conformal equivalence between cones and cylinders, it makes sense to require that $r^k P_\infty$ acting on basic sections is invariant under translations in $t$. 

\begin{definition}\label{def:Admissible:Operators}
Let $P\co C^\infty(B;E)\ra C^\infty (B;F)$ be a transversally elliptic operator of order $k$ between sections of admissible vector bundles over a transversally AC principal Seifert circle bundle $\pi\co M\ra B$ asymptotic to $\BC (N)$. Let $f\co \BC (N)\cap \{ r>R\} \ra M\setminus K$ be the identification of Definition \ref{def:Transv:AC}. Let $P_\infty\co C^\infty (f^\ast E)\ra C^\infty (f^\ast F)$ be a transversally elliptic operator on $\BC(N)$ such that $r^k P_\infty$ acting on basic sections is a translation invariant operator. We say that $P$ is an admissible operator asymptotic to $P_\infty$ if
\[
|\nabla _\infty^j \left( f^\ast (Pu) - P_\infty f^\ast u\right) |_{g_\BC\otimes h_\infty}=O(r^{-k+\nu-j}) 
\]
for some $\nu<0$ for every $j\geq 0$ and every smooth basic section $u$ of $E$ on $M\setminus K$.  
\end{definition}

By Definition \ref{def:AC:admissible:connection}, if $P\co C^\infty (E)\ra C^\infty (F)$ is an elliptic operator of order $k$ between admissible vector bundles defined as the composition of $\nabla ^k\co C^\infty (E)\ra C^\infty \left( \bigotimes ^k T^\ast M\otimes E\right)$ with a constant coefficient bundle map $\bigotimes ^k T^\ast M\otimes E\ra F$, then $P$ is admissible. In particular, the operator $d+d^\ast(=d_\nabla+d_\nabla^\ast)$ of \eqref{eq:Adapted:d+dstar}, the Laplacian $\triangle (= d_\nabla d_\nabla^\ast + d_\nabla^\ast d_\nabla)$ and the Dirac operator \eqref{eq:Adapted:Dirac} acting on basic spinors and differential forms on a transversally AC principal Seifert bundle are all admissible operators.

With this definition the theory of admissible transversally elliptic operators on transversally AC principal Seifert bundles acting on basic sections is identical to the standard theory of elliptic operators on AC manifolds. We will now briefly state the main analytic results we are going to use in the paper. We omit proofs as they are identical to the standard case of AC manifolds, for which we refer to \cite{Lockhart:McOwen,Lockhart, Melrose}, \cite[\S 4.3]{Marshall} and the brief summary in \cite[Appendix B]{FHN:ALC:G2:from:AC:CY3}.

First of all, the following integration-by-parts formula will be repeatedly used throughout the paper.
\begin{lemma}\label{lem:integration:parts}
Let $P\co C^\infty (B;E)\ra C^\infty (B;F)$ be an admissible operator of order $1$ and let $P^\ast$ its formal adjoint. Then for every $u\in L^2_{1,\nu}(B)$ and $v\in L^2_{1,\nu'}(B)$ with $\nu+\nu'\leq -n+1$ we have
\[
\langle Pu,v\rangle_{L^2} = \langle u, P^\ast v\rangle_{L^2}.
\]
\end{lemma}

The following elliptic regularity estimates can be proved using the scaling and covering technique discussed earlier.

\begin{theorem}\label{thm:Weighted:Regularity}
Let $P\co C^\infty (B;E)\ra C^\infty (B;F)$ be an admissible operator of order $k$. Then for every $l\geq 0$, $p\geq 1$, $\alpha\in (0,1)$ and $\nu\in \R$ there exists $C>0$ such that
\[
\begin{gathered}
\| u \|_{L^p_{l+k,\nu+k}} \leq C\left( \| Pu\|_{L^p_{l,\nu}} + \| u \|_{L^p_{0,\nu+k}}\right) ,\qquad \| u \|_{C^{l+k,\alpha}_{\nu+k}} \leq C\left( \| Pu\|_{C^{l,\alpha}_{\nu}} + \| u \|_{C^{0,\alpha}_{\nu+k}}\right),\\
\| u \|_{C^{l+k,\alpha}_{\nu+k}} \leq C \left( \| Pu\|_{C^{l,\alpha}_{\nu}} + \| u \|_{L^2_{\nu+k}}\right)
\end{gathered}
\]
for all basic section $u\in C^\infty_{c}(B)$.
\end{theorem}

\subsubsection{Fredholm theory for transversally elliptic operators}

In order to proceed further it is necessary to study in more details the mapping properties of the model operator $P_\infty$. This can be done explicitly by separation of variables.

\begin{definition}\label{def:indicial:roots}
Let $P_\infty$ be a transversally elliptic operator on $\BC(N)$ acting on basic sections of a projectable bundle $E_\infty\ra \BC (N)$ and assume that $r^k P_\infty$ is a translation invariant operator.
\begin{enumerate}
\item We say that a basic section $u$ of $E_\infty$ is \emph{homogeneous of order $\lambda$} if $\mathcal{L}_{r\partial_r}u=\lambda u$. In the particular case of differential forms, however, we adopt the convention that a basic $k$-form $\gamma$ on $\BC (N)$ is homogeneous of order $\lambda$ if $\mathcal{L}_{r\partial_r}\gamma=(\lambda+k)\gamma$.
\item We say that $\lambda$ is an \emph{indicial root} of $P_\infty$ if there exists a homogeneous basic section $u$ of rate $\lambda$ such that $P_\infty u=0$. We denote the set of indicial roots of $P_\infty$ by $\mathcal{D}(P_\infty)$.
\item For $\lambda\in \mathcal{D}(P_\infty)$ 
let $d(\lambda)$ denote the dimension of the space of basic sections $u\in \ker P_\infty$ of the form $u = \sum_{j=0}^m{u_j (\log{r})^j}$ with $u_0,\dots, u_m$ homogeneous of order $\lambda$.
\item For $\nu,\nu'\notin\mathcal{D}(P_\infty)$ with $\nu<\nu'$ set $N(\nu, \nu') = \sum_{\lambda\in \mathcal{D}(P_\infty)\cap (\nu,\nu')}{d(\lambda)}$
\end{enumerate}
\end{definition}
Because of the translation invariance, the indicial roots of $P_\infty$ are completely determined by the spectrum of a transversally elliptic operator $P_\infty|_N$ on the compact Seifert bundle $N$. In particular, Proposition \ref{prop:Basic:Eigensections} implies that $\mathcal{D}(P_\infty)$ is discrete.

In the rest of the paper we are going to make extensive use of the following results. First of all, separation of variables on $\BC (N)$ and decomposition of basic sections on $N$ into eigenspaces of $P_\infty|_N$ using Proposition \ref{prop:Basic:Eigensections} allow one to prove the following result about the asymptotic behaviour of solutions to $Pu=v$. 

\begin{prop}\label{prop:Decay:Solutions}
Let $P\co C^\infty (B;E)\ra C^\infty (B;F)$ be an admissible operator of order $k$ and fix $l\geq 0$, $\alpha\in (0,1)$ and $\nu,\nu'\in\R$ with $\nu<\nu'$ and $\nu,\nu'\notin\mathcal{D}(P_\infty)$. Set $N=N(\nu,\nu')$. Let $u_1,\dots, u_N$ be a basis of the space of basic sections $u\in \ker P_\infty$ of the form $u = \sum_{j=0}^m{\overline{u}_j (\log{r})^j}$ for basic sections $\overline{u}_0,\dots, \overline{u}_m$ of $E_\infty$ homogeneous of rate $\lambda\in (\nu,\nu')$.

There exists a compact set $K\subset M$ such that for every $v\in C^{0,\alpha}_{\nu-k}(B)$ with $v=Du'$ for some $u'\in C^{k,\alpha}_{\nu'}(B)$ there exist $a=(a_1,\dots,a_N)\in\R^N$ and a basic section $u$ on $M\setminus K$ of class $C^{k,\alpha}_{\nu}$ such that $u'|_{M\setminus K}=u+\sum_{i=1}^N{a_i\,u_i}$. Moreover, there exists a constant $C>0$ independent of $f,u,u',a$ such that
\[
\| u \|_{C^{k,\alpha}_{\nu}} + \|a\| \leq C\left( \| f\| _{C^{0,\alpha}_{\nu-k}} + \| u' \|_{C^{k,\alpha}_{\nu'}}\right).
\]
\end{prop}

The analysis of the asymptotic operator $P_\infty$ by separation of variables is then used to derive the main results about the Fredholm property and index of admissible operators acting on weighted H\"older spaces.

\begin{theorem}\label{thm:Fredholm}
Let $P\co C^\infty (B;E)\ra C^\infty (B;F)$ be an admissible operator of order $k$ and fix $l\geq 0$, $\alpha\in (0,1)$ and $\nu,\nu'\in\R$ with $\nu<\nu'$.
\begin{enumerate}
\item If $\nu \in \R \setminus \mathcal{D}(P_\infty)$ then there exists a compact set $K\subset M$ and a constant $C>0$ such that
\[
\| u \| _{C^{l+k,\alpha}_{\nu}} \leq C \left( \| Pu \| _{C^{l,\alpha}_{\nu-k}} + \| u \| _{L^2 (K)}\right)
\]
for all smooth basic section $u$. In particular, $P\co C^{l+k,\alpha}_{\nu}(B)\ra C^{l,\alpha}_{\nu-k}(B)$ is a Fredholm operator.
\item Assume that $\nu,\nu'\notin\mathcal{D}(P_\infty)$ and denote by $i(\nu)$ and $i(\nu')$ the indices of $P\co C^{k,\alpha}_{\nu}(B)\ra C^{0,\alpha}_{\nu-k}(B)$ and $P\co C^{k,\alpha}_{\nu'}(B)\ra C^{0,\alpha}_{\nu'-k}(B)$ respectively. Then 
\[
i(\nu')-i(\nu)=N(\nu,\nu'),
\]
where $N(\nu, \nu') = \sum_{\lambda\in \mathcal{D}(P_\infty)\cap (\nu,\nu')}{d(\lambda)}$.
\end{enumerate}
\end{theorem}

Finally, we describe the obstructions to solve the equation $Pu=v$.

\begin{prop}\label{prop:Obstructions:Pu=v}
Let $P\co C^\infty (B;E)\ra C^\infty (B;F)$ be an admissible operator of order $k$ and fix $\alpha\in (0,1)$ and $\nu\in\R$. Then for every $v\in C^{0,\alpha}_{\nu-k}(B)$ such that
\[
\langle v , \overline{u}\rangle _{L^2}=0
\]
for all $\overline{u}\in \ker P^\ast \cap C^\infty_{-n-\nu+k}(B)$, there exists $u\in C^{k,\alpha}_\nu (B)$ such that $Pu=v$ and
\[
\| u \|_{C^{k,\alpha}_{\nu}} \leq C\| Pu\|_{C^{0,\alpha}_{\nu-k}} .
\]
\end{prop}

\subsection{Basic weighted $L^2$--cohomology}\label{sec:L2:cohomology}

As an application of the theory we have introduced, we conclude this section with a calculation of the weighted basic $L^2$--cohomology of a transversally AC Riemannian principal Seifert circle bundle $\pi\co M\ra B$. The $L^2$--cohomology of a smooth AC manifold is well known, see \cite[Theorem 1.A]{HHM} and \cite[Example 0.15]{Lockhart}. For our applications later in the paper it will be important to understand the space of basic closed and coclosed forms that do not necessarily have an $L^2$--rate of decay; more precisely, we will need to determine all topological contributions to the space of basic closed and coclosed forms. This is not widely known even in the case of AC manifolds and we provide an elementary proof using the tools we have introduced.


Let $\pi\co M\ra B$ a transversally AC Riemannian principal Seifert circle bundle asymptotic to $\BC (N)$. Recall that the basic de Rham cohomology (with compact support) of $M$ and $N$ is defined as the cohomology of the de Rham complex of smooth basic differential forms (with compact support). As in \eqref{eq:Transv:AC} we identify the complement of a compact set in $M$ with an exterior region in $\BC (N)$ using a diffeomorphism $f$ such that $f_\ast\xi=\xi_\infty$. In particular, basic differential forms on $M$ pull back to basic differential forms on $\BC (N)$.

We introduce the main piece of topological data we will use. Regard $M$ as a manifold with boundary $N$, and similarly $B$ as a topological space with boundary $\Sigma$. Despite we do not require that the singularities of $B$ are contained in a compact set, the circle action on $M$ still has some finiteness properties: outside of a compact set the circle action on $M$ is determined by the circle action on the compact manifold $N$. In particular, the finiteness assumption in Proposition \ref{prop:Basic:Singular:Coho} is certainly satisfied and we can therefore deduce a long exact sequence in basic cohomology from the long exact sequence in singular cohomology (with real coefficients) for the pair $(B,\Sigma)$:
\begin{equation}\label{eq:Exact:Sequence:Boundary}
\cdots \ra H^{k-1}(\Sigma)\ra H^k_c(B)\ra H^k(B)\ra H^k(\Sigma)\ra \cdots
\end{equation}
We can in fact be completely explicit. The map $H^k_c(B)\ra H^k(B)$ is induced by the natural inclusion of compactly supported basic forms in the space of all smooth basic forms. The map $H^k(\Sigma)\ra H^{k+1}_c (B)$ is $[\tau]\mapsto [d(\chi\tau)]$ where $\chi$ is a basic cut-off function with $\chi\equiv 1$ outside a compact set. In order to define the map $H^k(B)\ra H^k(\Sigma)$ we use the following representation for basic cohomology classes on $M$.

\begin{lemma}\label{lem:Asymptotics:closed:forms}
Every basic cohomology class $[\sigma]\in H^k(B)$ can be represented by a smooth closed form $\sigma$ that decays as $r^{-k}$. More precisely, fix $R>0$ sufficiently large and embed $N$ in $M$ as $f(\{ r=R\})$. Let $\beta_0\in\mathcal{H}^k(\Sigma)$ be the basic harmonic representative on $N$ of the image of $[\sigma]$ in the basic cohomology of $N$ via the restriction map $H^k(B)\ra H^k(\Sigma)$. Then we can take $\sigma$ to be a smooth closed form with $\sigma = \beta_0$ outside a compact set.
\end{lemma}
Using the lemma, the map $H^k(B)\ra H^k(\Sigma)$ is explicitly defined by $[\sigma]\mapsto [\beta_0]$.

\proof
Choose a closed smooth basic representative $\sigma'$ of $[\sigma]$. Outside a compact set we can think of $\sigma'$ as a basic form defined on an exterior domain in $\tu{BC}(N)$. We can then write
\[
\sigma' = dr\wedge\alpha + \beta
\]
where $(\alpha,\beta)$ is a curve in $\Omega^{k-1}(\Sigma)\oplus\Omega^k(\Sigma)$ parametrised by $r\in [R,\infty)$ for some $R>0$. A straightforward calculation shows that $d\sigma'=0$ if and only if
\begin{equation}\label{eq:Asymptotics:closed:forms}
\partial_r\beta = d_\Sigma\alpha, \qquad d_\Sigma\beta =0.
\end{equation}
Here $d_\Sigma$ denoted the restriction of $d_\nabla$ on $N$ to basic forms, where $\nabla$ is the adapted connection of $\pi_\infty\co N\ra \Sigma$. We now use basic Hodge theory Proposition \ref{prop:Basic:Hodge} on $N$ to write $\beta = \beta_0 + d_\Sigma\gamma$, where $\beta_0$ is basic harmonic on $N$ and $\gamma$ is a coclosed basic $(k-1)$-form on $N$ depending smoothly on $r\geq R$. Note that $\beta_0$ is independent of $r$ because of the first equation in \eqref{eq:Asymptotics:closed:forms}. This same equation now yields $d_\Sigma (\partial_r\gamma -\alpha)=0$. Hence $\alpha = \alpha_0 + \partial_r\gamma$, for a smooth curve $\alpha_0$ in $\mathcal{H}^{k-1}(\Sigma)$ parametrised by $r\geq R$. We then consider the basic $(k-1)$-form $\Gamma$ on $[R,\infty)\times N$ defined by
\[
\Gamma = \int_R^r{\alpha_0(s)\, ds} + \gamma.
\]
Fix a basic cut-off function $\chi$ which vanishes for $r\leq R$ and is equal to $1$ on $r\geq R+1$ and define
\[
\sigma = \sigma' -d(\chi\Gamma).
\]
Then $d\sigma=0$ and $[\sigma]=[\sigma']\in H^k(B)$, $\sigma = \beta_0$ on $r\geq R+1$ since $d\Gamma = \sigma'-\beta_0$ and $[\beta_0]$ is the image of $[\sigma]$ in $H^k(\Sigma)$. 
\endproof

Fix $\nu\in\R$ and let $\mathcal{H}_\nu^k(B)$ be the space of basic closed and coclosed forms of class $L^2_\nu$. If $\nu$ is not an indicial root for the Laplacian on basic $k$-forms, then by Theorem \ref{thm:Weighted:Regularity} every form in $\mathcal{H}_\nu^k(B)$ is in fact of class $C^\infty_\nu$. Note that $-k$ is always an indicial root of the Laplacian acting on basic $k$-forms. Indeed, every basic closed and coclosed $k$-form on $N$ pulls back to a basic harmonic form on $\BC (N)$ homogeneous of order $-k$. Our main result is the identification of $\mathcal{H}_\nu^k(B)$ as we cross the indicial root $-k$.

\begin{theorem}\label{thm:L2:cohomology}
Let $\pi\co M^{n+1}\ra B$ be a transversally AC principal Seifert circle bundle asymptotic to $\BC(N)$, where $\pi_\infty\co N\ra \Sigma$ is a closed principal Seifert bundle. In the following statement $\delta>0$ is chosen sufficiently small so that the only indicial root in $[-k-\delta, -k+\delta]$ is $-k$.
\begin{enumerate}
\item If $k<\tfrac{n}{2}$ there are natural isomorphisms
\[
\mathcal{H}^k_{-k-\delta}(B)\simeq H^k_c(B), \qquad \mathcal{H}^k_{-k+\delta}(B)/\mathcal{H}^k_{-k-\delta}(B) \simeq \tu{im} \left( H^k(B)\ra H^k(\Sigma)\right).
\]
\item If $k=\tfrac{n}{2}$ (then $n$ is even) there are natural isomorphisms
\[
\mathcal{H}^k_{-k-\delta}(B)\simeq \tu{im}\left( H^k_c(B)\ra H^k(B)\right), \qquad \mathcal{H}^k_{-k+\delta}(B)/\mathcal{H}^k_{-k-\delta}(B) \simeq \tu{im} \left( H^k(B)\ra H^k(\Sigma)\right)^{\oplus 2}.
\]
\item If $k>\tfrac{n}{2}$ there are natural isomorphisms
\[
\mathcal{H}^k_{-k-\delta}(B)\simeq \tu{im}\left( H^k_c(B)\ra H^k(B)\right), \qquad \mathcal{H}^k_{-k+\delta}(B)\simeq H^k(B).
\]
\end{enumerate}
\end{theorem}
The rest of the section contains a proof of this theorem, which involves various steps.

We begin with the following calculation of excluded indicial roots. 
\begin{lemma}\label{lem:Excluded:indicial:roots:kforms}
Let $\pi_\infty\co N^{n+1}\ra\Sigma$ be a closed, connected, oriented Riemannian principal Seifert circle bundle and consider homogeneous basic forms on $\BC(N)$ in the sense of Definition \ref{def:indicial:roots}.
\begin{enumerate}
\item If $k\leq \tfrac{n}{2}-1$ there are no basic harmonic $k$-forms homogeneous of order $-n+k +2 < \lambda <-k$. Moreover every basic harmonic $k$-form homogeneous of order $\lambda=-k$ is closed and coclosed. 
\item If $k=\tfrac{n}{2}$ every basic harmonic $k$-form homogeneous of order $\lambda=-n+k=-k$ is closed and coclosed.
\item If $k< \tfrac{n}{2}$ there are no basic closed and coclosed $k$-forms homogeneous of order $-n+k< \lambda <-k$.

\item If $k\neq \tfrac{n}{2}$ every basic closed and coclosed $k$-form homogeneous of order $\lambda=-k$ is the pull-back of a basic harmonic $k$-form on $N$. If $k= \tfrac{n}{2}$ every basic closed and coclosed $k$-form homogeneous of order $\lambda=-k=-n+k$ is the pull-back of a basic harmonic $k$-form on $N$ or its image under the basic Hodge-star operator $\ast$.
\item Let $\gamma$ be a polynomial in $\log{r}$ with coefficients in the space of basic $k$-forms homogeneous of order $\lambda$. If $\gamma$ is basic harmonic then either $\gamma=\gamma_0$ is constant in $\log{r}$ or $\lambda=-\tfrac{n}{2}-1$ and $\gamma = \gamma_1 \log{r} + \gamma_0$ with $\gamma_0, \gamma_1$ basic harmonic $k$-forms homogeneous of order $\lambda$.
\end{enumerate}
\proof
With our definitions, $d_\nabla, d_\nabla^\ast$ and $d_\nabla d_\nabla^\ast + d_\nabla^\ast d_\nabla$ acting on basic forms on $\BC (N)$ are equivalent to the operators $d, d^\ast$ and $dd^\ast + d^\ast d$ acting on differential forms on a Riemannian cone $\tu{C}(\Sigma)$. A characterisation of harmonic and closed and coclosed homogeneous forms on a cone is given in \cite[Appendix A]{FHN:ALC:G2:from:AC:CY3}. The lemma follows immediately from these results using the fact that the Laplacian on $\Sigma$ is a non-negative operator.
\endproof
\end{lemma}
Since $\triangle\ast=\ast\triangle$, one deduces similar statements for $k\geq \tfrac{n}{2}$ by replacing $k$ with $n-k$.

\begin{remark*}
In view of Lemma \ref{lem:Excluded:indicial:roots:kforms} (iii) and (iv), for all $\delta>0$ sufficiently small we have $\mathcal{H}^k_{-k-\delta}(B)=L^2\mathcal{H}^k (B)$ if $k\leq \tfrac{n}{2}$ and $\mathcal{H}^k_{-k+\delta}(B)=L^2\mathcal{H}^k (B)$ if $k>\tfrac{n}{2}$. Therefore Theorem \ref{thm:L2:cohomology} includes the calculation of the basic $L^2$--cohomology of a transversally AC principal Seifert circle bundle.
\end{remark*}

\begin{lemma}\label{lem:Harmonic:CClosed}
Let $\gamma$ be a basic harmonic $p$-form of class $C^\infty_\lambda$. Assume that either
\begin{enumerate}
\item $\lambda< -\tfrac{n}{2}+1$, or
\item $p< \tfrac{n}{2}-1$ and $\lambda<-p$.
\end{enumerate}
Then $\gamma$ is closed and coclosed.
\proof
If $\gamma$ is a basic harmonic $p$-form of class $C^\infty_\lambda$ for some $\lambda<-\tfrac{n}{2}+1$, then the integration by parts formula Lemma \ref{lem:integration:parts} implies
\[
0=\langle \triangle \gamma, \gamma\rangle _{L^2} = \| d_\nabla\gamma\|^2_{L^2} + \| d_\nabla^\ast\gamma\|^2_{L^2}.
\]
If $p< \tfrac{n}{2}-1$ then by Lemma \ref{lem:Excluded:indicial:roots:kforms} (i) and elliptic regularity every basic harmonic $p$-form in $C^\infty_\lambda$ with $\lambda<-p$ is in fact in $C^\infty_{-n+p+2+\epsilon}$ for every $\epsilon>0$. Choose $\epsilon$ sufficiently small so that $2p\leq n-2-2\epsilon$. Then $2(-n+p+2+\epsilon)-1\leq -n+1$ and we can apply the first part of the proof.
\endproof
\end{lemma}

We can now prove the first and third part of Theorem \ref{thm:L2:cohomology}.

\proof[Proof of Theorem \ref{thm:L2:cohomology} (i) and (iii)]
Fix $k<\tfrac{n}{2}$ and $\delta>0$ as in the statement of the theorem and set $\nu=-k-\delta$. Since $\nu$ is not an indicial root, every $\sigma\in \mathcal{H}^k_\nu (B)$ is in fact of class $C^\infty_\nu(B)$ by weighted elliptic regularity. As in the proof of Lemma \ref{lem:Asymptotics:closed:forms}, outside a compact set we write $\sigma = dr\wedge\alpha + \beta$, where $(\alpha,\beta)$ is a smooth curve in $\Omega^{k-1}(\Sigma)\oplus\Omega^k(\Sigma)$ parametrised by $r\in [R,\infty)$ for some $R>0$ and satisfying $r|\alpha|+|\beta|=O(r^{\nu+k})$. In particular, observe that $\gamma = -\int_{r}^\infty{\alpha\, ds}$ is well defined and satisfies $d\gamma=\sigma$. Then, for $\chi$ a basic cut-off functions with $\chi\equiv 1$ for $r\geq R+1$, $\sigma-d(\chi\gamma)$ is closed and compactly supported. We then define $\Phi^k_-\co \mathcal{H}^k_\nu(B)\ra H^k_c(B)$ by $\Phi^k_- (\sigma) = [\sigma - d(\chi\gamma)]$.

\begin{enumerate}
\item $\Phi^k_-$ is injective. By Lemma \ref{lem:Excluded:indicial:roots:kforms} (i) if $\sigma\in \mathcal{H}^k_{\nu}(B)$ then $\sigma\in C^\infty_{-n+k+\epsilon}(B)$ for every $\epsilon>0$ and therefore the form $\gamma$ defined earlier by radial integration lies in $C^\infty_{-n+k+1+\epsilon}$. Hence if $\Phi^k_- (\sigma)=0$, \ie $\sigma - d(\chi\gamma)$ is the differential of a basic compactly supported form, then $\sigma = d\gamma'$ with $\gamma'\in C^\infty_{-n+k+1+\epsilon}(B)$. Since $2k < n$, as in Lemma \ref{lem:Harmonic:CClosed} integration by parts is justified and we obtain
\[
\|\sigma \|^2_{L^2}=\langle d\gamma, \sigma\rangle_{L^2}=\langle \gamma, d^\ast\sigma\rangle_{L^2}=0.
\]
\item $\Phi^k_-$ is surjective. If $\sigma$ is a closed basic smooth compactly supported form fix $\alpha\in (0,1)$ and consider the equation $\triangle\gamma = d^\ast\sigma$ for $\gamma\in C^{3,\alpha}_{\nu+1}(B)$. By Proposition \ref{prop:Obstructions:Pu=v} the obstructions to solve this equation lie in the space of basic harmonic $(k-1)$-forms in $C^\infty_{-n-\nu+1}(B)$. By Lemma \ref{lem:Harmonic:CClosed} every such form is closed and therefore all obstructions to solve $\triangle\gamma = d^\ast\sigma$ vanish. Moreover, $d^\ast\gamma$ is a basic harmonic $(k-2)$-form in $C^{2,\alpha}_{\nu}(B)$ and therefore a second application of Lemma \ref{lem:Harmonic:CClosed} implies that $d d^\ast\gamma=0$. We conclude that $\sigma - d\gamma\in\mathcal{H}^k_\nu (B)$ and $\Phi^k_- (\sigma-d\gamma) = \Phi^k_-(\sigma)=[\sigma]$ as desired. 
\end{enumerate}

We will now study $\mathcal{H}^k_\nu (B)$ where $\nu = -k +\delta$ and $k<\tfrac{n}{2}$. By Proposition \ref{prop:Decay:Solutions} and Lemma \ref{lem:Excluded:indicial:roots:kforms} (iv) every $\sigma\in \mathcal{H}^k_\nu (B)$ can be written in the form $\sigma = \tau + \sigma'$ outside a compact set, where $\tau\in\mathcal{H}^k(\Sigma)$ and $\sigma'\in C^{\infty}_{-k-\epsilon}(B)$ for every $\epsilon>0$ sufficiently small. It is also clear that $\tau$ represents the image of $[\sigma]\in H^k(B)$ in $H^k(\Sigma)$. We then define a map $\Phi^k_+\co\mathcal{H}^k_\nu (B)\ra H^k(\Sigma)$ by $\Phi^k_+ (\sigma)=[\tau]$. By basic Hodge theory Proposition \ref{prop:Basic:Hodge} on $N$ the kernel of $\Phi^k_+$ is $\mathcal{H}^k_{-k-\epsilon}$. The image of $\Phi^k_+$ is clearly contained in $\tu{im}\left( H^k(B)\ra H^k(\Sigma)\right)$: we have to prove it coincides with this subspace.

Fix $\tau\in\mathcal{H}^{k}(\Sigma)$ with $[\tau]\in \tu{im}\left( H^k(B)\ra H^k(\Sigma)\right)$. By assumption and Lemma \ref{lem:Asymptotics:closed:forms} there exists a basic closed $k$-form $\sigma$ on $M$ with $\sigma=\tau$ outside a compact set. Moreover, since $M$ is asymptotic to $\BC (N)$ and $\tau$ is basic closed and coclosed on $\BC (N)$, we have $d^\ast\sigma\in C^{0,\alpha}_{-k-1-\epsilon}$ for all sufficiently small $\epsilon>0$. We now study the equation $\triangle\gamma = d^\ast\sigma$ for a basic $(k-1)$-form $\gamma\in C^{3,\alpha}_{-k+1-\epsilon}(B)$. As before, if a solution $\gamma$ exists then $d d^\ast\gamma=0$ and the obstructions to solve the equation lie in the space of basic harmonic $(k-1)$-forms in $C^\infty_{-n+k+1+\epsilon}(B)$, which are all closed by Lemma \ref{lem:Harmonic:CClosed}. Moreover, taking $\epsilon>0$ smaller if necessary, we can assume that $2k<-n-\epsilon$; then $-n+k+1+\epsilon < -k+1$ and therefore every harmonic $(k-1)$-forms $\overline{\gamma}\in C^\infty_{-n+k+1+\epsilon}(B)$ actually lies in $C^\infty_{-n+k-1+\epsilon'}(B)$ for every $\epsilon'>0$. Then we can integrate by parts
\[
\langle d^\ast\sigma, \overline{\gamma}\rangle _{L^2} = \langle \sigma, d\overline{\gamma}\rangle_{L^2}=0
\]
and conclude that all obstructions to solve $\triangle\gamma = d^\ast\sigma$ vanish. It follows that $\sigma - d\gamma \in\mathcal{H}^k_\nu(B)$ and $\Phi^k_+(\sigma-d\gamma)=[\tau]$ as desired.

The case $k>\tfrac{n}{2}$ can be understood by duality. Namely, there is a natural pairing
\[
\mathcal{H}^{k}_{-k+\delta}(B)\times\mathcal{H}^{n-k}_{-n+k-\delta}(B)\ra\R, \qquad (\alpha,\beta)\mapsto \int_M {\theta\wedge\alpha\wedge\beta}.
\]
This pairing is non-degenerate since by Lemma \ref{lem:Excluded:indicial:roots:kforms} (iii) for every $\beta\in\mathcal{H}^{n-k}_{-n+k-\delta}(B)$ we have $\ast\beta\in \mathcal{H}^{k}_{-k+\delta}(B)$. Thus we have an isomorphism $\mathcal{H}^{k}_{-k+\delta}(B)\simeq \mathcal{H}^{n-k}_{-n+k-\delta}(B)^\ast \simeq H^{n-k}_c(B)^\ast\simeq H^k(B)$, which is simply the map that assigns to a basic closed and coclosed form its basic cohomology class. Moreover, using this isomorphism, Proposition \ref{prop:Decay:Solutions} and Lemma \ref{lem:Excluded:indicial:roots:kforms} (iv), it is also clear that $\mathcal{H}^k_{-k+\delta}(B)/\mathcal{H}^k_{-k-\delta}(B) \simeq \tu{im}\left( H^k(B)\ra H^k(\Sigma)\right)$ and therefore $\mathcal{H}^k_{-k-\delta}(B)\simeq \tu{im}\left( H^k_c(B)\ra H^k(B)\right)$ by exactness of \eqref{eq:Exact:Sequence:Boundary}.
\endproof

The case $n=2k$ is more challenging, but parts of Theorem \ref{thm:L2:cohomology} can be proved with very similar arguments to the ones we used.

\begin{lemma}\label{lem:L2:cohomology:n/2}
Let $2k=n$ and fix $\delta>0$ sufficiently small. There is a natural isomorphism
\[
\mathcal{H}^k_{-k+\delta}(B)/\mathcal{H}^k_{-k-\delta}(B)\ra \tu{im}\left( H^k(B)\ra H^k(\Sigma)\right)^{\oplus 2}
\]
and the natural map $\mathcal{H}^k_{-k-\delta}(B)\ra \tu{im}\left( H^k_c(B)\ra H^k(\Sigma)\right)$ that assigns to each closed and coclosed form its cohomology class is injective.
\proof
Set $k=\tfrac{n}{2}$. By Proposition \ref{prop:Decay:Solutions} and Lemma \ref{lem:Excluded:indicial:roots:kforms} (iv) outside a compact set every $\sigma\in \mathcal{H}^k_{-k+\delta}(B)$ can be written as $\sigma = \tau_1 + \ast\tau_2 + \sigma'$ with $\sigma'\in C^\infty_{-k-\delta}(B)$ and $\tau_1,\tau_2\in \mathcal{H}^{k}(\Sigma)$. Define $\Phi^k_+\co \mathcal{H}^k_{-k+\delta}(B)\ra H^k(\Sigma)^{\oplus 2}$ by $\Phi^k_+ (\sigma) = \left( [\tau_1],[\tau_2]\right)$. Now, clearly $[\tau_1]\in \tu{im}\left( H^k(B)\ra H^k(\Sigma)\right)$ and $\Phi^k_+ (\ast\sigma)=\left( \pm [\tau_2],[\tau_1]\right)$ hence $\Phi^k_+$ induces a map from $\mathcal{H}^k_{-k+\delta}(B)$ to $\tu{im}\left( H^k(B)\ra H^k(\Sigma)\right)^{\oplus 2}$ with kernel $\mathcal{H}^k_{-k-\delta}(B)$.

In order to prove that $\Phi^k_+$ is surjective onto $\tu{im}\left( H^k(B)\ra H^k(\Sigma)\right)^{\oplus 2}$, one can show that for every $\sigma\in C^\infty_{-k+\delta}(B)$ it is always possible to solve $\triangle\gamma = d^\ast\sigma$ for $\gamma\in C^\infty_{-k+1+\delta}(B)$ with $d d^\ast\gamma=0$. This is done exactly as in the proof of the surjectivity of $\Phi^k_-$ in the proof of Theorem \ref{thm:L2:cohomology} (i), see Proposition \ref{prop:Laplacian:n/2} below.

In order to prove the second part of the lemma, define $\Phi^k_-\co \mathcal{H}^k_{-k-\delta}(B)\ra H^k(B)$ by $\sigma\mapsto [\sigma]$. The image of $\Phi^k_-$ is clearly contained in the kernel of $H^k(B)\ra H^k(\Sigma)$, which coincides with $\tu{im}\left( H^k_c(B)\ra H^k(\Sigma)\right)$ by exactness of \eqref{eq:Exact:Sequence:Boundary}. The proof of the injectivity of $\Phi^k_-$ is analogous to the one in the proof of Theorem \ref{thm:L2:cohomology} (i).
\endproof
\end{lemma}

It remains to prove that $\Phi^k_-\co \mathcal{H}^k_{-k-\delta}(B)\ra \tu{im}\left( H^k_c(B)\ra H^k(\Sigma)\right)$ is surjective for $k=\tfrac{n}{2}$. This requires a refined analysis of the equation $\triangle\gamma=d^\ast\sigma$ for $\gamma$ a basic $(k-1)$-form of class $C^\infty_{-k+1-\delta}$.

\begin{prop}\label{prop:Laplacian:n/2}
Let $n=2k$ and $\delta>0$ sufficiently small. Let $\sigma$ be a basic smooth $k$-form such that $\sigma\in C^\infty_{-k+\delta}(B)$ and $d^\ast\sigma\in C^\infty_{-k-1-\delta}(B)$. Then the equation $\triangle\gamma=d^\ast\sigma$ has a solution $\gamma\in C^\infty_{-k+1+\delta}(B)$ with $d\gamma\in C^\infty_{-k-\delta}$.
\proof
The obstructions to solve $\triangle\gamma=d^\ast\sigma$ with $\gamma\in C^\infty_{-k+1\mp \delta}(B)$ lie in the space of basic harmonic $(k-1)$-forms of class $C^\infty_{-k+1\pm \delta}(B)$ which are not closed. In particular, we can always solve $\triangle\gamma=d^\ast\sigma$ with $\gamma\in C^\infty_{-k+1+\delta}(B)$ since every basic harmonic $(k-1)$-form of class $C^\infty_{-k+1-\delta}(B)$ is closed by Lemma \ref{lem:Harmonic:CClosed}. Note also that by Lemma \ref{lem:Harmonic:CClosed} any solution satisfies $d d^\ast\gamma=0$ since $d^\ast\gamma$ is a basic harmonic $(k-2)$-form of class $C^\infty_{-k+\delta}(B)$ and we can always assume that $\delta$ is small enough so that $-k+\delta=-\tfrac{n}{2}+\delta<-\tfrac{n}{2}+1$. We need to show that we can take $\gamma$ with $d\gamma\in C^\infty_{-k-\delta}$. The proof of this fact follows the lines of the proof of \cite[Proposition 5.16]{FHN:ALC:G2:from:AC:CY3}. We sketch the key ideas.

The main point is to understand exactly the space of basic harmonic $(k-1)$-forms as we cross the indicial root $-k+1$. Denote by $\triangle^{p}_{\nu}$ the Laplacian restricted to basic $p$-forms of class $C^{l,\alpha}_{\nu}$ for some $l\geq 2$ and $\alpha\in (0,1)$ and by $i(\triangle^{p}_{\nu})$ its index. By Lemma \ref{lem:Excluded:indicial:roots:kforms} (iv) and (v) we have
\[
i(\triangle^{k-1}_{-k+1+\delta}) - i(\triangle^{k-1}_{-k+1-\delta}) = 2\tu{dim}\, H^{k-1}(\Sigma).
\]
Moreover, $\tu{coker}\, \triangle^{k-1}_{-k+1\pm \delta} \simeq \tu{ker}\,\triangle^{k-1}_{-k+1\mp \delta}$ and therefore
\[
\tu{ker}\, \triangle^{k-1}_{-k+1+\delta} \simeq \tu{ker}\, \triangle^{k-1}_{-k+1-\delta} \oplus \mathcal{H}^{k-1}(\Sigma).
\]
Now, decompose $\mathcal{H}^{k-1}(\Sigma)$ into a subspace isomorphic to $\tu{im}\left( H^{k-1}(B)\ra H^{k-1}(\Sigma)\right)$ and a complementary subspace $W$, which is isomorphic to $\tu{im}\left( H^{k}(B)\ra H^{k}(\Sigma)\right)$ via the basic Hodge-star operator on $N$, see \cite[Lemma 5.11]{FHN:ALC:G2:from:AC:CY3}. By Theorem \ref{thm:L2:cohomology} (i) we have $\tu{ker}\, \triangle^{k-1}_{-k+1+\delta}/\mathcal{H}^{k-1}_{-k+1+\delta}(B)\simeq W$.

Note however that every basic harmonic $(k-1)$-form $\overline{\gamma}$ of class $C^\infty_{-k+1+\delta}(B)$ is coclosed. Indeed, outside a compact set we can write
\begin{equation}\label{eq:Harmonic:k-1:not:closed}
\overline{\gamma} = \tau_1 \log{r} + \tau_2 + \overline{\gamma}',
\end{equation}
with $\tau_1,\tau_2\in\mathcal{H}^{k-1}(\Sigma)$ and $\overline{\gamma}'\in C^\infty_{-k+1-\delta}(B)$. Since $\tau_1\log{r} + \tau_2$ is basic coclosed on $\tu{BC}(N)$, we conclude that, up to taking $\delta$ smaller if necessary, $d^\ast\overline{\gamma}$ is a basic harmonic $(k-2)$-form of class $C^\infty_{-k-\delta}(B)$. By Lemma \ref{lem:Harmonic:CClosed}, $d d^\ast\overline{\gamma}=0$ and then an integration by parts shows that
\[
0=\langle d d^\ast\overline{\gamma},\overline{\gamma}\rangle_{L^2} = \| d^\ast\overline{\gamma}\|_{L^2}^2.
\]
Using these facts we conclude that $\ast d\overline{\gamma}$ is a closed and coclosed $k$-form in $C^\infty_{-k+\delta}(B)$. More precisely, using \eqref{eq:Harmonic:k-1:not:closed}, we have
\[
\ast d\overline{\gamma} = \ast_\Sigma \tau_1 + C^\infty_{-k-\delta}(B).
\]
By Lemma \ref{lem:L2:cohomology:n/2} we conclude that $[\tau_1]\in W\subset\mathcal{H}^{k-1}(\Sigma)$. Lemma \ref{lem:L2:cohomology:n/2} also implies that $d\overline{\gamma}=0$ if $\tau_1=0$: indeed, if $\tau_1=0$ then $d\overline{\gamma}\in \mathcal{H}^k_{-k-\delta}(B)$ and $[d\overline{\gamma}]=0$.

Fix an $L^2$--orthonormal basis $\tau_1,\dots, \tau_m$ of $W$. The discussion above implies that for each $j=1,\dots, m$ there exists a basic harmonic form $\overline{\gamma}_j$, unique up to the addition of an appropriately decaying basic closed and coclosed basic form, such that
\[
\overline{\gamma}_j = \tau_j \log{r} + \sum_{i=1}^m{\alpha^i_j \tau_i} +  \overline{\gamma}'
\] 
for some $\alpha^i_j\in\R$ and $\overline{\gamma}'\in C^\infty_{-k+1-\delta}(B)$. The collection $\overline{\gamma}_1, \dots, \overline{\gamma}_m$ forms a basis of the space of obstructions to solve $\triangle\gamma=d^\ast\sigma$ with $\gamma\in C^{\infty}_{-k+1-\delta}(B)$. Now, fix a basic cut-off function $\chi$ with $\chi\equiv 1$ outside of a compact set. An integration by parts shows that
\[
\langle \triangle (\chi \tau_i), \overline{\gamma}_j \rangle_{L^2} = \delta_{ij}.
\]
Thus we can always solve the equation $\triangle\gamma=d^\ast\sigma$ with $\gamma\in C^{\infty}_{-k+1-\delta}(B)$ modulo the span of $\chi\tau_1, \dots, \chi\tau_m$. In order to conclude now note that $d (\chi \tau_i)$ is of class $C^\infty_{-k-\delta}(B)$ for $\delta$ sufficiently small since $\tau_i$ is closed on $\BC (\Sigma)$.
\endproof 
\end{prop}

\proof[Proof of Theorem \ref{thm:L2:cohomology} (ii)]
In view of Lemma \ref{lem:L2:cohomology:n/2} it remains only to prove that $\Phi^k_-\co \mathcal{H}^k_{-k-\delta}(B)\ra \tu{im}\left( H^k_c(B)\ra H^k(\Sigma)\right)$ is surjective for $k=\tfrac{n}{2}$. Consider then a closed compactly supported basic $k$-form $\sigma$. By Proposition \ref{prop:Laplacian:n/2} we can solve $\triangle\gamma=d^\ast\sigma$ with $d\gamma\in C^\infty_{-k-\delta}(B)$ and $d d^\ast\gamma=0$. Then $\sigma - d\gamma\in\mathcal{H}^k_{-k-\delta}(B)$ and $\Phi^k_-(\sigma-d\gamma)=\Phi^k_-(\sigma)=[\sigma]\in H^k(B)$. 
\endproof

\section{Highly collapsed $\spins$--metrics on principal Seifert circle bundles}\label{sec:Adiabatic}

In this section we prove our main existence result Theorem \ref{thm:Main}: the existence of highly collapsed $\spins$--metrics on the total space of a principal Seifert circle bundle over an AC $\gtwo$--orbifold. In Section \ref{sec:AC:G2} we establish properties of AC $\gtwo$--orbifolds we use in Section \ref{sec:Adiabatic:proof} to prove Theorem \ref{thm:Main}.

\subsection{Asymptotically conical $\gtwo$--orbifolds}\label{sec:AC:G2}

In this section we study $8$-dimensional transversally AC principal Seifert circle bundles $\pi\co M\ra B$ carrying a transverse $\gtwo$--structure which is parallel with respect to the adapted connection. We introduce the necessary notation and prove results about harmonic and closed and coclosed basic forms on such manifolds.

\subsubsection{Basic torsion-free $\gtwo$--structures}

Let $\pi\co M^8\ra B$ be a principal Seifert circle bundle endowed with a connection $1$-form $\theta$. A \emph{basic $\gtwo$--structure} is a reduction of the structure group of the horizontal subbundle $\mathcal{H}=\ker\theta$ to $\gtwo$. Equivalently, a basic $\gtwo$--structure is the choice of a basic $3$-form $\varphi\in\Omega^3(B)$ that is point-wise equivalent, under an appropriate identification of the fibres of $\mathcal{H}$ with $\R^7$, to the standard flat $\gtwo$--structure $(u,v,w)\mapsto \langle uv, w\rangle$ on $\R^7=\Imag\mathbb{O}$. Every basic $\gtwo$--structure determines a horizontal Riemannian metric $g_B=g_\varphi$ on $M$ and $(M,\theta, g_B)$ is then a Riemannian principal Seifert circle bundle. We will denote by $\psi$ the basic $4$-form $\psi = \ast\varphi$. The triple $(M,\pi,\theta,\varphi)$ will be called a \emph{$\gtwo$ principal Seifert circle bundle}.

We will now collect important identities for basic forms on $\gtwo$ principal Seifert circle bundles. We refer the reader to \cite{Bryant:Rmks:G2} for their proof.

\begin{lemma}\label{lem:decomp:forms}
Let $(M,\pi,\theta,\varphi)$ be a $\gtwo$ principal Seifert circle bundle. Then there are point-wise orthogonal decompositions
\[
\Omega^2(B) = \Omega^2_7(B) \oplus \Omega^2_{14}(B), \qquad \Omega^3(B) = \Omega^3_1(B)\oplus \Omega^3_7(B) \oplus \Omega^2_{14}(B)
\]
of basic forms according to irreducible representations of $\gtwo$, where: 
\[
\begin{gathered}
\Omega^2_7(B) = \{ X\lrcorner\varphi=\ast(X^\flat\wedge\psi)\, |\, X^\flat\in\Omega^1(B)\}, \qquad \Omega^3_7(B) = \{ X\lrcorner\psi=-\ast(X^\flat\wedge\varphi)\, |\, X^\flat\in\Omega^1(B)\},\\
\Omega^2_{14}(B)=\{ \tau\in\Omega^2(B)\, |\, \tau\wedge\psi=0\} = \{ \sigma\in\Omega^2(B)\, |\, \ast\sigma=-\sigma\wedge\varphi\},\\
\Omega^3_1(B)=\{ f\varphi\, |\, f\in\Omega^0(B)\}, \qquad 
\Omega^3_{27}(B) = \{ \rho\in\Omega^3(B)\, |\, \rho\wedge\varphi=0=\rho\wedge\psi\}.
\end{gathered}
\]
\end{lemma}
By acting with the basic Hodge-star operator $\ast$, Lemma \ref{lem:decomp:forms} implies similar decompositions of the space of basic $4$-forms and $5$-forms. The decomposition of $4$-forms can be used to describe the linearisation of the map $\varphi\mapsto \psi$, see \cite[Proposition 4]{Bryant:Rmks:G2}.

\begin{lemma}\label{lem:Linearisation:Hitchin:dual:3form:7d}
The linearisation of the map $\varphi\mapsto \psi$ is
\[
\rho \longmapsto \ast\left( \tfrac{4}{3}\pi_1\rho + \pi_7 \rho - \pi_{27}\rho\right).
\]
\end{lemma} 

\begin{lemma}\label{lem:Torsion:G2:structures}
Let $(M,\pi,\theta,\varphi)$ be a $\gtwo$ principal Seifert circle bundle. Then there exists a basic function $\tau_0$, a basic $1$-form $\tau_1$, $\tau_2\in \Omega^2_{14}(B)$ and $\tau_3\in\Omega^3_{27}(B)$ such that
\begin{align}
&d\varphi = \tau_0\psi + 2\tau_1\wedge\varphi + \tau_3,\\
&d\psi = 4\tau_1\wedge\psi + \tau_2\wedge\varphi.
\end{align}
\end{lemma}
For a proof, see \cite[Proposition 1]{Bryant:Rmks:G2}. If $d\varphi=0=d\psi$, \ie $\tau_i=0$ for all $i$, we say that $\varphi$ is \emph{torsion-free}. Linear algebra and representation theory of $\gtwo$ then imply that $\varphi$ is parallel for the adapted connection. Since $\gtwo$ is the stabiliser in $\tu{GL}(7,\R)$ of the flat $\gtwo$--structure on $\R^7$ one then concludes that the existence of a torsion-free transverse $\gtwo$--structure implies that the holonomy of the adapted connection reduces from $\sorth{7}$ to $\gtwo$. We also note that a crucial consequence of the torsion-free condition is the vanishing of the transverse Ricci-curvature, see \cite[\S 4.5.3]{Bryant:Rmks:G2}. From now on we will concentrate on torsion-free basic $\gtwo$--structures. In this case we will refer to $(M,\pi,\theta,\varphi)$ as a \emph{$\gtwo$--holonomy} principal Seifert circle bundle.

We collect some useful identities for basic functions and $1$-forms on a $\gtwo$--holonomy principal Seifert circle bundle. In the following lemma the $\tu{curl}$ operator acting on basic $1$-forms is defined by
\begin{equation}\label{eq:curl}
\tu{curl}\,\gamma =\ast(d\gamma\wedge \psi).
\end{equation}

\begin{lemma}\label{lem:identities:0:1:forms:torsion:free}
For $f\in\Omega^0(B)$ and $\gamma\in\Omega^1(B)$ with $X=\gamma^\sharp$, the following identities hold:
\begin{enumerate}
\item\label{lem:identities:0:1:forms:diff} $\pi_1 d (X\lrcorner\varphi) = -\tfrac{3}{7}(d^\ast\gamma)\varphi$ and $\pi_7d (X\lrcorner\varphi) = \tfrac{1}{2}\ast\big( \tu{curl}\,\gamma \wedge \varphi\big)$.
\item\label{lem:identities:0:1:forms:7:component} $d \ast(\gamma\wedge\psi) - d^\ast (\gamma\wedge\varphi) = \ast(\tu{curl}\gamma\wedge\varphi )- (d^\ast\gamma)\varphi\in\Omega^3_{1\oplus 7}(B)$.
\item\label{lem:identities:0:1:forms:Dirac} The basic Dirac operator $\slashed{D}$ can be identified with either of the following two operators:
\[
\begin{gathered}
\slashed{D}\co \Omega^0(B)\oplus \Omega^1(B)\ra \Omega^0(B)\oplus \Omega^1(B), \qquad (f,\gamma)\mapsto (d^\ast\gamma, d f + \tu{curl}\, \gamma),\\
\slashed{D}\co \Omega^0(B)\oplus \Omega^1(B)\ra \Omega^3_{1\oplus 7}(B), \qquad (f,\gamma)\mapsto \pi_{1\oplus 7}\Big(\ast d (f\varphi) + d \ast(\gamma\wedge\psi)\Big).
\end{gathered}
\]
\end{enumerate}
\proof
The lemma can be deduced from \cite[Proposition 3]{Bryant:Rmks:G2} and the identification of the Dirac operator of a $\gtwo$--manifold $B$ via the isomorphism between the spinor bundle with $\R\oplus TB$, see for example \cite[Equation (6.2)]{Nordstrom:Dirac}. 
\endproof
\end{lemma}
 
\subsubsection{Nearly K\"ahler orbifolds and transversally AC $\gtwo$--holonomy Seifert bundles}

Let $N$ be a closed oriented $7$-manifold. Assume that $\pi_\infty\co N\ra \Sigma$ is a principal Seifert circle bundle over a closed $6$-orbifold $\Sigma$ and fix a connection $\theta_\infty$ on $\pi_\infty$. We say that $N$ admits a basic $\sunitary{3}$--structure if there exist basic forms $\omega\in\Omega^2(\Sigma)$ and $\Real\Omega,\Imag\Omega\in\Omega^3 (\Sigma)$ such that $\omega$ and $\Omega$ are, respectively, a non-degenerate $2$-form and a complex volume form on the horizontal subspace $\mathcal{H}$ satisfying $\omega\wedge\Omega=0$ and $2\omega^3 = 3\Real\Omega\wedge\Imag\Omega$. Every basic $\sunitary{3}$--structure defines a horizontal metric $g_\Sigma$ and, together with $\theta_\infty$, the structure of a Riemannian Seifert bundle on $N$. A basic $\sunitary{3}$--structure $(\omega,\Omega)$ is called \emph{nearly K\"ahler} if
\begin{equation}\label{eq:nK}
d\omega = 3\Real\Omega, \qquad d\Imag\Omega = -2\omega^2.
\end{equation}
If $(\omega,\Omega)$ is a basic nearly K\"ahler structure on $N$ then $\BC (N)$ has a basic torsion-free $\gtwo$--structure
\begin{equation}\label{eq:G2:cone}
\varphi_{\tu{C}} = r^2 dr\wedge\omega + r^3\Real\Omega.
\end{equation}
In particular, since the metric on $\BC (N)$ has vanishing transverse Ricci-curvature, every basic nearly K\"ahler structure induces a transversally Einstein metric $g_\Sigma$ with Einstein constant $5$.

Since the horizontal metric induced by \eqref{eq:G2:cone} is of the form $g_{\tu{C}}=dr^2 + r^2 g_\Sigma$, it makes sense to talk of a transversally AC $\gtwo$--holonomy principal Seifert circle bundle $(M,\pi,\theta,\varphi)$ asymptotic to $\left( \BC (N), \pi_\infty,\theta_\infty,\varphi_{\tu{C}}\right)$. Since $\varphi$ and $\varphi_{\tu{C}}$ determine the horizontal metrics $g_B$ and $g_{\tu{C}}=dr^2 + r^2 g_\Sigma$, we can replace \eqref{eq:Transv:AC} in Definition \ref{def:Transv:AC} with a similar polynomial decay for $\theta$ and $\varphi$ (and their covariant derivatives) to $\theta_\infty$ and $\varphi_{\tu{C}}$.

In the discussion so far the choice of the connection $1$-form $\theta_\infty$ on $\pi_\infty\co N\ra \Sigma$ (and, by radial extension, on $\BC (N)$) remained free. In fact, up to gauge transformations there is a canonical choice. Indeed, if $(\omega,\Omega)$ is a basic nearly K\"ahler structure on $N$ then it follows from \cite[Theorem 3.8]{Foscolo:nK:def} that every basic closed and coclosed $2$-form $\kappa$ on $N$ satisfies
\begin{equation}\label{eq:Primitive:(1,1)}
\kappa\wedge\omega^2 =0=\kappa\wedge\Omega,
\end{equation}
\ie $\kappa$ is a basic \emph{primitive $(1,1)$--form}. Conversely, every basic closed primitive $(1,1)$--form $\kappa$ is necessarily coclosed since $\ast\kappa = -\kappa\wedge\omega$. In particular, by adding a basic $1$-form to $\theta_\infty$ so that $d\theta_\infty$ is closed and coclosed, we can always assume that $d\theta_\infty$ satisfies \eqref{eq:Primitive:(1,1)} and this requirement uniquely determines $\theta_\infty$ up to gauge transformations. In the rest of the section we say that $(N,\pi_\infty,\theta_\infty,\omega,\Omega)$ is a nearly K\"ahler principal Seifert circle bundle if $(\omega,\Omega)$ is a basic nearly K\"ahler structure and $\theta_\infty$ is normalised so that $d\theta_\infty$ is a primitive $(1,1)$--form.

In the rest of the paper we will also work under the following assumption.
\begin{assumption}
The principal Seifert circle bundle $\pi_\infty\co N\ra \Sigma$ is not flat.
\end{assumption}
Since we assume that $d\theta_\infty$ is the unique basic closed and coclosed representative of its cohomology class, the assumption is equivalent to the requirement that the image of $c^{orb}_1(N)$ in $H^2_{orb}(\Sigma;\R)$ does not vanish. Note that if $(M,\pi,\theta,\varphi)$ is asymptotic to $\tu{BC}(N)$, then necessarily $d\theta\neq 0$ also.

\begin{prop}\label{prop:fund:gp}
Let $(M,\pi,\theta,\varphi)$ be a transversally AC $\gtwo$--holonomy principal Seifert circle bundle asymptotic to $\tu{BC} (N)$ and assume that $d\theta_\infty\neq 0$. Then
\begin{enumerate}
\item $M$ has finite fundamental group;
\item $M$  is an \emph{irreducible} Seifert circle bundle, \ie there are no basic $1$-forms on $M$ that are parallel with respect to the adapted connection.
\end{enumerate}
\proof
The proof of part (i) is analogous to the one of \cite[Proposition 5.9]{FHN:ALC:G2:from:AC:CY3}. Observe that, since the horizontal metric $g_\Sigma$ is Einstein with uniform positive Einstein constant, the orbifold $\Sigma$ has finite orbifold fundamental group by Remark \ref{rmk:Myers}. Since $d\theta_\infty\neq 0$, the homotopy exact sequence \eqref{eq:Homotopy} then presents $\pi_1(N)$ as an extension of $\pi_1^{orb}(\Sigma)$ by a finite group. We therefore conclude that $N$ has finite fundamental group.

Now, the map $\pi_1 (M\setminus K)\ra \pi_1 (M)$ is surjective, since otherwise, as in \cite[Lemma 2.18]{ACyl:CY}, a finite cover of $M$ would be a Riemannian submersion over a Ricci-flat orbifold with at least two asymptotically conical ends: this is impossible by the orbifold version of Cheeger--Gromoll's Splitting Theorem \cite{Splitting:Orbifolds} or, considering a sequence of metrics on $M$ that collapse to the orbifold $B$, Cheeger--Colding's Almost Splitting Theorem \cite[Theorem 6.64]{Almost:Splitting}.

For part (ii) it is enough to observe that there are no parallel basic $1$-forms on $\tu{BC}(N)$. In order to show that this is true, observe that there is a correspondence between basic parallel $1$-forms on $\BC (N)$ and basic functions $f$ on $N$ satisfying $\nabla df=f g_\Sigma$ (this follows from an explicit calculation of the Levi-Civita connection of a Riemannian cone). In particular, any such $f$ would satisfy $\triangle f=6f$ and by Remark \ref{rmk:Obata} $\Sigma$ would then be isometric to a finite quotient of the round $6$-sphere. In this case however $\pi_\infty\co N\ra\Sigma$ would be forced to be flat.
\endproof
\end{prop}

\begin{remark*}
In view of part (i), the irreducibility assumption of part (ii) is equivalent to the requirement that the holonomy of the adapted connection is the whole group {\gtwo} by \cite[Lemma 1]{Bryant:1987}.
\end{remark*} 

Although we have outlined some of its important consequences, at this stage the assumption $d\theta_\infty\neq 0$ may seem unmotivated. In fact, it is a necessary condition for obtaining non-trivial examples from Theorem \ref{thm:Main}, see Remark \ref{rmk:Assumption} below. 

\subsubsection{Basic closed and coclosed forms on AC {\gtwo}--orbifolds}

Let $(N,\pi_\infty,\theta_\infty,\omega,\Omega)$ be a nearly K\"ahler principal Seifert circle bundle and $(M,\pi,\theta,\varphi)$ be a transversally AC $\gtwo$--holonomy principal Seifert circle bundle asymptotic to $\BC (N)$. In this section we prove a number of results about basic harmonic and closed and coclosed forms on $M$. The main results we are going to use in the next section are Proposition \ref{prop:Laplacian:2:forms} and Corollary \ref{cor:Exact:5:forms} which describe the obstructions to solve the Poisson equation $\triangle\sigma=\kappa$ for a basic $2$-form $\kappa$ and provide a normal form for exact basic $5$-forms with appropriate decay.

We first analyse basic harmonic functions and $1$-forms on $M$. Using the fact that every nearly K\"ahler principal Seifert circle bundle is transversally Einstein with Einstein constant $5$, Proposition \ref{prop:Eigenvalue:Estimates} yields an improvement of the indicial roots computations of Lemma \ref{lem:Excluded:indicial:roots:kforms}. 

\begin{lemma}\label{lem:Excluded:indicial:roots:0:1:forms} 
Let $N$ be a transversally nearly K\"ahler principal Seifert circle bundle.
\begin{enumerate}
\item There are no basic harmonic functions on $\BC(N)$ homogeneous of order $\lambda\in [-6,1]$ except for constant multiples of $1$ and $r^{-5}$ (of rate $\lambda=0$ and $\lambda=-5$ respectively).
\item There are no basic harmonic $1$-forms on $\BC(N)$ homogeneous of order $\lambda\in (-5,0)$.
\end{enumerate}
\end{lemma}

The irreducibility in Proposition \ref{prop:fund:gp} (ii) and the fact that $g$ has vanishing transverse Ricci-curvature have the following important consequence.

\begin{lemma}\label{lem:Vanishing:Harmonic:0:1:Forms}
There are no basic harmonic functions and $1$-forms on $M$ in $C^\infty_\nu (B)$ for $\nu< 0$.
\proof
Let $u$ and $\gamma$ be a basic harmonic function and $1$-form in $C^\infty_{\nu}(B)$. If $\nu<-\frac{5}{2}$ then we can integrate by parts:
\[
0=\langle \triangle u , u\rangle _{L^2} = \|\nabla u \|^2_{L^2}, \qquad 0=\langle \triangle \gamma, \gamma\rangle_{L^2} = \|\nabla\gamma\|^2_{L^2},
\]
where, since $M$ is transversally Ricci-flat, we used the fact that $\triangle = \nabla^\ast\nabla$ on basic $1$-forms. We conclude that $u$ is constant and therefore vanishes since it decays at infinity and $\gamma=0$ since $M$ is irreducible. On the other hand, by Lemma \ref{lem:Excluded:indicial:roots:0:1:forms} there are no indicial roots for the Laplacian acting on functions and $1$-forms in the interval $(-5,0)$.
\endproof 
\end{lemma}

\begin{remark}\label{rmk:Vanishing:Harmonic:0:1:Forms}
In particular, if $\sigma\in C^\infty_\nu(B)$ is a basic harmonic $2$-form and $\rho\in C^\infty_\nu (B)$ is a basic harmonic $3$-form on $M$ for $\nu<0$ then $\sigma\in\Omega^2_{14}(B)$ and $\rho\in\Omega^3_{27}(B)$. Indeed, the basic Laplacian preserves the decomposition of Lemma \ref{lem:decomp:forms} and its restriction to $\Omega^\bullet_1(B)$ and $\Omega^\bullet_7(B)$ coincides with the basic Laplacian on functions and $1$-forms, respectively.
\end{remark}

In view of this remark, we can always normalise the choice of $\theta$ on $M$ so that $d\theta\in\Omega^2_{14}(B)$. Indeed, by the exact sequence \eqref{eq:Exact:Sequence:Boundary} and Theorem \ref{thm:L2:cohomology} the basic cohomology class of $d\theta$ is represented by a unique basic closed and coclosed form decaying with rate $-2$. Therefore we can add a decaying basic $1$-form to $\theta$ so that $d\theta$ is closed and coclosed, thus of type $14$ by Remark \ref{rmk:Vanishing:Harmonic:0:1:Forms}, and $\theta$ is still asymptotic to $\theta_\infty$. In the rest of the section we will include this assumption in the definition of a transversally AC $\gtwo$--holonomy Seifert bundle $(M,\pi,\theta,\varphi)$. 

\begin{prop}\label{prop:Iso:Dirac}
The basic Dirac operator $\slashed{D}\co C^{1,\alpha}_\nu\Omega^0(B)\oplus C^{1,\alpha}_\nu\Omega^1(B)\ra C^{0,\alpha}_{\nu-1}\Omega^0(B)\oplus C^{0,\alpha}_{\nu-1}\Omega^1 (B)$ is surjective if $\nu>-6$ and injective if $\nu<0$.
\proof
Using the vanishing of the transverse Ricci-curvature we conclude that every basic $(f,\gamma)$ in the kernel of $\slashed{D}$ is actually harmonic. Then Lemma \ref{lem:Vanishing:Harmonic:0:1:Forms} implies that $\slashed{D}$ is injective if $\nu<0$ and, by duality, surjective if $\nu>-6$.
\endproof
\end{prop}

\begin{prop}\label{prop:Laplacian:2:forms}
Fix $\alpha\in (0,1)$, $k\geq 1$ and a generic $\nu\in (-4,0)$. Then there exists $C>0$ with the following significance. Let $\kappa$ be a basic $2$-form of class $C^{k-1,\alpha}_{\nu-1}(B)$. If $\langle \kappa, \overline{\kappa}\rangle_{L^2}=0$ for all $\overline{\kappa}\in \mathcal{H}^2_{-6-\nu}(B)$ then there exists a unique $\sigma\in C^{k+1,\alpha}_{\nu+1}\Omega^2 (B)$ which is $L^2_{\nu+1}$--orthogonal to basic harmonic $2$-forms of class $C^\infty_{\nu+1}(B)$ and satisfies $\triangle\sigma=\kappa$ and
\[
\| \sigma \|_{C^{k+1,\alpha}_{\nu+1}(B)}\leq C \| \kappa\|_{C^{k-1,\alpha}_{\nu-1}(B)}.
\]
Moreover if $d^\ast\kappa=0$ then $d^\ast\sigma=0$ and $\kappa = d^\ast d\sigma$ and if $d\kappa=0$ then $d\sigma\in\Omega^3_{27}(B)$.
\proof
Since $2<\tfrac{7}{2}-1$ and $\nu>-4$, Lemma \ref{lem:Harmonic:CClosed} shows that every basic harmonic $2$-form of class $C^\infty_{-6-\nu}(B)$ is closed and coclosed. The first part of the Proposition follows immediately. Since $\nu<0$, the last two statements follow immediately from Lemma \ref{lem:Vanishing:Harmonic:0:1:Forms} and Remark \ref{rmk:Vanishing:Harmonic:0:1:Forms} since $d^\ast\sigma$ and $d\sigma$ are then a basic (weakly) harmonic $1$-form and $3$-form, respectively, of class $C^k_\nu(B)$.
\endproof
\end{prop}

We can now combine the two previous propositions to deduce a normal form for basic exact $5$-forms.

\begin{corollary}\label{cor:Exact:5:forms}
Fix $\alpha\in (0,1)$, $k\geq 1$ and a generic $\nu\in (-4,0)$. Then for every basic closed $5$-form $\tau$ of class $C^{k-1,\alpha}_{\nu-1}(B)$ which is $L^2(B)$--orthogonal to $\mathcal{H}^5_{-6-\nu}(B)$, there exist unique $\gamma\in\Omega^1(B)$ and $\sigma\in\Omega^2(B)$ of class $C^{k+1,\alpha}_{\nu+1}(B)$ such that $d\sigma\in\Omega^3_{27}(B)$ and
\[
\tau = d \big( \tu{curl}\gamma\wedge\varphi - (d^\ast\gamma)\psi - \ast d\sigma\big).
\]
Moreover, there exists a constant $C>0$ independent of $\tau, \gamma,\sigma$ such that
\[
\| \gamma \| _{C^{k+1,\alpha}_{\nu+1}(B)} + \|\sigma\|_{C^{k+1,\alpha}_{\nu+1}(B)}\leq C \| \tau\|_{C^{k-1,\alpha}_{\nu-1}(B)}
\]
\proof
We apply Proposition \ref{prop:Laplacian:2:forms} with $\kappa=-\ast\tau$. The assumptions on $\tau$ guarantee that (a) all obstructions to solve $\triangle\sigma'=\kappa$ vanish, and (b) $\kappa = d^\ast d\sigma'$, \ie $\tau = d\ast d\sigma'$.

Using Proposition \ref{prop:Iso:Dirac} and Lemma \ref{lem:identities:0:1:forms:torsion:free} \ref{lem:identities:0:1:forms:Dirac} we now write
\[
d\sigma' = \ast d (f\varphi) + d \ast(\gamma\wedge\psi) - \rho,
\]
for $f\in\Omega^0 (B)$, $\gamma\in\Omega^1(B)$ and $\rho\in\Omega^3_{27}(B)$ of class $C^{k+1,\alpha}_{\nu+1}(B)$ and $C^{k,\alpha}_{\nu}(B)$, respectively. Indeed, since $\nu+1>-3>-6$, the Dirac operator $\slashed{D}$ is surjective. Since $\nu+1<1$, Lemma \ref{lem:Excluded:indicial:roots:0:1:forms} (i) implies that the kernel of $\slashed{D}$ consists only of constant functions, and therefore we can make our choice unique by normalising $f$ by requiring that it decays at infinity. With this normalisation, we conclude that $f=0$. Indeed, since $\rho\wedge\varphi_0=0$ the vanishing of $d (d\sigma'\wedge\varphi_0)$ implies that $f$ is harmonic. In particular, $\rho = d\sigma$ with $\sigma= \ast(\gamma\wedge\psi)- \sigma' \in C^{k+1,\alpha}_{\nu+1}(B)$.

In order to conclude the proof we now use Lemma \ref{lem:identities:0:1:forms:torsion:free} \ref{lem:identities:0:1:forms:7:component} to rewrite
\[
\tau = d\ast \Big(  d\ast(\gamma\wedge\psi)-d^\ast (\gamma\wedge\varphi) - d\sigma \Big) = d \big( \tu{curl}\gamma\wedge\varphi - (d^\ast\gamma)\psi - \ast d\sigma\big).\qedhere
\] 
\end{corollary}

\begin{remark}\label{rmk:Exact:5:forms}
An integration by parts shows that $\tau$ is $L^2(B)$--orthogonal to $\mathcal{H}^5_{-6-\nu}(B)$ whenever $\tau=du$ with $u\in C^{1,\alpha}_{\nu}(B)$.
\end{remark}

\subsection{Adiabatic limit of circle invariant $\spins$--metrics}\label{sec:Adiabatic:proof}

In \cite{FHN:ALC:G2:from:AC:CY3} we developed a construction of complete $\gtwo$--holonomy metrics on appropriate principal circle bundles over AC Calabi--Yau $3$-folds. In \cite{Apostolov:Salamon} Apostolov--Salamon described the dimensional reduction of the nonlinear PDEs for $\gtwo$--holonomy in the presence of a Killing field. The resulting equations, called the Apostolov--Salamon equations in \cite{FHN:ALC:G2:from:AC:CY3}, form a complicated system of nonlinear PDEs. The strategy of \cite{FHN:ALC:G2:from:AC:CY3} is to construct solutions of the Apostolov--Salamon equations by studying the adiabatic limit of the equations as the orbits of the Killing field shrink to zero size. In this section we describe a similar construction of complete $\spins$--holonomy metrics on principal Seifert circle bundles over AC $\gtwo$--orbifolds.

We consider a $\spins$--manifold admitting an isometric circle action. We interpret the dimensional reduction of the equations for $\spins$--holonomy in terms of the intrinsic torsion of the $\gtwo$--structure induced on the $7$-dimensional orbifold quotient and a coupled abelian $\gtwo$--monopole. These equations can be thought of as a $\spins$ analogue of the Gibbons--Hawking construction of $4$-dimensional hyperk\"ahler metrics with a triholomorphic circle action. As for the Apostolov--Salalmon equations, however, the dimensional reduction of the $\spins$--holonomy equations to $7$ dimensions still consists of \emph{nonlinear} equations and in general it is not clear how to solve them directly. We therefore consider the \emph{adiabatic limit} of these equations when the circle orbits have uniformly small length. This formal picture is then turned into our existence result Theorem \ref{thm:Main} using the analysis of the previous sections. With respect to \cite{FHN:ALC:G2:from:AC:CY3} we are able to give a more streamlined argument by applying the Implicit Function Theorem directly instead of proving convergence of a formal power series solution. 

\subsubsection{Gibbons--Hawking-type ansatz for {\spins}--manifolds}

Let $\pi\co M^8\ra B$ be a principal Seifert circle bundle over a $7$--orbifold $B$. Denote by $\xi$ the vector field that generates the fibre-wise circle action, normalised to have period $2\pi$. A $\spins$--structure on $M$ is the choice of a $4$-form $\Phi$ with distinguished algebraic properties at each point (an admissible $4$-form in the language of \cite[Definition 10.5.1]{Joyce:Book}). Circle-invariant $\spins$--structures are completely determined by the choice of a connection $1$-form $\theta$, a basic positive function $h$ and a basic $\gtwo$--structure $\varphi$:
\begin{equation}\label{eq:circle:invariant:Spin(7)}
\Phi = \theta \wedge \varphi + h^{\frac{2}{3}}\psi,
\end{equation}
where $\psi=\ast\varphi$. The metric $g_\Phi$ induced by $\Phi$ on $M$ is
\begin{equation}\label{eq:circle:invariant:Spin(7):metric}
g_\Phi = h^{\frac{1}{3}}\, g_B + h^{-1}\theta^2,
\end{equation}
where $g_B$ is the horizontal metric induced by the basic $\gtwo$--structure $\varphi$.

By \cite[Definition 10.5.2]{Joyce:Book} a $\spins$--structure $\Phi$ is torsion-free, \ie the metric $g_\Phi$ has holonomy contained in $\spins$, if and only if $d\Phi=0$. We now express the torsion-free condition $d\Phi=0$ for an $S^1$--invariant $\spins$--structure $\Phi$ \eqref{eq:circle:invariant:Spin(7)} as a PDE system for the triple $(\varphi,h,\theta)$.

\begin{prop}
The $S^1$--invariant $\spins$--structure $\Phi$ on $M$ determined by the triple $(\varphi,h,\theta)$ via \eqref{eq:circle:invariant:Spin(7)} is torsion-free if and only if
\begin{equation}\label{eq:GH:Spin(7)}
d\varphi =0, \qquad d\left( h^\frac{2}{3}\psi\right) +d\theta\wedge\varphi=0.
\end{equation}
\end{prop}
The aim of this section is to construct solutions to  \eqref{eq:GH:Spin(7)}. The equations are nonlinear, so it is unclear how to find solutions in general. We make however two easy remarks.

First of all, there is a cohomological constraint in order to be able to solve \eqref{eq:GH:Spin(7)}. Indeed, if a solution exists then $[\varphi]$ is a well-defined basic cohomology class and $[d\theta]\cup [\varphi]=0$ in basic cohomology. Note also that $[d\theta]$ represents the (real) orbifold first Chern class of the orbibundle $\pi\co M\ra B$ in the orbifold cohomology $H^2_{orb}(B;\R)$ of Proposition \ref{prop:Basic:Singular:Coho}, so we can rewrite the constraint as
\begin{equation}\label{eq:top:obstr}
c^{orb}_1(M)\cup [\varphi] = 0 \in H^5_{orb}(B). 
\end{equation}
Secondly, we can use Lemma \ref{lem:Torsion:G2:structures} to reinterpret \eqref{eq:GH:Spin(7)} as coupled equations for the torsion of the basic $\gtwo$--structure $\varphi$ and for the pair $(h,\theta)$.

\begin{lemma}\label{lem:Torsion:GH:Spin(7)}
The basic $\gtwo$--structure $\varphi$ arising from a solution $(\varphi,h,\theta)$ of \eqref{eq:GH:Spin(7)} has torsion
\[
\tau_0=\tau_1=\tau_3=0, \qquad \tau_2 = -h^{-\frac{2}{3}}\kappa_0,
\]
where $\kappa_0$ is the component of the curvature $d\theta$ of $\theta$ in $\Omega^2_{14}(B)$. Moreover, $\left( h,\theta\right)$ satisfies
\begin{equation}\label{eq:G2:monopole:eq}
\ast d\left( \tfrac{3}{2}h^{\frac{2}{3}}\right) +d\theta\wedge\psi =0.
\end{equation}
\proof
Decompose $d\theta = U\lrcorner\varphi + \kappa_0$ for a basic vector field $U$ and $\kappa_0\in\Omega^2_{12}(B)$. By \cite[Equations (3.2) and (3.5)]{Bryant:Rmks:G2}, $d\theta\wedge\varphi=2U^\flat\wedge\psi + \kappa_0\wedge\varphi$. Then it is straightforward to check that \eqref{eq:GH:Spin(7)} is equivalent to the stated expressions for the torsion conponents. Moreover, the vanishing of $\tau_1$ forces $U^\flat=\frac{1}{3}h^{-\frac{1}{3}}dh$. By \cite[Equations (3.4) and (3.5)]{Bryant:Rmks:G2} and the definition of $\Omega^2_{14}$ in Lemma \ref{lem:decomp:forms}, $\ast\left( d\theta\wedge\psi\right) = 3 U^\flat$ and therefore we arrive at \eqref{eq:G2:monopole:eq}.
\endproof
\end{lemma}
The equation \eqref{eq:G2:monopole:eq} is a known gauge-theoretic equation that arises as the dimensional reduction of the (abelian) $\spins$--instanton equations to dimension $7$; its solutions are called abelian $\gtwo$ (or octonionic) monopoles.

\subsubsection{Formal analysis of the adiabatic limit equations}

Lacking a better understanding of \eqref{eq:GH:Spin(7)}, our strategy to produce solutions to \eqref{eq:GH:Spin(7)} is to degenerate the equations by introducing a small parameter $\epsilon>0$ and study solutions in the limit $\epsilon\ra 0$: we assume the existence of a solution to the formal limit of the equations when $\epsilon=0$ and prove that it can be perturbed to a solution of the system for small $\epsilon>0$. The particular degeneration we introduce is geometrically very natural: we consider a $1$-parameter family $\{ \Phi_\epsilon\} _{\epsilon>0}$ of $S^1$--invariant torsion-free $\spins$--structures on $M$ such that the circle orbits shrink to zero length as $\epsilon \ra 0$. By rescaling along the circle orbits we write
\[
\Phi_\epsilon = \epsilon\,\theta\wedge\varphi + h^\frac{2}{3}\psi, \qquad g_{\Phi_\epsilon} = h^{\frac{1}{3}}\, g_B + \epsilon^2 h^{-1}\theta^2,
\]
where $g_B$ is the horizontal metric induced by the basic $\gtwo$--structure $\varphi$. The PDE system \eqref{eq:GH:Spin(7)} for $\Phi_\epsilon$ then becomes
\begin{equation}\label{eq:GH:Spin(7):collapsing:sequence}
d\varphi =0, \qquad d\left( h^{\frac{2}{3}}\psi\right)+ \epsilon\, d\theta\wedge\varphi=0.
\end{equation}
For $\epsilon>0$ the system \eqref{eq:GH:Spin(7):collapsing:sequence} is equivalent to \eqref{eq:GH:Spin(7)}. In the limit $\epsilon\ra 0$, however, the equations simplify: Lemma \ref{lem:Torsion:G2:structures} implies that solutions to \eqref{eq:GH:Spin(7):collapsing:sequence} satisfy $dh_0=0$ and $d\varphi_0=0=d\psi_0$. Assume then the existence of a basic torsion-free $\gtwo$--structure $\varphi_0$ on $M$ and set $h_0=1$. We want to perturb this solution of the limiting equations \eqref{eq:GH:Spin(7):collapsing:sequence} with $\epsilon=0$ to a solution of the system with $\epsilon>0$.

To this end, we reinterpret \eqref{eq:GH:Spin(7):collapsing:sequence} as the vanishing of a nonlinear map $\Psi$ defined by
\begin{equation}\label{eq:Main:Equation}
\xi = (\varphi, h, \kappa) \longmapsto d\left( h^{\frac{2}{3}}\psi\right)+ \kappa\wedge\varphi.
\end{equation}
Here $\varphi$ is a closed basic $\gtwo$--structure, $\psi$ is its dual $4$-form, $h$ is a basic function and $\kappa$ is a basic closed $2$-form (satisfying additional conditions that we will impose below). In particular, note that $\Psi (\xi)$ is a closed $5$-form. 

The triple $\xi_0 = (\varphi_0, 1, 0)$ is a solution to \eqref{eq:Main:Equation}. In order to understand nearby solutions we are going to linearise \eqref{eq:GH:Spin(7):collapsing:sequence} at $\xi_0$. Consider a perturbation $\xi = \xi_0 + \zeta$ with $\zeta = (\rho, f, \kappa)$. We assume that $\rho$ is sufficiently small in $C^0(B)$--norm so that $\varphi = \varphi_0 + \rho$ stills defines a basic $\gtwo$--structure.  We write $\psi = \psi_0 + \hat{\rho} + Q_{\varphi_0}(\rho)$ for the dual $4$-form, where $\hat{\rho}$ is the image of $\rho$ under the linear map of Lemma \ref{lem:Linearisation:Hitchin:dual:3form:7d} and $Q_{\varphi_0}$ is a smooth map satisfying
\begin{equation}\label{eq:Estimate:quadratic:Hitchin:dual}
|Q_{\varphi_0}(\rho)|\leq C | \rho |^2, \qquad |\nabla Q_{\varphi_0}(\rho)|\leq C | \rho |\, |\nabla\rho|
\end{equation}
for a uniform constant $C$, see \cite[Proposition 10.3.5]{Joyce:Book}. Here the adapted connection and the norms are defined using the metric induced by $\varphi_0$. We then write $\Psi (\xi_0 + \zeta) = \mathcal{L}_0 (\zeta) + \mathcal{N}_0(\zeta)$, where
\begin{equation}\label{eq:Main:Linearisation}
\mathcal{L}_0 (\rho, f, \kappa) = d\left( \hat{\rho} + \tfrac{2}{3}f\,\psi_0\right) +\kappa\wedge\varphi_0
\end{equation}
is linear and
\begin{equation}\label{eq:Main:Nonlinear}
\mathcal{N}_0(\rho, f, \kappa) = d\left( (1+f)^{\frac{2}{3}}\big( \psi_0 + \hat{\rho} + Q_{\varphi_0}(\rho) \big) - \psi_0 - \hat{\rho} - \tfrac{2}{3}f\,\psi_0  \right) + \kappa\wedge\rho
\end{equation}
contains the nonlinearities.

Suppose we are given a \emph{bounded} solution $\zeta_0 = (\rho_0, f_0, \kappa_0)$ to the linearised equation $\mathcal{L}_0 (\zeta_0)=0$. Then for $\epsilon>0$ small $\xi_0 + \epsilon\zeta_0 = (\varphi_0 + \epsilon\rho_0, 1+\epsilon f_0, \epsilon\kappa_0)$ is an approximate solution to \eqref{eq:Main:Equation}. Here we assume that $\epsilon$ is sufficiently small to ensure that $\varphi_\epsilon = \varphi_0 + \epsilon\rho_0$ is still a basic $\gtwo$--structure. In order for $\xi_\epsilon$ to be a geometrically meaningful approximate solution we must require that the basic closed form $\kappa_0$ represents the orbifold first Chern class of the Seifert circle bundle $\pi\co M\ra B$, \ie $\kappa_0 = d\theta_0$ is the curvature of a connection $1$-form $\theta_0$ on the Seifert bundle. Note also that, linearising \eqref{eq:G2:monopole:eq}, $(f_0, \theta_0)$ satisfies $\ast df_0 + d\theta_0 \wedge\psi_0=0$, \ie $(f_0, \theta_0)$ is an abelian $\gtwo$--monopole. In particular, $f_0$ is harmonic and, in the context of this paper, the assumption that $\zeta_0$ is bounded forces $f_0$ to be constant. In this case it makes sense to further normalise $\zeta_0$ by requiring that $f_0\equiv 0$. We work under this assumption in the rest of the section.

\begin{remark*}
It is also interesting to consider the case where $\zeta_0$ is \emph{not} bounded, which corresponds to families of $\spins$--metrics collapsing in codimension $1$ with unbounded curvature. A natural assumption is then to assume that $(f_0, \theta_0)$ is an abelian $\gtwo$--monopole with Dirac-type singularities along a coassociative submanifold. This more challenging case is beyond the scope of this paper. 
\end{remark*}

We now aim to exponentiate the infinitesimal deformation $\zeta_0$ to an exact solution to \eqref{eq:Main:Equation}. To this end, we look for a further perturbation $\xi_0 + \epsilon\zeta_0 + \zeta$ such that
\begin{equation}\label{eq:Main:Equation:Perturbed}
\Psi (\xi_0 + \epsilon\zeta_0) + \mathcal{L}_\epsilon (\zeta) + \mathcal{N}_\epsilon (\zeta)=0,
\end{equation}
where (using $\mathcal{L}_0(\rho_0)=0$ and $f_0\equiv 0$)
\begin{equation}\label{eq:Error}
\Psi (\xi_0 + \epsilon\zeta_0) = dQ_{\varphi_0}(\epsilon\rho_0) + \epsilon^2 d\theta_0\wedge\rho_0
\end{equation}
and $\mathcal{L}_\epsilon$ and $\mathcal{N}_\epsilon$ are defined by the same formulas \eqref{eq:Main:Linearisation} and \eqref{eq:Main:Nonlinear} with $\varphi_\epsilon$ in place of $\varphi_0$. Since $\varphi_\epsilon = \varphi_0 + \epsilon\rho_0$ is arbitrarily $C^0$--close to $\varphi_0$ as $\epsilon\ra 0$, the solvability of \eqref{eq:Main:Equation:Perturbed} reduces to showing that (under additional conditions on $\zeta$) $\mathcal{L}_0(\zeta)= \Psi (\xi_0 + \epsilon\zeta_0)$ has a solution and $\mathcal{L}_0$ can be inverted on the image of $\mathcal{N}_\epsilon$.

\subsubsection{Implementation of the adiabatic limit strategy}

In the rest of the section we exploit the analysis on AC $\gtwo$--orbifolds we have developed earlier in the paper to implement this strategy. Let $\pi\co M^8\ra B$ be a principal Seifert circle bundle with connection $1$-form $\theta_0$. We assume that $M$ carries a basic torsion-free $\gtwo$--structure $\varphi_0$ inducing a transversally AC metric on $M$. There are three main points we need to address in order to apply the Implicit Function Theorem:
\begin{enumerate}
\item construct a bounded solution of the linearised problem $\mathcal{L}_0(\zeta_0)=0$;
\item understand the mapping properties of $\mathcal{L}_0$;
\item prove that $\mathcal{L}_0$ can be inverted on $\Psi (\xi_0 + \epsilon\zeta_0)$ and on the image of $\mathcal{N}_\epsilon$.
\end{enumerate}

\subsubsection*{The infinitesimal deformation}

As explained after Remark \ref{rmk:Vanishing:Harmonic:0:1:Forms} we can assume without loss of generality that $d\theta_0\in\Omega^2_{14}(B)$. By Lemma \ref{lem:decomp:forms} this means that $\ast d\theta_0 = -d\theta_0\wedge\varphi_0$ and therefore $d\theta_0$ is closed and coclosed. Moreover, $|\nabla^k (d\theta_0)|=O(r^{-2-k})$ for all $k\geq 0$.

We want to find a basic closed $3$-form $\rho_0$ such that $\mathcal{L}_{0}(\rho_0,0,d\theta_0)=0$, where $\mathcal{L}_0$ is the linear operator \eqref{eq:Main:Linearisation}. Since $d\theta_0\in\Omega^2_{14}(B)$, the equation $\mathcal{L}_{0}(\rho_0,0,d\theta_0)=0$ is equivalent to
\begin{equation}\label{eq:infinit:def}
d\rho_0=0, \qquad d\hat{\rho}_0 -\ast d\theta_0=0,
\end{equation}
which is not obviously elliptic as $\hat{\rho}_0$ involves the basic Hodge-star operator and the type decomposition of basic forms.

Fix $\mu=-1+\delta$ for an arbitrarily small $\delta>0$ and consider instead the elliptic equation $\triangle\sigma_0=d\theta_0$ for a basic $2$-form $\sigma_0$ of class $C^{\infty}_{\mu+1}(B)$. By Proposition \ref{prop:Laplacian:2:forms}, if a solution $\sigma_0$ exists then $d\theta_0=d^\ast d\sigma_0$ and $d\sigma_0\in\Omega^3_{27}(B)$. We would therefore conclude that $\rho_0=d\sigma_0$ solves \eqref{eq:infinit:def}. 

By Proposition \ref{prop:Laplacian:2:forms} the equation $\triangle\sigma_0=d\theta_0$ has a solution if and only if $d\theta_0$ is $L^2(B)$--orthogonal to basic closed and coclosed $2$-forms on $M$ in $\mathcal{H}^{2}_{-6-\mu}(B)$. Now, since $-6-\mu=-5+\delta<-2$ we have $\mathcal{H}^{2}_{-6-\mu}(B)\simeq H^2_c(B)$ by Theorem \ref{thm:L2:cohomology} (i) and Lemma \ref{lem:Excluded:indicial:roots:kforms} (iii). Therefore by duality a solution $\sigma_0$ exists if and only if $[\ast d\theta_0] = -[d\theta_0\wedge\varphi_0]=0\in H^5(B)$. Note that this condition coincides with the necessary topological constraint \eqref{eq:top:obstr}.

\begin{remark}\label{rmk:Assumption}
Observe that our standing assumption $d\theta_\infty\neq 0$ is in fact necessary for \eqref{eq:top:obstr} to be satisfied in a non-trivial way. For otherwise $d\theta_0$ would represent an $L^2$ basic cohomology class by Theorem \ref{thm:L2:cohomology} (i) and therefore $\ast d\theta_0=-d\theta_0\wedge\varphi_0$ could never be exact unless $d\theta_0=0$. However, if $d\theta_0=0$ then $\Phi_\epsilon = \epsilon \theta_0 \wedge\varphi_0 + \psi_0$ is already torsion-free for all $\epsilon>0$, albeit locally reducible and therefore not interesting from a $\spins$ geometry point of view. In fact, in this case one can argue that $M=(B\times S^1)/\Gamma$ for a freely-acting finite group $\Gamma$ and a \emph{smooth} AC $\gtwo$--manifold $B$.
\end{remark}

\subsubsection*{Mapping properties of the linearisation}

We now consider the linear operator $\mathcal{L}_0$ acting on the space of triples $(\rho,f,\kappa)$ of a basic closed $3$-form $\rho$, a basic function $f$ and a basic closed $2$-form $\kappa$. In fact we must further restrict $\kappa$ to be exact, \ie we vary $d\theta_0$ in the fixed basic cohomology class $c^{orb}_1(M)$.

We now fix $\nu=-2+\delta$ for an arbitrarily small $\delta>0$ and $\alpha\in (0,1)$. Given a basic $5$-form $\tau$ of class $C^{0,\alpha}_{\nu-1}(B)$, we look for a basic closed $3$-form $\rho$, basic $1$-form $\eta$ and basic function $f$ of class $C^{1,\alpha}_\nu(B)$ such that
\[
\mathcal{L}_0(\rho ,f, d\eta) = d\left( \hat{\rho} + \tfrac{2}{3}f\psi_0 + \eta\wedge\varphi_0\right) = \tau.
\]
Clearly $\tau$ must be closed (in fact, exact) if a solution exists.

Corollary \ref{cor:Exact:5:forms} describes sufficient conditions to solve this equation. If $\tau$ is closed and $L^2$--orthogonal to $\mathcal{H}^5_{-6-\nu}(B)$, then there exist a basic $2$-form $\sigma$ and a basic $1$-form $\gamma$ of class $C^{2,\alpha}_{\nu+1}(B)$ such that $\tau=\mathcal{L}_0 \left( d\sigma,-\tfrac{3}{2}d^\ast\gamma , d\,\tu{curl}\, \gamma \right)$.

\subsubsection*{The nonlinear equation}

This discussion suggests we now write $\zeta = \left( d\sigma,-\tfrac{3}{2}d^\ast\gamma , d\,\tu{curl}\, \gamma \right)$ and consider \eqref{eq:Main:Equation:Perturbed} as an equation for a basic $2$-form $\sigma$ and a basic $1$-form $\gamma$ of class $C^{2,\alpha}_{\nu+1}(B)$, where $\nu=-2+\delta$ with $\delta>0$ small.

In order to control the nonlinearities we use Remark \ref{rmk:Weighted:Embedding} (iv): if $\sigma$ and $\gamma$ are of class $C^{2,\alpha}_{\nu+1}(B)$ with $\nu = -2+\delta$, then
\[
\mathcal{N}_\epsilon \left( d\sigma,-\tfrac{3}{2}d^\ast\gamma , d\,\tu{curl}\, \gamma \right) = d\left( (1+f)^{\frac{2}{3}}\big( \psi_\epsilon + \hat{\rho} + Q_{\varphi_\epsilon}(\rho) \big) - \psi_\epsilon - \hat{\rho} - \tfrac{2}{3}f\,\psi_\epsilon + \tu{curl}\, \gamma\wedge d\sigma\right) 
\]
lies in $d\, C^{1,\alpha}_{\nu}(B)$. Thus $\mathcal{L}_0$ can be inverted on the image of $\mathcal{N}_\epsilon$ by Remark \ref{rmk:Exact:5:forms}.

Similarly, using the fact that $\rho_0=d\sigma_0$ is exact, we observe that
\[
\Psi (\xi_0 + \epsilon\zeta_0) = d\left( Q_{\varphi_0}(\epsilon\rho_0) + \epsilon^2 d\theta_0\wedge\sigma_0\right).
\]
Moreover, $|Q_{\varphi_0}(\epsilon\rho_0) + \epsilon^2 d\theta_0\wedge\sigma_0| \leq C r^{-2+\delta}$ for $\delta>0$ sufficiently small by \eqref{eq:Estimate:quadratic:Hitchin:dual} and Theorem \ref{thm:Weighted:Embedding} (v). Therefore the error \eqref{eq:Error} is also in the image of $\mathcal{L}_0$ by Remark \ref{rmk:Exact:5:forms}.

\subsubsection*{The existence result}

We have now all the ingredients to apply the Implicit Function Theorem and guarantee the existence of solutions to \eqref{eq:Main:Equation:Perturbed} for $\epsilon>0$ sufficiently small. We now summarise the resulting existence result for solutions to \eqref{eq:GH:Spin(7)}.

Let $(N,\pi_\infty,\theta_\infty,\omega,\Omega)$ be a nearly K\"ahler principal Seifert circle bundle and $(M,\pi,\theta_0,\varphi_0)$ be a transversally AC $\gtwo$--holonomy principal Seifert circle bundle asymptotic to $\BC (N)$. Given $\epsilon>0$ we modify the metric on $\BC (N)$ by rescaling $g_\BC$ in the direction of the circle orbits:
\[
g_{\tu{BC},\epsilon} = dr^2 + r^2 g_\Sigma + \epsilon^2 \theta_\infty^2.
\]
A Riemannian metric $g$ on $M$ is called \emph{ALC (asymptotically locally conical)} if there exists a compact set $K\subset M$, $\epsilon>0$, $R>0$ and a diffeomorphism $f\co (R,\infty)\times N\ra M\setminus K$ such that
\[
|\nabla^k(f^\ast g- g_{\tu{BC},\epsilon})| = O(r^{-k+\mu})
\] 
for all $k\geq 0$ and some $\mu<0$. Here covariant derivatives and norms are computed using the Riemannian metric $g_{\tu{BC},\epsilon}$. 

\begin{theorem}\label{thm:Main:Technical}
Let $\pi\co M^8\ra B$ be a principal Seifert circle bundle endowed with a transversally AC basic torsion-free $\gtwo$--structure $\varphi_0$ and let $\theta_0$ be the (unique up to diffeomorphisms) connection $1$-form on $\pi$ such that $d\theta_0\in\Omega^2_{14}(B)$.

If $\pi\co M\ra B$ is non-trivial and in basic cohomology 
\[
[d\theta_0\wedge\varphi_0]=0\in H^5(B)
\] 
then there exists $\epsilon_0=\epsilon_0 (M,\varphi_0)$ with the following significance. For all $\epsilon \in (0, \epsilon_0)$ there exists an $S^1$--invariant torsion-free $\spins$--structure $\Phi_\epsilon$ on $M$ such that
\begin{enumerate}
\item the induced metric $g_{\Phi_\epsilon}$ has holonomy $\spins$ and ALC asymptotics;
\item as $\epsilon\ra 0$, $(M,g_{\Phi_\epsilon})$ is arbitrarily close to $g_{\varphi_0} + \epsilon^2 \theta_0^2$ in $C^{k,\alpha}_{loc}$ for every $k\geq 0$. In particular, $(M,g_{\Phi_\epsilon})$ collapses with bounded curvature to the orbifold $(B,g_{\varphi_0})$ as $\epsilon\ra 0$.
\end{enumerate}
\proof

The existence of torsion-free $S^1$--invariant $\spins$--structures $\Phi_\epsilon$ on $M$ for all $\epsilon>0$ follows from the previous discussion. We start with the solution $(\varphi_0,1,0)$ of the limiting equation \eqref{eq:Main:Equation:Perturbed} with $\epsilon=0$. We have explained how to construct a solution $(d\sigma_0,0,d\theta_0)$ to the linearised problem $\mathcal{L}_0 (d\sigma_0,0,d\theta_0)$. We then consider \eqref{eq:Main:Equation:Perturbed} as an equation for a basic $2$-form $\sigma$ and a basic $1$-form $\gamma$ of class $C^{2,\alpha}_{\nu+1}(B)$, where $\nu=-2+\delta$ with $\delta>0$ small. An application of the Implicit Function Theorem (in the quantitative version stated for example in \cite[Lemma 1.3]{Biquard:Minerbe}) yields the existence of a torsion-free $\spins$--structure $\Phi_\epsilon$ on $M$ for every $\epsilon>0$ sufficiently small. The induced $\spins$--metric $g_{\Phi_\epsilon}$ is ALC (with $\mu$ the maximum between the rate of decay of $(\pi,\theta,\varphi)$ to $(\pi_\infty,\theta_\infty,\varphi_{\tu{C}})$ and $-1+\delta$) and satisfies the limiting behaviour in (ii) only in a $C^{1,\alpha}$--sense. Higher regularity follows from elliptic regularity since $\Phi_\epsilon$ is a closed and coclosed form with respect to a metric which differs from the model metric $g_{\BC,\epsilon}$ by $C^{1,\alpha}_{\mu}$--terms.

Finally, in order to prove that $g_{\Phi_\epsilon}$ has holonomy $\spins$ we use \cite[Lemma 2]{Bryant:1987}, which states that the holonomy of a metric induced by a torsion-free $\spins$--structure on a simply connected $8$-manifold $M$ is equal to $\spins$ if and only if there are no parellel $1$-forms and $2$-forms on $M$. By Proposition \ref{prop:fund:gp} (i), up to a finite cover we can assume that $M$ is simply connected. We first consider parallel forms of degree $1$ and $2$ on $\tu{BC}(N)$, since these determine the asymptotic behaviour of parallel forms on $M$. In fact, since $|d\theta_\infty|_{g_\BC}=O(r^{-2})$ we can consider forms on $\tu{BC}(N)$ parallel with respect to the adapted connection. Since the holonomy of the adapted connection of $\tu{BC} (N)$ is $\gtwo$, there are no parallel $2$-forms and the only parallel $1$-form is $\theta_\infty$. Since we assume that $d\theta_\infty\neq 0$, however, $\theta_\infty$ cannot extend to a parallel $1$-form on $M$. Hence every parallel form on $M$ of degree $1$ and $2$ must decay and therefore vanish.   
\endproof
\end{theorem}

\section{Examples}\label{sec:Examples}

In this final section we use our existence result Theorem \ref{thm:Main:Technical} to construct concrete examples of complete $\spins$--metrics. Given the currently limited knowledge of AC $\gtwo$--manifolds, using orbifolds is essential. We are able to use Theorem \ref{thm:Main:Technical} to produce infinitely many different topological types of complete $\spins$--metrics and examples of $8$-manifolds carrying infinitely many distinct families of ALC $\spins$--metrics (Theorems \ref{thm:Examples:Betti} and \ref{thm:Examples:Families} in the Introduction). Previously only a handful of complete non-compact $\spins$--metrics was known.

In Section \ref{sec:ALC:G2} we use the analysis on AC orbifolds we have developed in this paper to extend the construction in \cite{FHN:ALC:G2:from:AC:CY3} of ALC $\gtwo$--manifolds from AC Calabi--Yau $3$-folds to the case of AC Calabi--Yau orbifolds. As an illustrative example, we use this extension to produce infinitely many distinct families of ALC $\gtwo$--metrics on $S^3\times\R^4$. 

\subsection{Self-dual Einstein $4$-orbifolds and Bryant--Salamon's AC $\gtwo$--metrics}

The starting point for the construction of Theorem \ref{thm:Main:Technical} is an AC $\gtwo$--manifold or orbifold satisfying appropriate topological conditions. AC $\gtwo$--manifolds are hard to construct. All currently known examples of AC $\gtwo$--manifolds admit a symmetry group that acts with cohomogeneity one, \ie with generic orbits of codimension 1. The large symmetry group affords a reduction of the PDE for the holonomy reduction to $\gtwo$ to a system of nonlinear ODEs. Studying solutions to these ODE systems is still non-trivial: in 1989 Bryant--Salamon \cite{Bryant:Salamon} constructed three explicit AC $\gtwo$--metrics, but only very recently \cite[Theorem C]{ALC:G2:coho1} have further examples of AC $\gtwo$--manifolds been found; these examples are not explicit and their existence is based on the qualitative analysis of the relevant ODE system.

From a different point of view, Bryant--Salamon's examples of AC $\gtwo$--manifolds in \cite{Bryant:Salamon} fall into the class of constructions, pioneered by Calabi \cite{Calabi:AC:CY} in the Calabi--Yau and hyperk\"ahler case, of complete Ricci-flat metrics on total spaces of vector bundles over compact manifolds satisfying appropriate curvature conditions. In particular, Bryant--Salamon's construction yields an AC $\gtwo$--metric (unique up to scale) on the total space of the bundle of anti-self-dual $2$-forms over a self-dual Einstein $4$-manifold $Q$ with positive scalar curvature. By a theorem of Hitchin \cite{Hitchin:SDE+} the only such manifolds are $S^4$ and $\C\PP^2$ with their standard metrics. On the other hand, there are infinitely many self-dual Einstein $4$-orbifolds with positive scalar curvature and Bryant--Salamon's construction extends immediately (as observed in the introduction of \cite{Bryant:Salamon}) to the case where $Q$ is an orbifold and $B$ is the total space for the orbibundle of anti-self-dual $2$-forms on $Q$.

The most powerful known method of construction of self-dual Einstein $4$-orbifolds is the quaternionic K\"ahler quotient construction of Galicki--Lawson \cite{Galicki:Lawson}. We briefly recall this here. The construction is based on the tight relationship between quaternionic K\"ahler and hyperk\"ahler geometry \cite{Swann:Bundle}. Let $\tu{C}$ be a hyperk\"ahler cone acted upon by a group $K$ of triholomorphic isometries. Using the conical structure one can always find a hyperk\"ahler moment map $\mu\co\tu{C}\ra\Lie{k}^\ast\otimes\Imag\HH$, see \cite[Proposition 13.6.1]{Boyer:Galicki}. Besides the triholomorphic action of $K$, there is an action of $\HH^\ast$ on $\tu{C}$ generated by the Euler vector field $r\partial_r$. Letting $\HH^\ast$ act on $\Lie{k}^\ast\otimes\Imag\HH$ via conjugation on $\Imag\HH$, the moment map $\mu$ is equivariant with respect to the action of $K\cdot\HH^\ast$, where $K\cdot\HH^\ast=(K\times\HH^\ast)/\Z_2$ if $-1\in K$ and $K\cdot\HH^\ast=K\times\HH^\ast$ otherwise. The hyperk\"ahler quotient construction \cite{HK:Quotient} yields the existence of a hyperk\"ahler stucture on (the smooth part of) $\mu^{-1}(0)/K$. Using the $\HH^\ast$--equivariance property of the moment map one can show that $\mu^{-1}(0)/K$ is a new hyperk\"ahler cone $\tu{C}'$.

Now, to each hyperk\"ahler cone $\tu{C}$ there is a naturally associated positive quaternionic K\"ahler ``space'' $Q$, obtained as the quotient of $\tu{C}$ by the $\HH^\ast$ action generated by $r\partial_r$. Here we say that $Q$ is a positive quaternionic K\"ahler space if its smooth part carries a Riemannian metric with holonomy contained in $\Sp{1}\cdot\Sp{n}$ (in particular the dimension of $Q$ must be a multiple of $4$); any such metric is Einstein and the qualification `positive' refers to the sign of the Einstein constant. In dimension $4$ this definition must be modified: we say that $Q^4$ is quaternionic K\"ahler if it is self-dual and Einstein. Because of the $\HH^\ast$--equivariance of the moment map $\mu\co \tu{C}\ra\Lie{k}^\ast\otimes\Imag\HH$, the construction of the hyperk\"ahler cone $\tu{C}'$ as a hyperk\"ahler quotient of the cone $\tu{C}$ yields a quaternionic K\"ahler structure on $Q'=\tu{C}'/\HH^\ast$ as a quaternionic K\"ahler quotient of $Q=\tu{C}/\HH^\ast$ \cite{Galicki:Lawson}.

Many interesting self-dual Einstein $4$-orbifolds with positive scalar curvature can be obtained in this way even from the simplest hyperk\"ahler cone, $\tu{C}=\HH^{n+1}$. (Here we regard $\HH ^{n+1}$ as the hyperk\"ahler cone over the round sphere $\Sph^{4n+3}$; the associated quaternionic K\"ahler manifold is $\HH\PP^n = \Sph^{4n+3}/\tu{Sp}(1)$.) For example, all \emph{toric} self-dual Einstein $4$--orbifolds with positive scalar curvature, \ie those self-dual Einstein $4$-orbifolds with a $T^2$--symmetry, arise in this way \cite{Calderbank:Singer}, see \cite[\S\S 12.4-12.5 and 13.7]{Boyer:Galicki} for further details.

\begin{remark*}
There are also self-dual Einstein $4$-orbifolds with positive scalar curvature that are \emph{not} obtained via quaternionic K\"ahler reduction, for example the families of $\sorth{3}$--invariant self-dual Einstein orbifold metrics on $S^4$ and $\C\PP^2$ constructed by Hitchin in \cite{Hitchin:AH:SDE+}. We will not use these metrics in this paper.
\end{remark*}

Now, in view of the assumptions of Theorem \ref{thm:Main:Technical}, amongst all self-dual Einstein $4$-orbifolds with positive scalar curvature we are interested in those that satisfy the following additional property.

\begin{definition}\label{def:Admissible:SDE+:orbifolds}
Let $Q$ be a self-dual Einstein $4$-orbifold with positive scalar curvature. We say that $Q$ is \emph{\spins--admissible} if there exists a principal Seifert circle bundle $S\ra Q$.
\end{definition}

The relevance of this assumption is explained by the following lemma.

\begin{lemma}\label{lem:Spin7:admissible}
Let $Q$ be a self-dual Einstein $4$-orbifold with positive scalar curvature and denote by $B$ the total space of the orbibundle of anti-self-dual $2$-forms on $Q$ endowed with Bryant--Salamon AC torsion-free $\gtwo$--structure $\varphi_0$. Then $Q$ is \spins--admissible if and only if there exists a principal Seifert circle bundle $\pi\co M\ra B$ such that $c_1^{orb}(M)\cup[\varphi_0]=0\in H^5_{orb}(B)$.
\proof
If $\pi\co S\ra Q$ is a principal Seifert circle bundle and $p\co B\ra Q$ is the orbibundle of anti-self-dual $2$-forms then $M=p^\ast S$ is a principal circle orbibundle over $B$. Over a small enough uniformising chart $U/\Gamma\subset Q$ we can trivialise both $\pi$ and $p$. Then there exist representations of $\Gamma$ in $\unitary{1}$ and $\sorth{3}$ such that $\pi^{-1}(U/\Gamma)= (U\times S^1)/\Gamma$ and locally $M$ can be described as $(U\times\R^3\times S^1)/\Gamma$. If $S$ is smooth then $\Gamma$ acts freely on $S^1$ and therefore $M$ is smooth as well. Moreover, the topological constraint $c_1^{orb}(M)\cup[\varphi_0]=0$ is automatically satisfied since $H^5_{orb}(B)\simeq H^5_{orb}(Q)=0$.

Conversely, if $\pi\co M\ra B$ is a principal Seifert circle bundle over the total space of the orbibundle of anti-self-dual $2$-forms on $Q$, then $Q$ is \spins--admissible since the restriction $S$ of the orbibundle $M\ra B$ to the zero-section $Q$ in $B$ yields a principal Seifert circle bundle $S\ra Q$.
\endproof
\end{lemma}

\subsection{Concrete examples of ALC $\spins$--metrics}

In order to apply our existence result Theorem \ref{thm:Main:Technical} we need to understand to what extent known constructions of self-dual Einstein $4$-orbifolds with positive scalar curvature yields examples that are \spins--admissible. In the rest of the section we will limit ourselves to discuss families of examples that give a sense of the rich variety of ALC $\spins$--metrics that can be obtained using Theorem \ref{thm:Main:Technical} and defer a more systematic study to elsewhere. We discuss three sets of examples. Our first theorem provides a proof of the existence of a $1$-parameter family of ALC $\spins$--metrics conjectured by Gukov--Sparks--Tong \cite{Gukov:Sparks:Tong}; currently this is the only example that can be obtained from Theorem \ref{thm:Main:Technical} starting from a smooth AC \gtwo--manifold. Using AC $\gtwo$--orbifolds, we then prove that the same smooth $8$-manifold in fact carries infinitely many distinct families of ALC \spins--metrics. Finally, we construct infinitely many smooth $8$-manifolds carrying complete ALC $\spins$--metrics.

\subsubsection{An example from a smooth AC $\gtwo$--manifold}

Amongst the known examples of smooth AC $\gtwo$--manifolds, only $\Lambda^-T^\ast\C\PP^2$ endowed with Bryant--Salamon's AC $\gtwo$--metric can be used in Theorem \ref{thm:Main:Technical}. Indeed, the other two Bryant--Salamon examples of AC $\gtwo$--manifolds \cite{Bryant:Salamon} have vanishing second cohomology, while all the infinitely many examples constructed in \cite{ALC:G2:coho1} only have compactly supported second cohomology. Thus the topological constraint \eqref{eq:top:obstr} can be satisfied in a non-trivial way only in the case of $\Lambda^-T^\ast\C\PP^2$.

\begin{theorem}\label{thm:Gukov:Sparks:Tong}
The total space of the non-trivial rank-$3$ real vector bundle over $S^5$ carries a $1$-parameter family of ALC $\spins$--metrics. The Bryant--Salamon AC $\gtwo$--metric on $\Lambda^-T^\ast\C\PP^2$ arises as a collapsed limit of this family.
\proof
The Bryant--Salamon AC $\gtwo$--manifold $B=\Lambda^-T^\ast\C\PP^2$ has torsion-free $1$-dimensional second cohomology and contains no $5$-cycle. Hence up to a change of orientation and finite quotients there is a unique non-trivial circle bundle $M$ over $B$ and the topological constraint \eqref{eq:top:obstr} is automatically satisfied. The Bryant--Salamon AC $\gtwo$--metric is $\sunitary{3}$--invariant and this $\sunitary{3}$--action lifts to an isometric action of the ALC $\spins$--metrics produced by Theorem \ref{thm:Main:Technical}. We give a description of all the manifolds involved in terms of this $\sunitary{3}$--action: $\C\PP^2 = \sunitary{3}/\unitary{2}$ and therefore $B=\sunitary{3}\times_{\unitary{2}}\Lie{su}_2$ (here the action of $\unitary{2}$ on $\Lie{su}_2$ is induced by the adjoint representation); the unique simply connected circle bundle over $\C\PP^2$ is $S^5=\sunitary{3}/\sunitary{2}$ and the $8$-manifold carrying ALC $\spins$--metrics by Theorem \ref{thm:Main:Technical} is $M=\sunitary{3}\times_{\sunitary{2}}\Lie{su}_2$.
\endproof
\end{theorem}

\begin{remark*}
The existence of this $1$-parameter family of ALC $\spins$--metrics was conjectured by Gukov--Sparks--Tong \cite{Gukov:Sparks:Tong}. The family is expected to be part of a geometric transition in $\spins$--geometry which physically corresponds to a duality between Type IIA String Theory on $\Lambda^-T^\ast\C\PP^2$ with D6 branes/Ramond--Ramond fluxes. The family of ALC metrics of Theorem \ref{thm:Gukov:Sparks:Tong} provides the M theory lift of Type IIA Theory on $\Lambda^-T^\ast\C\PP^2$ with fluxes, while the lift of Type IIA Theory on $\Lambda^-T^\ast\C\PP^2$ with a D6-brane wrapping the coassociative $\C\PP^2$ corresponds to an explicit ALC $\spins$--metric found by Gukov--Sparks \cite{Gukov:Sparks} and its conjectural $1$-parameter family of deformations up to scale. A complete geometric explanation of the physical duality would involve constructing AC $\spins$--metrics arising as limits of the two $1$-parameter families of ALC $\spins$--metrics as well as an ALC $\spins$--metric on $\R_+\times\sunitary{3}/\unitary{1}$ with an isolated conical singularity. This is completely analogous to the analysis of \cite{ALC:G2:coho1} in the $\gtwo$ setting. In fact, since the metrics in Theorem \ref{thm:Gukov:Sparks:Tong} admit a cohomogeneity one action of $\sunitary{3}$, it is likely that the methods of \cite{ALC:G2:coho1} can address these conjectures.
\end{remark*}

\subsubsection{Infinitely many families of ALC $\spins$--metrics on the same smooth $8$-manifold}

The first examples of self-dual Einstein $4$-orbifolds obtained by Galicki--Lawson \cite{Galicki:Lawson} via the quaternionic K\"ahler construction were weighted projective planes arising from quotients of $\HH\PP^2$ by a circle. In this section we use these examples to produce infinitely many distinct families of ALC \spins--metrics on the smooth $8$-manifold of Theorem \ref{thm:Gukov:Sparks:Tong}.

In order to describe the main features of Galicki--Lawson's self-dual Einstein metrics we follow \cite[\S\S 7 and 8]{Boyer:Galicki:circle:reduction} (that generalises the construction to circle quotients of quaternionic projective spaces of arbitrary dimension), see also \cite[Proposition 12.5.3]{Boyer:Galicki}. Fix $p_1,p_2,p_3\in\Z_{>0}$ with $\gcd{(p_1,p_2,p_3)}=1$ and consider the circle action $e^{i\theta}\cdot [u_1:u_2:u_3]= [e^{ip_1\theta}u_1:e^{ip_2\theta}u_2:e^{ip_3\theta}u_3]$ on $\HH\PP^2$. The quaternionic K\"ahler quotient of $\HH\PP^2$ by the circle action is
\[
Q = \{ [u_1:u_2:u_3]\in\HH\PP^2\, |\, p_1 \overline{u}_1 i u_1 + p_2 \overline{u}_2 i u_2 + p_3 \overline{u}_3 i u_3=0\}/S^1.
\]
By \cite[Proposition 7.5]{Boyer:Galicki:circle:reduction} the equation $p_1 \overline{u}_1 i u_1 + p_2 \overline{u}_2 i u_2 + p_3 \overline{u}_3 i u_3=0$ cuts out in $\Sph^{11}\subset\HH^3$ a smooth $8$-manifold diffeomorphic to the complex Stiefel manifold $V_2 (\C^3) \simeq \sunitary{3}$. There is an action of $S^1\cdot\sunitary{2}$ on $V_2 (\C^3)$ arising from the circle action on $\HH\PP^2$ and the standard $\tu{Sp}(1)$--action on the $3$--Sasakian $\Sph^{11}$. The quotient $Q=V_2 (\C^3)/S^1\cdot\sunitary{2}$ then has a quaternionic K\"ahler metric. Note that in general only the maximal torus of $\sunitary{3}$ commutes with the $S^1$--action on $V_2 (\C^3)$ and therefore, contrary to the example of Theorem \ref{thm:Main:Technical}, $Q$ is only toric and not $\sunitary{3}$--invariant. We also obtain two natural orbibundles over $Q$: the Konishi bundle $V_2 (\C^3)/S^1$ and the principal circle orbibundle $S^5=V_2 (\C^3)/\sunitary{2}$. Studying the fibre-wise $S^1$--action on $S^5$ one can show that $Q$ is isomorphic to the weighted projective plane $\mathbb{W}\C\PP^2[q_1,q_2,q_3]$, with $q_i = p_j+p_k$ if $p_1+p_2+p_3$ is odd and $2q_i = p_j+p_k$ otherwise. Here $(ijk)$ runs through cyclic permutations of $(123)$. It is now clear that $Q$ is \spins--admissible and that the $8$-dimensional principal Seifert circle bundle over $\Lambda^-T^\ast Q$ is $M=V_2 (\C^3)\times_{\sunitary{2}}\Lie{su}_2$.

\begin{theorem}\label{thm:Weighted:CP2}
The total space of the non-trivial rank-$3$ real vector bundle over $S^5$ carries infinitely many distinct families of ALC \spins--metrics.
\proof
The existence of highly collapsed ALC \spins--metrics follows from Theorem \ref{thm:Main:Technical} applied to the \spins--admissible self-dual Einstein $4$-orbifolds $Q$ obtained as quotients of $\HH\PP^2$ by circles labelled by $(p_1,p_2,p_3)$. The fact that, up to obvious symmetries, different choices of $(p_1,p_2,p_3)$ give rise to non-isometric families of \spins--metrics follows from the fact that their tangent cones at infinity, the Bryant--Salamon AC orbifold \gtwo--metrics on $\Lambda^-T^\ast Q$, are distinct.
\endproof
\end{theorem}

It is likely that many more examples of toric self-dual Einstein $4$-orbifolds are \spins--admissible. For example, an infinite family of orbifolds $Q$ with $2$-dimensional $H^2_{orb}(Q)$ (therefore not weighted projective planes) can be obtained by quaternionic K\"ahler quotient of $\tu{Gr}_2(\C^4)$ by a circle (these examples can also be described as particular reductions of $\HH\PP^3$ by a $2$-torus, since $\tu{Gr}_2(\C^4)$ is itself a circle quotient of $\HH\PP^3$). The generic $Q$ constructed in this way is \spins--admissible: indeed, the zero-level set of the $3$-Sasakian moment map in the Konishi bundle of $\tu{Gr}_2(\C^4)$, a circle orbibundle over the $4$-orbifold, is smooth for a generic choice of embedding of $S^1$ into the symmetry group $\sunitary{4}$ of $\tu{Gr}_2(\C^4)$. In general, however, we do not know how to recognise which toric self-dual Einstein $4$-orbifolds are \spins--admissible. In \cite{Boyer:Galicki:Inventiones} very clear combinatorial conditions for such an orbifold to admit a smooth $3$-Sasaki Konishi bundle are given and it is likely that similar combinatorial conditions characterise \spins--admissibility. Instead of pursuing such a systematic combinatorial approach, in the next section we construct by hand an explicit family of examples with unbounded second orbifold Betti number.

\begin{remark*}
Let $Q$ be a toric self-dual Einstein $4$-orbifold, which must be obtained as a quotient of $\HH\PP^n$ by an $(n-1)$-torus $T^{n-1}$ by \cite{Calderbank:Singer}. Denote by $\mu\co\HH^{n+1}\ra \R^{n-1}\otimes\Imag\HH$ the associated hyperk\"ahler moment map and recall that, thinking of $\mu$ as a section of the bundle $\Sph^{4n+3}\times_{\sunitary{2}}\Lie{su}_2\ra\HH\PP^n$, $Q=\mu^{-1}(0)/T^{n-1}$. The orbifold $Q$ is \spins--admissible if there exists a subtorus $T^{n-2}\subset T^{n-1}$ such that $S=\mu^{-1}(0)/T^{n-2}$ is a smooth $5$-manifold. If this were the case, then we could consider the principal $G$--bundle $P=\{ \triple{u}\in\Sph^{4n+3}\, |\, \mu(\triple{u})=0\}/T^{n-2}\ra S$, with $G=\sunitary{2}$ or $\sorth{3}$ depending on whether $-1\in T^{n-2}$ or not, and $M=P\times_{\sunitary{2}}\Lie{su}_2$ would carry ALC \spins--metrics by Theorem \ref{thm:Main:Technical}. In general $P$ is an orbifold. In \cite{Boyer:Galicki:Mann:Hypercomplex}, motivated by the fact that $P$ carries a natural hypercomplex structure, Boyer--Galicki--Mann determine conditions on $T^{n-2}\subset T^{n-1}$ under which $P$ is a smooth $8$-manifold, but they do not study which additional conditions guarantee that the action of $G$ on $P$ is free.
\end{remark*}

\subsubsection{Examples with arbitrarily large second Betti number from $A_n$ ALE spaces}

In this section we prove the existence of infinitely many diffeomorphism-types of simply connected $8$-manifolds carrying complete {\spins} metrics. The examples we will consider arise from an extension of Kronheimer's construction of ALE spaces \cite{Kronheimer,Kronheimer:Classification} to the quaternionic K\"ahler setting due to Galicki--Nitta \cite{Galicki:Nitta}.

Let $\Gamma$ be a finite subgroup of $\sunitary{2}$ acting freely on $\C^2\setminus\{0\}$. Kronheimer \cite{Kronheimer} constructs ALE hyperk\"ahler metrics on the minimal resolution of $\C^2/\Gamma$ using the hyperk\"ahler quotient construction. Let $R_0,\dots, R_r$ be the irreducible representations of $\Gamma$, with $R_0$ the trivial representation. Set $n_i=\dim R_i$. The regular representation $R$ of $\Gamma$ decomposes as $R=\bigoplus_{i=0}^r{\C^{n_i}\otimes R_i}$. Kronheimer considers $\tu{Hom}_\Gamma (R,R\otimes\HH)\simeq \HH^n$, where $\Gamma$ acts on $\HH\simeq\C^2$ via its embedding in $\sunitary{2}$. The McKay correspondence implies that $n=|\Gamma|$. Define $\hat{K}$ as the group of unitary transformations of $R$ that commute with the action of $\Gamma$; by the Schur Lemma $\hat{K}=\prod_{i=0}^r{\unitary{n_i}}$. Then $K=\hat{K}/\triangle\unitary{1}\simeq\prod_{i=1}^r{\unitary{n_i}}$ acts effectively and triholomorphically on $\HH^n$.  Let $\mu\co \HH^n\ra \Lie{k}^\ast\otimes\Imag\HH$ denote the hyperk\"ahler moment map for the action of $K$ on $\HH^n$ and let $\Lie{z}$ denote the Lie algebra of the centre of $K$. By the hyperk\"ahler quotient construction, for each $\triple{\zeta}\in\Lie{z}^\ast\otimes\Imag\HH$ the smooth part of the quotient $X_{\triple{\zeta}}=\mu^{-1}(\triple{\zeta})/K$ carries a natural hyperk\"ahler structure. Kronheimer shows that for generic $\triple{\zeta} \in\Lie{z}^\ast\otimes\Imag\HH$, $X_{\triple{\zeta}}$ is a smooth manifold diffeomorphic to the minimal resolution of $\C^2/\Gamma$ and that its natural hyperk\"ahler structure is asymptotic at infinity to the flat hyperk\"ahler structure on $\C^2/\Gamma$ (with rate $-4$). Conversely, in \cite{Kronheimer:Classification} Kronheimer shows that any asymptotically locally Euclidean (ALE) hyperk\"ahler $4$-manifold is obtained from this quotient construction.

For the extension of this construction to the quaternionic K\"ahler setting in \cite{Galicki:Nitta}, Galicki--Nitta think of $\triple{\zeta}\in \Lie{k}^\ast\otimes\Imag\HH$ as a map $\Lie{k}\ra\Imag\HH\simeq \Lie{su}_2$ and assume that there exists a group homomorphism $\rho\co K\ra \sunitary{2}$ with $\rho_\ast = -\triple{\zeta}$. We use $\rho$ to define an action of $K$ on $\HH\PP^{n}$ by
\[
g\cdot [u_0:\triple{u}] = [\rho(g)u_0: g\triple{u}].
\]
The quaternionic K\"ahler quotient of $\HH\PP^n$ by $K$ is $Q=\hat{\mu}^{-1}(0)/K$, where
\[
\hat{\mu}([u_0:\triple{u}]) = -\overline{u}_0 \triple{\zeta}u_0 + \mu (\triple{u}).
\]
Galicki--Nitta show \cite[Theorem 3.2]{Galicki:Nitta} that if $\triple{\zeta}$ is generic in the sense of Kronheimer (\ie $\mu^{-1}(\triple{\zeta})/K$ is smooth) then $\hat{\mu}^{-1}(0)/K$ is a quaternionic K\"ahler $4$-orbifold.

Now, denote by $K_\rho$ the kernel of $\rho$ and assume that $K/K_\rho\simeq\unitary{1}$, \ie $\rho\co K\ra\unitary{1}\subset\sunitary{2}$. In this case $S=\hat{\mu}^{-1}(0)/K_\rho\ra \hat{\mu}^{-1}(0)/K$ is a principal circle orbibundle over $Q$. Up to rotations and using a $K$--invariant metric to identify $\Lie{k}$ with its dual, we must have $\triple{\zeta}=2\pi i\zeta$, where $\zeta\in\Lie{z}$. Since $K=\prod_{i=1}^{r}{\unitary{n_i}}$ we can identify $\Lie{z}$ with $\R^r$ and since $\triple{\zeta}$ integrates to a group homomorphism $K\ra \unitary{1}$ we must have $\zeta\in\Z^r\subset\R^r$. We then define a $1$-dimensional representation $R_\zeta$ of $\Gamma$ by $R_\zeta = {\bigotimes_{i=1}^r\det (R_i)^{\zeta_i}}$.

\begin{prop}\label{prop:ALE:Spin7:admissible}
Assume that $\triple{\zeta}=2\pi i\zeta$ is generic in the sense of Kronheimer and that the homomorphism $\Gamma\ra\unitary{1}$ corresponding to the representation $R_\zeta$ is injective. Then $S=\hat{\mu}^{-1}(0)/K_\rho$ is smooth, \ie $Q=\hat{\mu}^{-1}(0)/K$ is a \spins--admissible quaternionic K\"ahler $4$-orbifold.
\proof
Write $\HH\PP^n$ as $\HH^n\cup\HH\PP^{n-1}$, where $\HH^n$ is identified with the open set where $u_0\neq 0$ and $\HH\PP^{n-1} = \{ [0:\triple{u}]\in\HH\PP^n\}$.

First work on the open set $\HH^n$. We introduce affine quaternionic coordinates $\triple{v}=\triple{u}\overline{u}_0$. Note that the open set $\hat{\mu}^{-1}(0)\cap\HH^n/K$ of $Q$ is identified with $\{ \triple{v}\in\HH^n\, |\, \mu (\triple{v}) = \triple{\zeta}\}/K$, where $K$ acts by $g\cdot\triple{v}=g\triple{v}\overline{\rho(g)}$. Thus a dense open set of $Q$ and the ALE manifold $X_{\triple{\zeta}}$ are obtained as quotients of $\mu^{-1}(\triple{\zeta})$ by $K$, but the $K$--action is different in the two cases. However, since the $K$--actions agree when restricted to $K_\rho$, $\hat{\mu}^{-1}(0)\cap \HH^n/K_\rho$ coincides with the principal circle bundle $\mu^{-1}(\triple{\zeta})/K_\rho\ra X_{\triple{\zeta}}$. In particular, the dense open set $\hat{\mu}^{-1}(0)\cap \HH^n/K_\rho$ of $S$ is smooth.

Consider now a point $[0:\triple{u}]\in\HH\PP^{n-1}$ such that $\mu(\triple{u})=0$. We must show that the stabiliser of $[0:\triple{u}]$ in $K_\rho$ is trivial. Now, in \cite[Lemma 3.1]{Kronheimer} Kronheimer shows that there exists a copy of $\C^2$ in $\mu^{-1}(0)\subset \HH^n$ such that every orbit of the $K$--action on $\mu^{-1}(0)$ meets a single orbit of the $\Gamma$--action on $\C^2$. Even if not explicitly mentioned in \cite{Kronheimer}, the identification of $\mu^{-1}(0)/K$ with $\C^2/\Gamma$ can be made equivariant with respect to the action of $\Sp{1}$ given by the (diagonal) right quaternionic multiplication on $\C^2=\HH$ and $\mu^{-1}(0)\subset\HH^n$. Thus $Q$ is obtained by adding a single point $\infty$ with isotropy $\Gamma$ to the open set $\hat{\mu}^{-1}(0)\cap\HH^n/K$. To show that $S$ is smooth, we need to show that the induced action of $\Gamma$ on the fibre $S^1$ over $\infty$ of the orbibundle $S\ra Q$ is free. This follows from the assumption that $R_\zeta$ is an effective $\Gamma$--representation, since the action of $\Gamma$ on the fibre of $S$ over $\infty$ is precisely given by $R_\zeta$ by \cite[Proposition 2.2 (ii)]{Kronheimer:Nakajima}.

In other words, $S$ is obtained by compactifying the principal circle bundle $\mu^{-1}(\triple{\zeta})/K_\rho\ra X_{\triple{\zeta}}$ by adding the circle fibre over the orbifold point in the natural orbifold compactification of $X_{\triple{\zeta}}$. The circle bundle $\mu^{-1}(\triple{\zeta})/K_\rho\ra X_{\triple{\zeta}}$ carries a natural anti-self-dual connection which is asymptotic to the flat connection on $\C^2/\Gamma$ with monodromy $R$. Hence if $R$ is an effective representation of $\Gamma$, $S$ is smooth.
\endproof
\end{prop}

The condition that $\Gamma\ra\unitary{1}$ is injective forces us to restrict to the abelian case $\Gamma=\Z_n$ for some $n\geq 2$. We can then be completely explicit. The irreducible representations of $\Gamma$ are all $1$-dimensional and labelled by an integer $0\leq i\leq n-1$. The group $K=T^{n-1}$ acts on $\HH^n$ by
\[
\left( e^{i\theta_1},\dots, e^{i\theta_{n-1}} \right) \cdot (u_1,\dots, u_n) = \left( e^{i\theta_1}u_1, e^{i(\theta_2-\theta_1)}u_2,\dots, e^{i(\theta_{n-1}-\theta_{n-2})}u_{n-1}, e^{-i\theta_{n-1}}u_n \right).
\]
The moment map $\mu\co\HH^n\ra \R^{n-1}\otimes\Imag\HH$ is therefore
\[
\mu (u_1,\dots, u_n)= (\overline{u}_1 i u_1 - \overline{u}_2 i u_2,\dots, \overline{u}_{n-1} i u_{n-1}- \overline{u}_n i u_n).
\]
Fix $\zeta = (\zeta_1,\dots,\zeta_{n-1})\in\Z^{n-1}$. Kronheimer's genericity conditions are, see for example \cite[Example 2.22]{Boyer:Galicki:Inventiones}, 
\begin{equation}\label{eq:An:genericity}
\zeta_i+\zeta_{i+1}+\dots + \zeta_{i+j}\neq 0,\qquad \text{for all } 1\leq i \leq n-1,\quad 0\leq j\leq n-1-i.
\end{equation}
Moreover, the representation $R_\zeta$ induces an injective homomorphism $\Gamma\ra\unitary{1}$ if and only if
\begin{equation}\label{eq:An:admissibility}
\gcd{\left( |\zeta|,n \right)}=1, \qquad |\zeta|=\zeta_1 + 2\,\zeta_2 + \dots + (n-1)\,\zeta_{n-1}
\end{equation}
We can also explicitly see that \eqref{eq:An:admissibility} guarantees that $K_\rho$ acts freely on points in $\HH\PP^n$ of the form $[0:\triple{u}]$ with $\mu (\triple{u})=0$. Consider the action of $T^n\times\HH$ on $\HH^n$ defined by $(e^{i\psi_1},\dots, e^{i\psi_n},u)\cdot\triple{u}=(e^{i\psi_1}u_1 \overline{u},\dots, e^{i\psi_n}u_n \overline{u})$. It is immediate to check that $\mu^{-1}(0)$ is the orbit of $(1,\dots, 1)$ and therefore $\{ u_0=0\}\cap\hat{\mu}^{-1}(0)$ reduces to a single point $[0:1:\dots:1]$. We now choose $(\tilde{\zeta}_1,\dots,\tilde{\zeta}_n)\in\Z^n$ such that $\tilde{\zeta}_i-\tilde{\zeta}_{i+1}=\zeta_i$ for all $i=1,\dots, n-1$. Note that $\tilde{\zeta}_1+\dots +\tilde{\zeta}_n \equiv |\zeta|$ modulo $n$. Then $K_\rho$ can be described as the subgroup of $T^n$ cut out by the constraints $e^{i(\psi_1+\dots + \psi_n)}=1=e^{i(\tilde{\zeta}_1\psi_1 + \dots + \tilde{\zeta}_n\psi_n)}$. We conclude that the stabiliser of $[0:1:\dots:1]$ in $K_\rho\times\HH$ consists of elements of the form $(\lambda,\dots, \lambda)$ with $\lambda\in S^1$ such that $\lambda^n=1=\lambda^{\tilde{\zeta}_1+\dots+\tilde{\zeta}_n}=\lambda^{|\zeta|}$. If $n$ and $|\zeta|$ are coprime then necessarily $\lambda=1$.

Note there exists a suitable choice for $\zeta$ for all $n\geq 2$. Indeed, if $n$ is odd consider $\zeta = (2,1,\dots, 1)$, while if $n\geq 4$ is even consider $\zeta = (2,1,\dots, 1, 2, 1, \dots, 1)$, where the second $2$ is the $\frac{n}{2}$th coordinate of $\zeta$. Since $\zeta_i>0$ for all $i$ the genericity conditions \eqref{eq:An:genericity} are certainly satisfied. Condition \eqref{eq:An:admissibility} is also satisfied since $|\zeta| = \frac{n(n-1)}{2}+1$ if $n$ is odd and $|\zeta|=\frac{n^2}{2}+1$ if $n$ is even. When $n=2$ we choose $\zeta=1$. Note that in each case our choice for $\zeta$ satisfies the additional constraint $\gcd{(\zeta_1,\dots, \zeta_{n-1})}=1$.

We now consider the $8$-manifold $M=\pi^\ast \Lambda^-T^\ast Q$, where $\pi\co S\ra Q$ is the orbibundle map. By abuse of notation, think of $\hat{\mu}$ as a map $\hat{\mu}\co \Sph^{4n+3}\ra\Lie{k}^\ast\otimes\Imag\HH$. Then $\hat{\mu}^{-1}(0)/K$ is a $3$-Sasaki orbifold. Despite $\hat{\mu}^{-1}(0)/K$ might be singular, $P=\hat{\mu}^{-1}(0)/K_{\rho}$ is a smooth $8$-manifold. Indeed, note that the action of $K_\rho$ on $\HH^{n+1}$ coincides with the restriction of the action of $K$ on $\HH\oplus\HH^n$ trivial on the first factor. Since $K$ acts freely on $\mu^{-1}(\overline{u}_0\triple{\zeta}u_0)$ when $u_0\neq 0$ by the genericity assumption on $\triple{\zeta}$ and on $\mu^{-1}(0)\setminus\{ 0 \}$ since $\mu^{-1}(0)/K=\C^2/\Gamma$, we conclude that $K_\rho$ acts freely on $P$. Note that $P$ is a principal $\sunitary{2}$--bundle over $S$; the fact that the structure group is certainly $\sunitary{2}$ rather than $\sorth{3}$ is because $K_\rho$ does not contain $-1\in\Sp{n+1}$ (since $\gcd{\left( |\zeta|,n \right)}=1$). Then the $8$-manifold $M$ is the total space of the associated adjoint bundle $P\times_{\sunitary{2}}\Lie{su}_2\ra S$.

\begin{prop}\label{prop:ALE:Topology}
Let $\Gamma=\Z_n$ and assume that $\triple{\zeta}=2\pi i\zeta\in 2\pi i\Z^{n-1}$ satisfies \eqref{eq:An:genericity} and \eqref{eq:An:admissibility} together with the additional constraint $\gcd{(\zeta_1,\dots, \zeta_r)}=1$. Then $S$ is simply connected, spin and $b_2 (S)=n-2$. It follows that $S$ is diffeomorphic to $\sharp_{n-2}(S^2\times S^3)$.
\proof
As in the proof of Proposition \ref{prop:ALE:Spin7:admissible}, we write $S=S_0\cup S_\infty$, where $S_0$ is a principal circle bundle over the $4$-manifold $X_{\triple{\zeta}}$ and $S_\infty = (\C^2\times S^1)/\Z_n$. The intersection $S_0\cap S_\infty = (S^3\times S^1)/\Z_n$.

Consider first the fundamental group of $S$. We have $\pi_1 (S_\infty)= \pi_1 (S_0\cap S_\infty)\simeq \Z$. Moreover, if $\gcd{(\zeta_1,\dots, \zeta_r)}=1$ then the first Chern class of $S_0\ra X_{\triple{\zeta}}$ is primitive in $H^2(X_{\triple{\zeta}},\Z)$ and therefore $S_0$ is simply connected. Van Kampen's Theorem then implies that $S$ is simply connected.

Taking into account the isomorphism $H^1(S_0)\oplus H^1(S_\infty)\ra H^1(S_0\cap S_\infty)$ and the fact that $H^2(S_\infty)=H^2(S_0\cap S_\infty)=0$, the Mayer--Vietoris sequence for $S=S_0\cup S_\infty$ yields an isomorphism $H^2(S)\simeq H^2(S_0)$. Since $c_1 (S_0)\neq 0$, the Gysin sequence for the circle fibration $S_0\ra X_{\triple{\zeta}}$ immediately yields $H^2(S_0)\simeq \R^{n-2}$.

In order to prove that $S$ is spin, note that $w_2(S)\in H^2(S,\Z_2)$ must be the image of $w^{orb}_2(Q)\in H^2_{orb}(Q,\Z_2)$ since $\pi\co S\ra Q$ is a principal circle orbibundle and therefore $TS\simeq \pi^\ast TQ \oplus \R$. On the other hand, in dimension $4$ we have $w_2^{orb}(Q)=\varepsilon (Q)\in H^2_{orb}(Q,\Z_2)$, where $\varepsilon (Q)$ is the Marchiafava--Romani class of $Q$, the obstruction to lift the structure group of the standard $\R^3$--orbibundle $\hat{\mu}^{-1}(0)/K \times_{\sunitary{2}}\R^3$ over $Q$ from $\sorth{3}$ to $\sunitary{2}$; here we think of $\hat{\mu}$ as defined on $\Sph^{4n+3}$. As already noticed, since $K_\rho$ does not contain $-1\in\Sp{n+1}$ we have $\varepsilon (Q) = 0$.

The diffeomorphism-type of $S$ now follows from the Smale classification of simply connected spin $5$-manifolds \cite{Smale:spin:5}. 
\endproof 
\end{prop}

Since the $8$-manifold $M$ retracts onto $S$, applying Theorem \ref{thm:Main:Technical} to the \spins--admissible self-dual Einstein $4$-orbifolds we have constructed immediately yields the following result.

\begin{theorem}\label{thm:Spin7:ALE:An}
There exists infinitely many smooth non-compact simply connected $8$-manifolds carrying complete $\spins$--metrics.
\end{theorem}

\begin{remark}\label{rmk:HK/QK}
As an aside, we note that the self-dual Einstein $4$-orbifolds we considered in this section and their relationship with ALE gravitational instantons are instances of the so-called hyperk\"ahler/quaternionic K\"ahler correspondence \cite{Haydys,Hitchin:HK/QK}. The correspondence relates a hyperk\"ahler manifold $X$ with a circle action that fixes one complex structure and rotates the other two (a \emph{rotating} circle action) with a quaternionic K\"ahler manifold $Q$ (in general incomplete) with circle symmetry. Consider for example the triholomorphic circle action $e^{i\theta}\cdot (u_1,u_2) = (e^{i\theta}u_1,e^{i\theta}u_2)$ on $\HH^2$, the corresponding moment map $\mu (u)=\overline{u}_1 i u_1 + \overline{u}_2 i u_2$ and the Eguchi--Hanson space $X=\mu^{-1}(-i)/S^1$. The particular choice of level set of the moment map implies that $X$ admits a rotating circle action, denoted by $S^1_R$. We fix a choice of lift of $S^1_R$ to $\HH^2$ by $e^{i\theta}\cdot (u_1,u_2) = (e^{im\theta}u_1e^{-i\theta},e^{im\theta}u_2e^{-i\theta})$ for some $m\in\Z$. Haydys \cite{Haydys} and Hitchin \cite{Hitchin:HK/QK,Hitchin:triholo} then consider the principal circle bundle $\mu^{-1}(-i)$ over $X$ and construct a quaternionic K\"ahler metric on $\mu^{-1}(-i)/S^1_R$. When $m\neq 0$ we can realise this quaternionic K\"ahler metric as a quaternionic K\"ahler quotient of $\HH\PP^2$. Indeed, define an action of $S^1$ on $\HH\PP^2$ by $e^{i\theta}\cdot [u_0:u_1:u_2] = [e^{i\theta}u_0:e^{im\theta}u_1,e^{im\theta}u_2]$. The quaternionic K\"ahler quotient of $\HH\PP^1$ by this action is $Q=\{ [u_0:u_1:u_2]\in\HH\PP^2\, |\, m\,\mu(u) + \overline{u}_0 i u_0=0\}/S^1$ and therefore the open set $Q_0\subset Q$ where $u_0\neq 0$ is identified with $\mu^{-1}(-i)/S^1_R$ via the map $[u_0:u_1:u_2]\mapsto \frac{\sqrt{m}}{|u_0|^2}\left( u_1\overline{u}_0, u_2\overline{u}_0\right)$. In the limiting case $m=0$ we can identify $\mu^{-1}(-i)/S^1_R$ with (the complement of $\mu^{-1} (0)$ in) $\HH\PP^1\simeq S^4$: if $(v_1,v_2)\in \HH^2$ satisfies $\mu (v)\neq 0$ then there exists $v_0\in\HH^\ast$, unique up to a circle factor, such that $\mu (v_1\overline{v}_0,v_2\overline{v}_0)=-i$. According to \cite[Theorem 3]{Haydys}, the general hyperk\"ahler/quaternionic K\"ahler corrrespondence consists in replacing $\HH^2$ with an arbitrary hyperk\"ahler cone $\tu{C}$ admitting a triholomorphic circle action and $\HH\PP^2$ with the (singular) quaternionic K\"ahler space $(\HH\times\tu{C})/\HH^\ast$ associated with the split hyperk\"ahler cone $\HH\times\tu{C}$. Here $\HH^\ast$ acts diagonally on $\HH$ and $\tu{C}$. The hyperk\"ahler space $X$ is the hyperk\"ahler quotient of $\tu{C}$ by $S^1$ at level set $-i$, say, of the hyperk\"ahler moment map. The corresponding quaternionic K\"ahler space $Q$ can be realised as the quaternionic K\"ahler quotient of $(\HH\times\tu{C})/\HH^\ast$ by the circle action
\begin{equation}\label{eq:HK/QK:lift}
e^{i\theta}[u_0:u]=[e^{i\theta}u_0:e^{im\theta}\cdot u].
\end{equation}
Here $[u_0:u]$, where $u_0\in\HH$ and $u\in\tu{C}$, denotes a point in $(\HH\times\tu{C})/\HH^\ast$, $e^{i\psi}\cdot u$ denotes the triholomorphic circle action on $\tu{C}$ and $m\in\Z$. In particular, if $X$ is a hyperk\"ahler quotient $\mu^{-1}(\triple{\zeta})/K$ of $\HH^n$ by a subgroup $K$ of $\tu{Sp}(n)$, as in the cases we have considered, for $X$ to admit a rotating circle action we must require that $\triple{\zeta}$ exponentiates to a (non-trivial) group homomorphism $\rho\co K\ra \unitary{1}$. The cone $\tu{C}$ is then the hyperk\"ahler quotient of $\HH^n$ by $\ker{\rho}$ at the zero level-set (with induced triholomorphic action of $K/\ker{\rho}\simeq S^1$) and similarly $(\HH\times\tu{C})/\HH^\ast$ is the quaternionic K\"ahler quotient of $\HH\PP^n$ by the action of $\ker{\rho}$ induced by $K \stackrel{\rho\times\tu{id}}{\longrightarrow} \unitary{1}\times K\subset\tu{Sp}(1)\times\tu{Sp}(n)\subset\tu{Sp}(n+1)$. Choosing $m=1$ in \eqref{eq:HK/QK:lift}, the quaternionic K\"ahler space $Q$ corresponding to $X$ is the quaternionic K\"ahler quotient of $\HH\PP^n$ by $K$.
\end{remark}

\subsection{Complete $\gtwo$--manifolds from Calabi-Yau orbifolds}\label{sec:ALC:G2}
 
The framework introduced in this paper to do weighted analysis on AC orbifolds allows us to extend the results of \cite{FHN:ALC:G2:from:AC:CY3} to the construction of highly collapsed ALC $\gtwo$--holonomy metrics on principal Seifert circle bundles over AC Calabi--Yau orbifolds of complex dimension $3$. Given the language introduced in Section \ref{sec:AC:Orbifolds} and our calculation of the basic weighted $L^2$--cohomology of a Seifert circle bundle of arbitrary dimension in Theorem \ref{thm:L2:cohomology}, the main result of \cite{FHN:ALC:G2:from:AC:CY3} and its proof can be extended to the orbifold setting without any further complication.  

In contrast to the $\spins$--case, the construction of \cite{FHN:ALC:G2:from:AC:CY3}, which uses only circle bundles over smooth AC Calabi--Yau manifolds, already yielded infinitely many complete $\gtwo$--manifolds. The freedom to consider AC Calabi--Yau orbifolds is simply an addition of further examples to this already rich landscape. However, using Calabi--Yau orbifolds we will now construct infinitely many families of ALC $\gtwo$--metrics on a manifold as simple as $S^3\times\R^4$: given the wealth of examples arising in \cite{FHN:ALC:G2:from:AC:CY3} it seemed likely that many different families of ALC $\gtwo$--metrics would end up being defined on the same underlying smooth $7$-manifold, but no concrete example was given. 

\begin{theorem}\label{thm:ALC:G2}
$S^3\times\R^4$ carries infinitely many distinct families of ALC $\gtwo$--metrics.
\proof
We will describe $S^3\times\R^4$ as the total space of a Seifert circle bundle over an AC Calabi--Yau orbifold $B$ in infinitely many different ways. Forgetting for the moment the AC Calabi--Yau metric, we construct $B$ as a K\"ahler manifold with trivial canonical bundle as the K\"ahler quotient of $\C^4$ by a circle in $\sunitary{4}$. Up to conjugation we assume that the circle is embedded in $\sunitary{4}$ via $e^{i\theta}\mapsto \tu{diag}\left(e^{ip_1\theta},e^{ip_2\theta}, e^{-iq_1 \theta}, e^{-iq_2\theta}\right)$ for non-negative integers $p_1,p_2,q_1,q_2$ such that $p_1+p_2 = q_1+ q_2$. Then
\[
B=B_\zeta=\{ (z_1,z_2,z_3, z_4)\in\C^4 \textup{ such that } p_1 |z_1|^2 + p_2 |z_2|^2 - q_1 |z_3|^2 - q_2 |z_4|^2 =\zeta\}/S^1,
\]
where $\zeta\in\R$ is a parameter. If $\zeta\neq 0$ then $B$ is an orbifold; indeed, the level set
\[
M=\{ (z_1,z_2,z_3, z_4)\in\C^4 \textup{ such that } p_1 |z_1|^2 + p_2 |z_2|^2 - q_1 |z_3|^2 - q_2 |z_4|^2 =\zeta\}
\]
is a smooth manifold and $S^1$ acts on $M$ with nontrivial finite stabilisers: $M\ra M/S^1=B$ is a principal Seifert circle bundle over the orbifold $B$. 

If we further assume that $\gcd{(p_i,q_j)}=1$ for all $i,j=1,2$ (in particular, $p_i, q_j>0$) then $M$ is diffeomorphic to $S^3\times\R^4$ and $B$ is smooth outside a compact set. In order to see this, assume without loss of generality that $\zeta>0$. Then up to an anisotropic rescaling $(z_1,z_2,z_3,z_4)\mapsto \left( \sqrt{p_1}z_1,\sqrt{p_2}z_2,\sqrt{q_1}z_3,\sqrt{q_2}z_4\right)$, $M$ is cut out by the equation $|z_1|^2 + |z_2|^2 = |z_3|^2 + |z_4|^2+\zeta$ and therefore, identifying $\C^2$ and $\R^4$ with $\HH$ and $S^3$ with the unit sphere in $\HH$, can be parametrised by $(x,y)\mapsto \left( \sqrt{|y|^2 +\zeta}\, x, y\right)$ for $(x,y)\in S^3\times\R^4$. Moreover, the assumptions $\gcd{(p_i,q_j)}=1$ for all $i,j=1,2$ also imply that the only points in $M$ with non-trivial stabiliser are those with $z_1=0=z_2$ or $z_3=0=z_4$ and since $\zeta>0$ the former case is impossible. The circle action is therefore free on the complement of $S^3\times\{ 0\}\subset S^3 \times \R^4\simeq M$.

Still working under the assumption that $\zeta>0$ to fix ideas, $B$ is a $\C^2$--orbibundle over the weighted projective line $\mathbb{W}\C\PP^1 [p_1,p_2]$. In particular, $H^4_{orb}(B)=0$ so that $c_1^{orb}(M)\cup [\omega_0]=0$. Here $[\omega_0]$ is the K\"ahler class of the orbifold K\"ahler metric on $B$ induced by the K\"ahler quotient construction. In the main existence result of \cite{FHN:ALC:G2:from:AC:CY3} the topological constraint $c_1^{orb}(M)\cup [\omega_0]=0\in H^4_{orb}(B)$ plays the role of \eqref{eq:top:obstr} in Theorem \ref{thm:Main:Technical}.

It remains to show that $B$ carries an AC orbifold Calabi--Yau metric in the same K\"ahler class. First of all, note that under our assumptions the K\"ahler reduction $B_0$ of $\C^4$ at the zero-level set of the moment map has an isolated singularity at the origin and is in fact a Gorenstein toric K\"ahler cone. The existence of a Calabi--Yau cone metric on $B_0$ follows from a general result of Futaki--Ono--Wang \cite{Futaki:Ono:Wang}. The Calabi--Yau cone metric on $B_0$ is in fact explicit: the case where $p_1=p_2=p$, $q_1=p-q$ and $q_2=p+q$ where $p>q>0$ and $\gcd{(p,q)}=1$ coincides with the $Y^{p,q}$ Sasaki--Einstein manifolds of \cite{GMSW}; the general case was considered in \cite{CLPP:SE,Martelli:Sparks:SE}.

The existence of an AC Calabi--Yau metric on $B$ asymptotic to the Calabi--Yau cone metric on $B_0$ now follows from the general existence theory for AC Calabi--Yau metrics on crepant resolutions of Calabi--Yau cones, in particular \cite{Goto:Crepant} (since the orbifold K\"ahler class $[\omega_0]$ is not compactly supported). Indeed, we can regard $B$ as an orbifold partial small (therefore necessarily crepant) resolution of the cone $B_0$. Strictly speaking the existence result of \cite{Goto:Crepant} applies to a smooth manifold $B$; however, since the orbifold singularities of $B$ are contained in a compact set, the extension to the orbifold setting should pose no additional difficulty. In fact, in the special case where $p_1=p_2=p$, $q_1=p-q$ and $q_2=p+q$, Martelli--Sparks \cite[Theorems 1.1 and 1.3]{Martelli:Sparks:Partial:Resolutions} constructed an explicit AC Calabi--Yau metric (unique up to scale) on $B$ asymptotic to $B_0$.

In summary, $M=S^3\times\R^4$ is a principal Seifert circle bundle over the AC Calabi--Yau orbifold $B$ and the topological constraint $c_1^{orb}(M)\cup [\omega_0]=0\in H^4_{orb}(B)$ is automatically satisfied. The main existence result of \cite{FHN:ALC:G2:from:AC:CY3} guarantees the existence of a $1$-parameter family up to scale of highly collapsed ALC $\gtwo$--metrics on $M$. Up to obvious symmetries, families corresponding to different choices of $p_1,p_2,q_1,q_2$ cannot be isometric to each other since their (unique) tangent cones at infinity are distinct.
\endproof
\end{theorem}

\begin{remark*}
All the complete $\gtwo$--metrics constructed in the proof of the theorem are \emph{toric} in the sense of Madsen--Swann \cite{Madsen:Swann:toric}, \ie they admit a multi-Hamiltonian isometric action of a $3$-torus preserving the $\gtwo$--structure. Indeed, the AC Calabi--Yau orbifold metric on $B$ is itself toric (in the usual sense of K\"ahler and symplectic geometry), but only a $2$-dimensional sub-torus also preserves the holomorphic volume form. The $2$-torus symmetry lifts to a symmetry of the ALC $\gtwo$--metrics because of uniqueness results in the construction of \cite{FHN:ALC:G2:from:AC:CY3}. Finally, the $\gtwo$--metrics are also invariant under the circle action on the fibres of the Seifert bundle.
\end{remark*}

\begin{remark*}
Since in Theorem \ref{thm:ALC:G2} $B$ has no orbifold singularities outside a compact set, the application of \cite{FHN:ALC:G2:from:AC:CY3} does not really require the machinery of Section \ref{sec:AC:Orbifolds}. However, there are likely very many other examples with non-compact singular set.  
\end{remark*} 

\bibliographystyle{amsinitial}
\bibliography{Spin(7)_ALC}
\end{document}